\documentclass{usiinftr}
\usepackage{float}
\usepackage{amsmath}

\usepackage{algorithm}
\usepackage{algpseudocode}

\usepackage{pgfplotstable}
\usepackage{booktabs}
\usepackage{array}
\usepackage{colortbl}

\newcolumntype{C}{>{\arraybackslash}m{5em}}

\usepackage{color}
\usepackage{tikz}
\usepackage{pgfplots}
\pgfplotsset{compat=1.8}
\usetikzlibrary{pgfplots.groupplots} 
\usetikzlibrary{patterns}
\usepackage{subcaption}

\definecolor{mycolor0}{HTML}{2B83BA}
\definecolor{mycolor1}{RGB}{130,220,202}
\definecolor{mycolor2}{RGB}{79,122,142}  
\definecolor{mycolor3}{RGB}{170,35,3} 
\definecolor{mycolor4}{RGB}{207,170,114}
\definecolor{mycolor5}{RGB}{80,135,63}
\definecolor{mycolor6}{RGB}{255,140,190}
\definecolor{mycolor7}{RGB}{49,163,84}

\definecolor{tablecolor1}{HTML}{F0F0F0}
\definecolor{tablecolor2}{HTML}{FFFFFF}

\newcommand{\matpower}[0]{{\sc Matpower}}
\newcommand{\matlab}[0]{{\sc Matlab}}

\newcommand{\mips}[0]{{\small MIPS}}
\newcommand{\BELTISTOS}{\BELTISTOSSUITE{}}
\newcommand{\MATPOWER}{{\matpower{}}}
\newcommand{\PARDISO}{{\small PARDISO}}
\newcommand{\IPOPT}{{\small IPOPT}}
\newcommand{\FMINCON}{{\small FMINCON}}
\newcommand{\MIPS}{{\mips{}}}
\newcommand{\KNITRO}{{\small KNITRO}}

\newcommand{\BELTISTOSOPF}{{\small BELTISTOS-OPF}}
\newcommand{\BELTISTOSMP}{{\small BELTISTOS-MP}}
\newcommand{\BELTISTOSMEM}{{\small BELTISTOS-MEM}}
\newcommand{\BELTISTOSSUITE}{{\small BELTISTOS}}

\makeatletter
\newcommand*{\transpose}{%
	  {\mathpalette\@transpose{}}%
}
\newcommand*{\@transpose}[2]{%
	        \raisebox{\depth}{$\m@th#1\intercal$}%
}
\makeatother

\def\0{{ 0}}

\def\Ng{{N_\text{G}}}
\def\Ns{{N_\text{S}}}

\newcommand \discr[1]{{\boldsymbol{#1}}}   
\def\B{{\discr{B}}}
\def\Bs{{\mathbf{B}^\text{S}}}
\def\E{{\mathbf {E}}}
\def\Ezero{{\mathbf E}^{\text{0}}}

\newcommand{\Set}[1]{\mathbb{#1}}
\newcommand{\R}[1]{\Set{R}^{#1}}

\newcommand{\eref}[1]{(\ref{#1})}
\newcommand \Min[1]{{#1}^{\text{min}}}
\newcommand \Max[1]{{#1}^{\text{max}}}

\def\minim{\mathop{\hbox{\rm minimize}}}
\def\subject{\text{\rm subject to}}
\def\minimize#1{{\displaystyle\minim_{#1}}}

\newcommand{\pmat}[1]{\begin{bmatrix}#1\end{bmatrix}}

\def\bepsilon{{\boldsymbol \epsilon}}
\def\bepsilonmax{{\boldsymbol \epsilon}^\text{max}}
\def\bepsilonmin{{\boldsymbol \epsilon}^\text{min}}
\def\bepsilonmaxS{{\boldsymbol \epsilon}^\text{max}_\text{S}}
\def\bepsilonminS{{\boldsymbol \epsilon}^\text{min}_\text{S}}
\def\bepsilonmaxSi{{\boldsymbol \epsilon}^\text{max}_{\text{S},i}}
\def\bepsilonminSi{{\boldsymbol \epsilon}^\text{min}_{\text{S},i}}

\begin{document}


\title{\bf {C}omplete results for a numerical evaluation of interior point solvers for large-scale optimal power flow problems}

\author{Juraj Kardo\v s}{1}
\author{Drosos Kourounis}{1}
\author{Olaf Schenk}{1}
\author{Ray Zimmerman}{2}

\affiliation{1}{Institute of Computational Science, Universit\` a della Svizzera italiana, Switzerland}
\affiliation{2}{Charles  H.  Dyson  School  of  Applied  Economics and Management, Cornell University, NY, USA}

%
%
\TRnumber{2020-7}

%
%

\maketitle

\begin{abstract}
Recent advances in open source interior-point optimization methods and power system related software have provided researchers and educators with the necessary platform for simulating and optimizing power networks with unprecedented convenience.  Within the \matpower{} software platform a combination of several different interior point optimization packages are supported and four different AC optimal power flow (OPF) formulations are available. 
This study investigates robustness and performance of interior-point methods for different OPF formulations for minimizing the generation cost starting from different initial guesses, for a wide range of networks provided in the \matpower{} library ranging from 1951 buses to 193,000 buses. Performance profiles are presented for iteration counts, overall time, and memory consumption, revealing the most reliable optimization method for the particular metric.
Additionally,  novel structure-exploiting methods are investigated for the multi-period OPF problems, where energy storage devices couple the individual OPF problems defined at each subdivision of the time  period of interest, resulting in computationally intractable problems for standard general-purpose optimization software.
\end{abstract}
\section{Introduction}


\matpower{}~\cite{Zimmerman:2011, matpowerManual, Murillo:2013, Murillo-Sanchez2013} a package of free, open-source Matlab-language M-files has been available for power-system researchers and educators as a simulation tool for solving power flow (PF), and extensible optimal power flow (OPF) problems. It is packaged with a library of several power networks of increasing complexity. Interfaces to multiple, high-performance nonlinear optimizers such as \FMINCON{}, \IPOPT{}, \KNITRO{}, and its in-house default solver, \mips{}, are also available for its users. Recently, several different formulations of the standard AC-OPF problem were added, including polar and Cartesian representations of complex voltage variables and both current and power versions of the nodal mismatch equations.

We use performance benchmarking profiles to evaluate various optimization methods and software for power grid applications. 
We will focus on benchmarking metrics such as the overall solution runtime, memory requirements, or iteration count with a particular emphasis on power-grid application within the \matpower{} software framework. In pursuing these objectives, we focus on single-objective optimization algorithms that run in serial (i.e., that do not use parallel processing).
The reason motivating optimization benchmarking in  \matpower{} is twofold: to demonstrate the value of a novel algorithm and formulation versus more classical methods, and to evaluate the performance of an optimization algorithm and the related optimization software on networks of increasing complexity and sizes.
 Our key contribution is a detailed performance profile study  of the effects of different optimizers for large-scale single-period and multi-period optimal power flow problems that will assist users in making an informed decision about how and which software should be utilized. 
\section{OPF problems}\label{sec:OPFproblems}
The OPF problem is defined in terms of the conventional economic dispatch
problem, aiming at determining the optimal settings for optimization variables. The standard formulation of the OPF problem takes the form of a general non-linear programming problem, with the following form:
\begin{subequations}
	\label{NLP}
	\begin{align}
	\label{NLP:Objective}
	\minimize{x}    \quad &f(x) \\
	\label{NLP:EqualityConstraints}
	\subject \quad & c_E(x) = 0, \\
	\label{NLP:InequalityConstraints}
	& c_I(x) \le 0,      \\
	\label{NLP:VariableBounds}
	& x_\mathrm{min} \le x \le x_\mathrm{max}.
	\end{align}
\end{subequations}
The objective function $f(x)$ consists of polynomial costs of generator injections, the equality constraints $c_E(x)$ are the nodal balance equations, the inequality constraints $c_I(x)$ are the branch flow limits, and the $x_\mathrm{min}$ and $x_\mathrm{max}$ bounds include reference bus angles, voltage magnitudes, and active and reactive generator injections. For further details, comprehensive explanation of the variable, nodal balance equations and discussion regarding the modeling aspects, please refer to the \matpower{} OPF model~\cite{matpowerManual}.

Multi-period OPF problems couple the individual single-period OPF problems over multiple time periods $n = 1,2,\ldots,N$, for example through the equations modeling the storage device energy levels over the planning horizon. 
The evolution of the vector of storage energy levels $\bepsilon_n \in\R{\Ns}$, where $\Ns$ represents the number of storage devices, follows the update equation
\begin{align}
\label{Eq:linup}
\bepsilon_{n} = \bepsilon_{n-1} + \Bs \; {P_g^s}^{n}, \quad n = 1,\ldots,N,
\end{align}
which depends on the power output of the storage devices ${P_g^s} = \pmat{{P_g^{sd}}&  {P_g^{sc}}} \in\R{2\Ns}$ including discharging and charging powers. 
The initial storage level is denoted $\bepsilon_0$ and the constant matrix $\Bs \in\R{\Ns\times2\Ns}$ contains two diagonal matrices
\begin{align}
\Bs  = -\delta t
\pmat{\eta_{\text{d},1}^{-1}&  & &\eta_{\text{c},1} & &  \\
            & \ddots&  &  &\ddots & \\
            & & \eta_{\text{d},N_s}^{-1}& & & \eta_{\text{c},N_s} }
\end{align}
with the discharging and charging efficiencies $\eta_{\text{d},i}$ and $\eta_{\text{c},i}$, $i = 1,2,\ldots,\Ns$. The vector of storage levels has to respect the storage level bounds at each time period $n$, $\bepsilonminS \leq \bepsilon_n \leq \bepsilonmaxS$.

The linear inequality constraints introduced by storage devices, extending the constraint set of the OPF problem \eqref{NLP}, involving powers from all storage devices, generators, and time periods can be written in matrix form as 
\begin{equation}\label{Eq:storages_N}
\underbrace{ \pmat{\bepsilonminS\\\bepsilonminS\\\vdots 
 \\ \bepsilon_0}}_{\bepsilonmin}
\leq
  \underbrace{ \pmat{\bepsilon_0\\ \bepsilon_0\\ \vdots\\ \bepsilon_0  }}_{\Ezero} +
\underbrace{
        \pmat{\Bs &        &        &        \\
              \Bs &     \Bs &        &        \\
          \vdots & \vdots & \ddots &        \\
              \Bs &     \Bs & \cdots & \Bs}
}_{\E}
\underbrace{\pmat{{P_g^s}_1\\ {P_g^s}_2 \\ \vdots \\ {P_g^s}_N}}_{{P_g^s}}
\leq
\underbrace{ \pmat{\bepsilonmaxS\\ \bepsilonmaxS\\ \vdots\\ \bepsilonmaxS }}_{\bepsilonmax}.
\end{equation}
where $\B = \pmat{\0_{\Ns\times\Ng} & \Bs  } \in \R{\Ns \times \Ng+2 \Ns}$, and $\0_{k\times l}$ denotes the zero matrix of size $k\times l$. 
Storage models with known self-discharge rate can also be included through a modification of $\Ezero$. Furthermore, a storage degradation cost, modeled as affine or quadratic function of the storage powers $P_g^{s,n}$ and the state of charge $\bepsilon_n$, can be incorporated through the problem's objective function.  The resulting optimal control problem reads
\begin{subequations}
  \label{MPOPF}
\begin{align}
  \label{MPOPF:Objective}
  \minimize{x}    \quad & \sum_{n=1}^N f(x_n) \\
  \label{MPOPF:EqualityConstraints}
  \subject \quad & c_E^n(x_n) = 0, \\
  \label{MPOPF:InequalityConstraints}
          & c_I^n(x_n) \le 0,      \\
  \label{MPOPF:VariableBounds}
          & x_n^\mathrm{min} \le x_n \le x_n^\mathrm{max}, \\
  \label{MPOPF:StorageLimits} 
  & \Min{\bepsilon} \leq \Ezero + \E \, P_g^s \leq \Max{\bepsilon}.
\end{align}
\end{subequations}
Each one of the control parameter vectors $x_n = \pmat{\Theta^n & V_m^n &  P_g^n & P_g^{s,n} &  Q_g^n}$ stands for the variables from the time periods $n=1,2,\ldots,N$. We assume the following ordering for the angles $\Theta = \pmat{\Theta^1 & \Theta^2 & \ldots & \Theta^N}$.  The voltages,  active powers, reactive powers and their bounds are ordered similarly.  The MPOPF problem~\eref{MPOPF} is essentially a simultaneous formulation of $N$ OPF problems with standard power flow constraints \eref{MPOPF:EqualityConstraints}--\eref{MPOPF:VariableBounds}, coupled through the storage device limits~\eref{MPOPF:StorageLimits} introduced as linear inequality constraints.  For detailed analysis of the storage model and problem structure, see discussion in~\cite{beltistosMPOPF}. 

\section{Optimization software} 
In what follows we describe several different primal-dual interior point methods used by many practitioners for OPF problems and provided in the software package \matpower{}.

\begin{description}
  \item \IPOPT{}~\cite{IPOPT2005, IPOPT2006}
  is a software package for large-scale nonlinear optimization.  It implementes a primal-dual interior-point algorithm with a filter line-search method for nonlinear programming, second-order corrections, and inertia correction. Benchmarks have shown that MA57~\cite{HSL} and PARDISO~\cite{Pardiso:2006} are often the  most reliable direct sparse linear solvers within \IPOPT{}. It is written in C++ by Andreas W\"achter and Carl Laird and it is released as open source code\footnote{\url{https://github.com/coin-or/Ipopt}} under the Eclipse Public License (EPL). 

\item \MIPS{}~\cite{Wang:2007,Wang:Thesis, mipsManual}
is a primal-dual interior-point  solver introduced by Wang for OPF problems.  It is entirely implemeted in MATLAB code and distributed with \matpower{}. We assume that step control is enabled (not enabled by default), which implements additional step-size control in the MIPS algorithm.

\item \FMINCON{}~\cite{Byrd:1999,Byrd:2000}
is a gradient-based method, the default optimization method of the MATLAB optimization toolbox, and it is designed to work on problems where the objective and constraint functions are both continuous and have continuous first derivatives. In its default setting it uses an interior point solver that can exploit the Hessian of the Lagrangian. The IP method applies projected conjugate gradient (PCG) method to solve the linear systems in an iterative fashion.

\item \KNITRO{}~\cite{KNITRO}
is a commercial software package from Artelys\footnote{\url{https://www.artelys.com/solvers/knitro}} for solving large-scale mathematical optimization problems. \KNITRO{} is specialized for nonlinear optimization.  \KNITRO{} offers four different optimization algorithms for solving optimization problems. Two algorithms are of the interior point type, and two are of the active set type. \KNITRO{} provides both types of algorithm for greater flexibility in solving problems, and allows crossover during the solution process from one algorithm to another. The IP algorithm uses direct step computation, however iterative conjugate gradient (CG) solver may be used automatically if the direct step is suspected to be of poor quality.

\item \BELTISTOS{}~\cite{beltistosMPOPF,beltistos,patent,thesisJuraj}
is a suite of high-performance OPF algorithms\footnote{\url{https://www.beltistos.com}} including extremely scalable and low memory MPOPF and security-constrained OPF solvers. \BELTISTOSSUITE{} adopts selected algorithms implemented in \IPOPT{}, adjusted specifically for the nature of the OPF problems. Sufficiently accurate search directions are maintained by controlling pivoting and scaling schemes provided by the direct sparse solver \PARDISO{}.  Furthermore, special factorization schemes are applied for accelerating structured problems. \BELTISTOSMP{} implements structure exploiting and data compression algorithms designed for the particular structure of the MPOPF problems.
\BELTISTOSMEM{} is a variant of \BELTISTOSMP{} algorithm which performs redundant computation in favor of reducing memory footprint. 
\end{description}



\pgfplotstableread{
	{Optimizer}    {IP method}  {Version}  {SQP method}  {Gradients}  {Hessian} {License}   {Reference}         {Solver}         {Web}
	{\BELTISTOSSUITE{}} {yes}       {1.0}       {no}       {yes}        {yes}    {Free academic use}  {\cite{beltistosMPOPF,beltistos,patent,thesisJuraj}}  {\PARDISO{}}       {\href{http://www.beltistos.com}{beltistos.com}}
	{\KNITRO{}}        {yes}       {12.2.0}     {yes}      {yes}        {yes}    {Artelys}            {\cite{KNITRO}}  {MA57, CG}       {\href{https://www.artelys.com/solvers/knitro/}{https://www.artelys.com/solvers/knitro/}}
	{\IPOPT{}}         {yes}       {3.12.10}    {no}       {yes}        {yes}    {Open source (EPL)}  {\cite{IPOPT2005,IPOPT2006}}  {MA57}       {\href{https://github.com/coin-or/Ipopt}{https://github.com/coin-or/Ipopt}}
	{\MIPS{}}          {yes}       {1.3.1}      {no}       {yes}        {yes}    {Open source (BSD)}  {\cite{Wang:2007,mipsManual}}  {mldivide, `\textbackslash' (MA57)}       {\href{https://matpower.org/docs/MIPS-manual-1.3.1.pdf}{https://matpower.org/docs/MIPS-manual-1.3.1.pdf}}
	{\FMINCON{}}       {yes}       {2018b}      {yes}      {yes}        {yes}    {MATLAB}             {\cite{Byrd:2000}}  {PCG}       {\href{https://www.mathworks.com/help/optim/ug/fmincon.html}{https://www.mathworks.com/help/optim/ug/fmincon.html}}
}\optimizersTable

\begin{table}[t]
	\centering
	\caption{Open source and commercial (highlighted) optimizers. \label{tab:optimizer}}
	\pgfplotstabletypeset[
	columns={Optimizer, Version, Solver, License, Reference},
	every head row/.style={ before row=\toprule,after row=\midrule},
	columns/Optimizer/.style={string type,column type=c},
	columns/Version/.style={string type,column type=c},
	columns/IP method/.style={string type,column name={IP}, column type=c},
	columns/SQP method/.style={string type,column name={SQP}, column type=c},
	columns/Gradients/.style={string type,column type=c, column name={Gradient}},
	columns/Hessian/.style={string type,column type=c},
	columns/License/.style={string type,column type=c},
	columns/Solver/.style={string type,column type=c, column name={Linear solver}},
	columns/Web/.style={string type,column type=c},
	columns/Reference/.style={string type,column type=c},
	every last row/.style={ after row=\bottomrule},
	every row no 0/.style={ before row={\rowcolor{tablecolor2}}},
	every row no 1/.style={ before row={\rowcolor{tablecolor1}}},
	every row no 2/.style={ before row={\rowcolor{tablecolor2}}},
	every row no 3/.style={ before row={\rowcolor{tablecolor2}}},
	every row no 4/.style={ before row={\rowcolor{tablecolor1}}},
	]\optimizersTable
\end{table}

The solution of sparse linear systems is the cornerstone of an robust high-performance optimization package. Here we describe some sparse direct linear solvers. Specific results for the OPF problems are presented in section~\ref{sec:benchmarking}.
\begin{description}
  \item HSL 2002~\cite{HSLold, HSL} is an ISO Fortran library of packages for many areas in scientific computation. It is probably best known for its codes for the direct solution of sparse linear systems, including multifrontal algorithm with approximate minimum degree ordering (MA57). \IPOPT{} provides support for a wide variety of linear solvers, including HSL linear solvers MA27, MA57 and others. \matlab{} uses the MA57 routines for real sparse symmetric matrices (operator \texttt{ldl}).
  Evaluation of the individual solvers in terms of robustness and performance is provided in~\cite{HSLbenchmarks}.
%
  
  \item \PARDISO{}~\cite{Pardiso:2006} is a thread-safe, high-performance, robust, memory efficient software for solving large sparse symmetric and unsymmetric linear systems of equations on shared-memory and distributed-memory multiprocessors. \IPOPT{} and \MIPS{} contain ready to use interfaces to the solver.
\end{description}

We note that in our study we consider \IPOPT{} with \PARDISO{} and HSL MA57 solvers and \MIPS{} with the default backslash '\textbackslash' solver and \PARDISO{}. \KNITRO{} may utilize HSL routines MA27 or MA57 in order to solve linear systems arising at every iteration of the algorithm. 
\section{Performance profiles}
\label{sec:Performance}

In order to evaluate the quality of the different optimization methods for OPF problems we will use performance profiles for compact comparison of the benchmark problems using different optimization packages. 
In recent years, these performance profiles have become a very popular and widely used tool for benchmarking and evaluating the performance of several optimizers when run on a large test set. Performance profiles have been introduced in~\cite{dolan:2002} in 2002 and have rapidly become a standard in benchmarking of optimization algorithms. Comparative studies using performance profiles have been performed throughout the optimization literature~\cite{Maurer2001}, and in the evaluation of sparse linear solvers~\cite{scottgould:2004} but also pointing out some limitations~\cite{Gould:2016:NPP:2988256.2950048}. 

The profiles are generated by running the set of optimizers $\mathcal{M}$ on a set of OPF problems $\mathcal{S}$ and recording information of interest, e.g., time to solution or memory consumption. Let us assume that a power flow optimizer $ m \in \mathcal{M}$ reports a statistic $\theta_{ms} \ge 0$ for the OPF problem $s \in \mathcal{S}$; smaller statistics $\theta_{ms}$ indicates better solution strategies. We can further define $\tilde{\theta}_{s} = \min_{m \in \mathcal{M}} \{\;\theta_{ms}\}$, which represents the best statistic for a given OPF problem $m$. Then for $\alpha \ge 1$ and each $m \in \mathcal{M}$ and $s \in \mathcal{S}$ we define
\begin{equation}
      k ( \theta_{ms}, \tilde{\theta}_{s}, \alpha) = \left \{
          \begin{array}{cc}
              1 & \; \theta_{ms} \le \alpha \; \tilde{\theta}_{s}\\
              0 & \; \theta_{ms} > \alpha \; \tilde{\theta}_{s}.
          \end{array} \right.
\end{equation}
The performance profile $ p_m({\alpha}) $ of the  power flow  optimizer $m$ is then defined by
\begin{equation}
p_m({\alpha}) = 
\frac{\sum_{s \in \mathcal{S}}{k ( \theta_{ms}, \tilde{\theta}_{s}, \alpha)}}{  | \mathcal{S} |  }.
\end{equation}
Thus, in these profiles, the values of $p_m(\alpha)$ indicate the fraction of all examples, which can be solved within $\alpha$ times, the time the best solver needed, e.g., $p_m(1)$ gives the fraction of which optimizer $m$ is the most effective package and $p^*_i := \lim_{\alpha \rightarrow \infty} p_i (\alpha) $ indicates the fraction for which the algorithm succeeded. If we are just interested in the number of wins on  $\mathcal{S}$, we need only compare the values of $p_i(1)$ for all the solvers $i \in  \mathcal{M}$, but if we are interested in optimizers with a high probability of success on the set $\mathcal{S}$, we should choose those for which $p^*_i$ is largest.  Thereby, for a selected test set, performance profiles provide a very useful and convenient means of assessing the performance of optimizers relative to the best optimizer on each example from that set~\cite{Gould:2016:NPP:2988256.2950048}. 
When commenting, e.g., on a performance profile presented in their paper, Dolan and Mor{\'e} state that it ``gives a clear indication'' of the relative performance of each optimizer~\cite{dolan:2002} and one  can determine which optimizer has the highest probability $p_i(f)$ of being within a factor $f$ of the best optimizer for $f$ in a chosen interval.
In this paper we use performance profiles to compare various aspects of problem formulation, problem setup and performance of several optimizers on sets of smooth or piecewise-smooth power flow problems. Our results provide estimates for the best configuration of the problems and identification of the optimizer with the best possible performance.

\section{Numerical results\label{sec:benchmarking}} 
We proceed with the evaluation of various aspects for the set of benchmark cases with increasing complexity, listed in Table~\ref{tab:benchmarks}. The benchmarks are split into two groups, standard benchmarks used mainly to test robustness of optimization frameworks on wide spectrum of power grid networks and large-scale benchmarks used to test the performance. In collecting the test data we imposed only two conditions: The OPF problem be of order greater than 5'000 variables and the the data must be available to other users. The first condition was imposed because our interest in this study is in medium to large-scale scale problems. The second condition was to ensure that our tests could be repeated by other users and, furthermore, it enables other software developers to test their codes on the same set of examples and thus to make comparisons with other optimizers. 
Comparing  algorithms for multiobjective optimization, or optimization algorithms that use parallel processing issues are outside of the scope of this paper since it would introduce another level of complexity to the benchmarking process and most of the users are using \matpower{} in default single core mode. The memory usage was collected using the Linux utility \texttt{time -v}, which reports maximum resident set size of the process during its lifetime.

Simulations are performed on a workstation equipped with an Intel Xeon CPU E7-4880 v2 at 2.50 GHz and 1 TB RAM using the \matpower{} version 7.0 (v20-Jun-2019). The results are presented from four different perspectives, each being a contributing factor to the complexity and behavior of the optimization procedure. These factors are (i) the initial guess provided to the optimizer, (ii) the OPF formulation, (iii) the underlying linear system solver, add finally (iv) the optimization package as a whole (iv). With this experimental design, the common factors influencing the performance of a given approach were separated from the remaining options. The other factors were fixed in order to reduced variance of the benchmarks. The fixed factors include: (i) common convergence criteria of the optimization solvers, and (ii) explicit sequential (single-core) execution. For all the optimizers we used optimality, feasibility and complementarity (where available) tolerance of $10^{-6}$ with the maximum number of iterations set to $500$ and maximum time limit of 5 hours. Where necessary, we explicitly selected the IP method with exact Hessian and direct solution method of the KKT system. The remaining options of the optimization software were using default values, or the values identified by the software itself as the optimal ones.

\begin{table}[ht!]
    \centering
    \footnotesize
    \caption{Benchmark cases statistics (number of buses, generators and lines) including the number of optimization variables and nonlinear equality and inequality constraints (considering polar-based OPF formulation). \label{tab:benchmarks}}
\begin{filecontents*}{./data/CasesTable.dat}
CaseStatistics        nb  ng nl  nlc  ndc nvar nnle   nnli   nlin ni nib pq PLQ
1951rte           1951   391  2596  2099   0   4634   3902   4198      0   1   0 p   L 
2383wp            2383   327  2896  2896   0   5420   4766   5792      0   1   0 p   L 
2868rte           2868   599  3808  2281   0   6858   5736   4562      0   1   0 p   L 
ACTIVSg2000     2000   544  3206  3206   0   4864   4000   6412      0   1   0 p   Q
2869pegase        2869   510  4582  2743   0   6758   5738   5486      0   1   0 p   L 
2737sop           2737   399  3506  3506   0   5912   5474   6538      0   1   0 p   L 
2736sp            2736   420  3504  3504   0   6012   5472   6538      0   1   0 p   L 
2746wop           2746   514  3514  3514   0   6354   5492   6614      0   1   0 p   L 
2746wp            2746   520  3514  3514   0   6404   5492   6558      0   1   0 p   L 
3012wp            3012   502  3572  3572   0   6794   6024   7144      0   1   0 p   L 
3120sp            3120   505  3693  3693   0   6836   6240   7386      0   1   0 p   L 
3375wp            3374   596  4161  4161   0   7706   6748   8322      0   1   0 p   L 
6468rte           6468  1295  9000  2313   0  13734  12936   4626      0   1   0 p   L 
6470rte           6470  1330  9005  3110   0  14462  12940   6220      0   1   0 p   L 
6495rte           6495  1372  9019  3109   0  14350  12990   6218      0   1   0 p   L 
6515rte           6515  1388  9037  3131   0  14398  13030   6262      0   1   0 p   L 
9241pegase        9241  1445 16049  6295   0  21372  18482  12590      0   1   0 p   L 
13659pegase      13659  4092 20467     0   0  35502  27318      0      0   1   0 p   L 
ACTIVSg10k     10000  2485 12706 10244   0  23874  20000  20488      5   1   0 p   Q
ACTIVSg25k     25000  4834 32230 23331   0  57558  50000  46660      0   1   0 p   Q
ACTIVSg70k     70000 10390 88207 68617   0 156214 140000 137234      0   1   0 p   Q
case21k              21084  2692 25001 25001   0  54091  42168  50002  19614   7   0 p   L
case42k              42168  5384 50001 50001   0 107027  84336 100002  38458  14   0 p   L
case99k              99396 12689 117860 117860   0 250703 198792 235720  89593  33 0 p   L
case193k             192768 24611 228574 228574   0 485135 385536 457148 173432  64 0 p   L
\end{filecontents*}
\pgfplotstabletypeset[
    columns={CaseStatistics,nb,ng,nl,nvar,nnle,nnli},
    every head row/.style={ before row=\toprule,after row=\midrule},
    precision=0, fixed, fixed zerofill, column type={r},
    columns/CaseStatistics/.style={column name=Benchmark,string type,column type=l},
    columns/nb/.style={column name=$n_b$},
    columns/nvar/.style={column name=$|x|$},
    columns/ng/.style={column name=$n_g$},
    columns/nl/.style={column name=$n_l$},
    columns/nnle/.style={column name=$|c_E(x)|$},
    columns/nnli/.style={column name=$|c_I(x)|$},
    every last row/.style={ after row=\bottomrule},
    every row no 0/.style={ before row={\rowcolor{tablecolor2}}},
    every row no 1/.style={ before row={\rowcolor{tablecolor2}}},
    every row no 2/.style={ before row={\rowcolor{tablecolor2}}},
    every row no 3/.style={ before row={\rowcolor{tablecolor2}}},
    every row no 4/.style={ before row={\rowcolor{tablecolor2}}},
    every row no 5/.style={ before row={\rowcolor{tablecolor2}}},
    every row no 6/.style={ before row={\rowcolor{tablecolor2}}},
    every row no 7/.style={ before row={\rowcolor{tablecolor2}}},
    every row no 8/.style={ before row={\rowcolor{tablecolor2}}},
    every row no 9/.style={ before row={\rowcolor{tablecolor2}}},
    every row no 10/.style={ before row={\rowcolor{tablecolor2}}},
    every row no 11/.style={ before row={\rowcolor{tablecolor2}}},
    every row no 12/.style={ before row={\rowcolor{tablecolor2}}},
    every row no 13/.style={ before row={\rowcolor{tablecolor2}}},
    every row no 14/.style={ before row={\rowcolor{tablecolor2}}},
    every row no 15/.style={ before row={\rowcolor{tablecolor2}}},
    every row no 16/.style={ before row={\rowcolor{tablecolor2}}},
    every row no 17/.style={ before row={\rowcolor{tablecolor2}}},
    every row no 18/.style={ before row={\rowcolor{tablecolor2}}},
    every row no 20/.style={ before row={\rowcolor{tablecolor1}}},
    every row no 21/.style={ before row={\rowcolor{tablecolor1}}},
    every row no 22/.style={ before row={\rowcolor{tablecolor1}}},
    every row no 23/.style={ before row={\rowcolor{tablecolor1}}},
    every row no 24/.style={ before row={\rowcolor{tablecolor1}}},
    every nth row={19}{before row=\midrule \multicolumn{6}{l}{Large-scale benchmarks}\\ \midrule \rowcolor{tablecolor1}},
]
{./data/CasesTable.dat}
\end{table}



\subsection{Benchmark cases and optimization problem properties}


Table~\ref{tab:benchmarks} lists the number of buses, generators and lines for each \matpower{} benchmark case used in this study. Additionally, we also show properties of the corresponding optimization problem such as number of variables, equality and inequality constraints. Size of the optimization problem in terms of number of nonlinear equality and inequality constraints, depends on the formulation. For the Cartesian coordinate voltage case, voltage magnitude constraints (which are simple variable bounds for the polar case) are now nonlinear inequality constraints. Presented problem sizes in Table~\ref{tab:benchmarks} consider the polar voltage representation. 
In addition to the standard \matpower{} cases, there are four larger cases, case21k -- case193k, built from the case3012wp considering the largest generator outage and line contingencies. 


\subsection{Convergence tolerance \label{sec:tolerance}}
Before proceeding with the benchmarks, we analyze the selection of the convergence tolerance. First of all, the convergence tests implemented by optimizers vary in some aspects, e.g. scaling of the residual errors or type of the norms, making the user specified tolerance not necessarily equivalent amongst the optimizers. Second, various stopping tolerances were used in previous OPF studies, ranging from $10^{-3}$ to $10^{-8}$~\cite{OPFbench-2a}. The selection of the convergence criteria has to consider multiple factors, including required optimality and feasibility errors, numerical issues associated with very tight tolerances and computation time. For very tight tolerances, the linear systems become very ill-conditionded and the effect of the round-off error is pronounced, thus influencing the numerical stability and progress of optimizer. 
Many optimizers implement "acceptable" termination criteria, such that it will terminate before the desired convergence tolerance is met (e.g. there is no improvement in the objective function or feasiblity norms over some specified number of iterations or the step becomes too small). This is useful in cases where the algorithm might not be able to achieve the desired level of accuracy. For the purpose of this section, we disabled such heuristics and allow optimizers to terminate only if the desired tolerance was reached. 
We observed that there is no significant improvement of the objective function value for tight tolerances. The absolute and relative errors of the objective function value between tolerances $10^{-4}$ and $10^{-9}$ are less than order of $10^{-3}$ and $10^{-9}$, respectively. In several cases, the objective function slightly increased for tighter tolerances, in similar orders of magnitude. The constraint violations (feasibility error) and the number of iterations for different tolerances using \BELTISTOSOPF{}, \KNITRO{} and other optimizers are shown in figures 
\ref{fig:convergenceBeltistos}--\ref{fig:convergenceFmincon}.
The figures illustrate that even for a modest tolerance $10^{-4}$ the constraint violations are usually much smaller (except \FMINCON{}) and in most cases tightening the tolerance involves only a few additional iterations (except some numerically ill-posed problems). For tight tolerances below $10^{-9}$, optimizers start to become numerically unstable and the number of iterations starts to significantly grow or optimizers terminate with an error message. Some optimizers were not able to reach the tight tolerances, e.g. \IPOPT{} and \MIPS{} failed to solve most of the benchmarks for tolerances below $10^{-9}$, while \FMINCON{} terminated with feasibility errors larger than the specified tolerance.
We thus use a modest value of the tolerance $10^{-6}$ in the following benchmarks to focus on the performance of the optimizer, and isolate issues related to the numerical errors.

\begin{figure*}[ht!]
	\centering
	\begin{subfigure}[b]{\columnwidth}
		\centering
		\begin{tikzpicture}[every mark/.append style={mark size=1.5pt}, every line/.append style={thick}]


\pgfplotstableread{
Id Tolerance case1951rte case2383wp case2736sp	case2737sop case2746wop case2746wp case2868rte case2869pegase case3012wp case3120sp case3375wp case6468rte case6470rte case6495rte case6515rte  case9241pegase caseACTIVSg2000 caseACTIVSg10k case13659pegase
1 $10^{-2}$ 25	26	15	16	19	16	28	19	37	30	33	29	34	45	56	26	17	22	53
2 $10^{-3}$ 26	27	15	17	20	16	28	21	38	31	33	29	34	47	64	26	19	22	53
3 $10^{-4}$ 26	28	15	17	22	19	31	21	41	32	35	30	38	47	64	29	19	23	92
4 $10^{-5}$ 27	28	18	20	23	19	32	21	41	32	35	31	42	48	66	29	22	24	151 
5 $10^{-6}$ 27	29	18	20	23	20	35	22	42	33	36	32	46	49	67	29	22	25	240 
6 $10^{-7}$ 28	31	19	20	26	22	37	23	45	36	44	33	50	50	67	29	23	26	293 
7 $10^{-8}$ 29	32	22	24	27	23	51	25	45	36	39	33	67	51	69	32	24	30	310 
8 $10^{-9}$ 35	33	22	24	27	24	51	28	45	37	39	35	99	58	78	35	27	30	467 
9 $10^{-10}$ NaN	34	22	24	27	25	144	29	46	38	40	117	160	168	136	37	28	32	NaN 
10 $10^{-11}$ NaN	90	200	215	174	251	314	234	61	105	89	200	199	NaN	NaN	NaN	35	NaN	NaN 
11 $10^{-12}$ NaN	319	361	420	191	NaN	NaN	NaN	271	150	172	NaN	323	NaN	NaN	NaN	127	NaN	NaN
}\datatableIters

\pgfplotstableread{
Id Tolerance case1951rte case2383wp case2736sp	case2737sop case2746wop case2746wp case2868rte case2869pegase case3012wp case3120sp case3375wp case6468rte case6470rte case6495rte case6515rte  case9241pegase caseACTIVSg2000 caseACTIVSg10k case13659pegase
1 $10^{-2}$ 	   1.4019E-05	7.8111E-08	1.2049E-07	2.0534E-07	1.6032E-07	6.0361E-08	1.4479E-04	7.8043E-07	1.0232E-06	2.2829E-07	8.1614E-08	7.9880E-06	5.3137E-06	2.5141E-05	2.4579E-04	2.1523E-05	2.6996E-05	3.6003E-07	2.6579E-05
2 $10^{-3}$ 	   2.3183E-06	6.4620E-09	1.2049E-07	7.5768E-09	1.4858E-08	6.0361E-08	1.4479E-04	5.1691E-08	9.1341E-08	2.1298E-08	8.1614E-08	7.9880E-06	5.3137E-06	1.2387E-06	2.7313E-05	2.1523E-05	1.0413E-06	3.6003E-07	2.6579E-05
3 $10^{-4}$ 	   2.3183E-06	1.2797E-09	1.2049E-07	7.5768E-09	6.2639E-11	3.0245E-10	4.5130E-05	5.1691E-08	2.5695E-11	1.1132E-09	1.5833E-10	4.7315E-07	3.7442E-07	1.2387E-06	2.7313E-05	1.0163E-08	1.0413E-06	7.4595E-08	1.4909E-05
4 $10^{-5}$ 	   1.4803E-07	1.2797E-09	3.5980E-10	9.6073E-12	1.6008E-11	3.0245E-10	2.1216E-06	5.1691E-08	2.5695E-11	1.1132E-09	1.5833E-10	5.4268E-08	3.2425E-08	4.1500E-07	1.3097E-06	1.0163E-08	1.8317E-09	1.5444E-08	1.3372E-06
5 $10^{-6}$ 	   1.4803E-07	1.1557E-10	3.5980E-10	9.6073E-12	1.6008E-11	1.4213E-11	6.0261E-08	3.0077E-07	9.3485E-12	1.2932E-10	1.1116E-11	2.6420E-09	4.2950E-09	2.1875E-08	2.6311E-08	1.0163E-08	1.8317E-09	2.0155E-09	1.9915E-07
6 $10^{-7}$ 	   6.9850E-09	6.0038E-12	7.3350E-12	7.6299E-12	7.4220E-12	1.6704E-11	5.2295E-09	8.5562E-09	1.0572E-11	9.2289E-12	2.5551E-11	4.1989E-10	4.9887E-09	9.9189E-10	2.6311E-08	1.0163E-08	8.4860E-12	1.1693E-10	1.4839E-09
7 $10^{-8}$ 	   9.8873E-12	6.1062E-12	9.4700E-12	1.3577E-11	8.6307E-12	1.5871E-11	7.1214E-10	6.5639E-12	9.5954E-12	1.2950E-11	1.6660E-11	1.6984E-10	5.7649E-10	2.5358E-10	1.1151E-11	5.5258E-12	4.7086E-12	2.2059E-11	1.6087E-09
8 $10^{-9}$ 	   1.1238E-11	7.2819E-12	7.9786E-12	1.0564E-11	2.3704E-11	1.5822E-11	1.1862E-11	6.2655E-12	9.2776E-12	9.6693E-12	1.1733E-11	2.4281E-10	1.9023E-11	1.3801E-11	8.5292E-12	5.4557E-12	6.6610E-13	4.4076E-11	6.6094E-11
9 $10^{-10}$	   NaN	6.1131E-12	5.9165E-12	1.0098E-11	1.9618E-11	1.5846E-11	1.5607E-11	4.8474E-12	9.7421E-12	1.0427E-11	1.4107E-11	1.1430E-11	1.1334E-11	9.7131E-12	1.3242E-11	7.1034E-12	1.4155E-12	3.2697E-12	5.0058E-12
10 $10^{-11}$      NaN	5.5736E-12	6.6774E-12	1.1299E-11	1.0794E-11	7.9353E-12	2.8086E-11	4.4406E-12	1.0845E-11	1.1490E-11	1.9557E-11	1.0900E-11	1.2774E-11	NaN	NaN	NaN	6.5837E-13	NaN	NaN
11 $10^{-12}$      NaN	7.3111E-12	6.5443E-12	8.9782E-12	1.5943E-11	NaN	NaN	NaN	8.0319E-12	8.7333E-12	8.5841E-12	NaN	1.2273E-11	NaN	NaN	NaN	8.3233E-13	NaN	NaN
}\datatableFeas

\begin{groupplot}[max space between ticks=15, group style = {group size = 2 by 1, horizontal sep = 25pt}, width = 8cm, height = 4.5cm]
    \nextgroupplot[ title = {Number of Iterations},
                	ymin=0,
                	ymax=160,
                	xmin=0.6,
                	xmax=11.4,
                	ymajorgrids, yminorgrids,
                    xmajorgrids, xminorgrids,
					minor xtick={2,4,...,10},
					xtick={1,3,...,11},
					xticklabels={$10^{-2}$,$10^{-4}$,$10^{-6}$,$10^{-8}$,$10^{-10}$,$10^{-12}$},
                	yticklabel style={font=\scriptsize},
                	xticklabel style={font=\scriptsize},
                	ylabel style={font=\scriptsize},
                	xlabel style={font=\scriptsize},
                	title style={font=\footnotesize,at={(0.5,0.95)}},
                	xlabel={Tolerance},
                	cycle list name=exotic,
                    legend style={legend cell align=left,align=right,draw=white!85!black,font=\tiny,legend columns=6, legend to name = grouplegend,}]
        \pgfplotstablegetcolsof{\datatableIters}
        \pgfmathparse{\pgfplotsretval-1}
          \foreach \i in {2,...,\pgfmathresult}{%
            \pgfplotstablegetcolumnnamebyindex{\i}\of{\datatableIters}\to{\colname};
            \addplot table[x index=0,y index=\i] {\datatableIters};
            \addlegendentryexpanded{\colname};
          }
    \nextgroupplot[ title = {Feasibility Error},
                	ymode=log, 
                	ymin=1e-13,
                	ymax=3.e-4,
                	xmin=0.6,
                	xmax=11.4,
                	ymajorgrids, 
                    xmajorgrids, xminorgrids,
					minor xtick={2,4,...,10},
                    xtick={1,3,...,11},
                   	xticklabels={$10^{-2}$,$10^{-4}$,$10^{-6}$,$10^{-8}$,$10^{-10}$,$10^{-12}$},
                	yticklabel style={font=\scriptsize},
                	xticklabel style={font=\scriptsize},
                	ylabel style={font=\scriptsize},
                	xlabel style={font=\scriptsize},
                	title style={font=\footnotesize,at={(0.5,0.93)}},
                	xlabel={Tolerance},
                	cycle list name=exotic,]
        \pgfplotstablegetcolsof{\datatableFeas}
        \pgfmathparse{\pgfplotsretval-1}
          \foreach \i in {2,...,\pgfmathresult}{%
            \addplot table[x index=0,y index=\i] {\datatableFeas};
          }
\end{groupplot}
\node at ($(group c2r1) + (-3.5,2.7cm)$) {\ref*{grouplegend}};

\end{tikzpicture}
	\end{subfigure}
	\vspace{-0.8cm}
	\caption{Statistics for different convergence tolerances for the standard benchmark cases --  \BELTISTOSOPF{}. \label{fig:convergenceBeltistos}}
	\begin{subfigure}[b]{\columnwidth}
		\centering
		\begin{tikzpicture}[every mark/.append style={mark size=1.5pt}, every line/.append style={thick}]

\pgfplotstableread{
	Id Tolerance case1951rte case2383wp case2736sp	case2737sop case2746wop case2746wp case2868rte case2869pegase case3012wp case3120sp case3375wp case6468rte case6470rte case6495rte case6515rte  case9241pegase caseACTIVSg2000 caseACTIVSg10k case13659pegase
	1 $10^{-2}$ 21	29	17	18	19	19	27	22	23	25	23	27	29	37	31	64	17	24	132
	2 $10^{-3}$ 21	29	17	18	19	19	27	22	23	25	23	27	29	37	31	64	17	24	132
	3 $10^{-4}$ 21	29	17	18	19	19	27	22	23	25	23	27	29	37	31	64	17	24	132
	4 $10^{-5}$ 21	30	17	18	19	19	42	22	23	25	23	27	36	37	31	64	19	25	132
	5 $10^{-6}$ 21	30	18	18	19	20	42	22	24	25	23	27	46	37	31	64	19	25	133
	6 $10^{-7}$ 22	30	18	18	20	20	42	24	26	26	23	30	67	37	37	64	22	25	133
	7 $10^{-8}$ 22	32	18	19	20	20	52	24	26	27	23	39	121	39	37	66	23	25	134
	8 $10^{-9}$ 22	34	20	20	24	20	93	147	33	31	24	33	NaN	53	66	212	23	29	136
	9 $10^{-10}$ 136	278	178	NaN	215	NaN	103	355	335	NaN	NaN	402	NaN	381	479	454	25	262	NaN
	10 $10^{-11}$ NaN	NaN	239	NaN	NaN	374	NaN	222	288	493	181	NaN	NaN	NaN	NaN	378	48	NaN	NaN
	11 $10^{-12}$ NaN	NaN	NaN	NaN	NaN	NaN	NaN	NaN	NaN	NaN	NaN	NaN	NaN	NaN	NaN	NaN	46	NaN	NaN
}\datatableIters

\pgfplotstableread{
	Id Tolerance case1951rte case2383wp case2736sp	case2737sop case2746wop case2746wp case2868rte case2869pegase case3012wp case3120sp case3375wp case6468rte case6470rte case6495rte case6515rte  case9241pegase caseACTIVSg2000 caseACTIVSg10k case13659pegase
	1 $10^{-2}$ 	1.3800E-07	8.9900E-08	3.9900E-08	7.0200E-09	1.3600E-08	1.9900E-09	3.0900E-07	2.6400E-08	8.6700E-09	5.6200E-08	2.0200E-09	4.2100E-07	1.6800E-05	2.7200E-08	6.6900E-07	4.4300E-09	8.3800E-07	2.8800E-07	3.1400E-07
	2 $10^{-3}$ 	1.3800E-07	8.9900E-08	3.9900E-08	7.0200E-09	1.3600E-08	1.9900E-09	3.0900E-07	2.6400E-08	8.6700E-09	5.6200E-08	2.0200E-09	4.2100E-07	1.6800E-05	2.7200E-08	6.6900E-07	4.4300E-09	8.3800E-07	2.8800E-07	3.1400E-07
	3 $10^{-4}$ 	1.3800E-07	8.9900E-08	3.9900E-08	7.0200E-09	1.3600E-08	1.9900E-09	3.0900E-07	2.6400E-08	8.6700E-09	5.6200E-08	2.0200E-09	4.2100E-07	1.6800E-05	2.7200E-08	6.6900E-07	4.4300E-09	8.3800E-07	2.8800E-07	3.1400E-07
	4 $10^{-5}$ 	1.3800E-07	9.2300E-09	3.9900E-08	7.0200E-09	1.3600E-08	1.9900E-09	3.9200E-08	2.6400E-08	8.6700E-09	5.6200E-08	2.0200E-09	4.2100E-07	8.0100E-08	2.7200E-08	6.6900E-07	4.4300E-09	9.9400E-09	8.5200E-09	3.1400E-07
	5 $10^{-6}$ 	1.3800E-07	9.2300E-09	1.3300E-09	7.0200E-09	1.3600E-08	1.4300E-11	3.9200E-08	2.6400E-08	5.250E-10	5.6200E-08	2.0200E-09	4.2100E-07	2.5900E-09	2.7200E-08	6.6900E-07	4.4300E-09	9.9400E-09	8.5200E-09	3.47E-08
	6 $10^{-7}$ 	5.8300E-10	9.2300E-09	1.3300E-09	7.0200E-09	1.3600E-09	1.4300E-11	3.9200E-08	6.0300E-10	1.0000E-10	2.2700E-09	2.0200E-09	1.8900E-08	2.3800E-08	2.7200E-08	8.4300E-09	4.4300E-09	1.1600E-12	8.5200E-09	3.47E-08
	7 $10^{-8}$ 	5.8300E-10	5.3800E-10	1.3300E-09	9.9900E-11	1.3600E-09	1.4300E-11	7.3200E-11	6.0300E-10	1.0000E-10	2.3100E-10	2.0200E-09	2.8200E-10	4.0300E-09	2.950E-09	8.4300E-09	1.7600E-10	8.8200E-13	8.5200E-09	3.08E-09
	8 $10^{-9}$ 	5.8300E-10	1.0000E-11	1.0000E-11	9.9900E-12	1.450E-11	2.4000E-11	7.5100E-11	4.6700E-10	1.0000E-11	1.0000E-11	1.2400E-11	8.2300E-10	NaN	4.650E-11	1.350E-11	9.4400E-12	8.0100E-13	6.1900E-12	4.04E-10
	9 $10^{-10}$	9.4200E-12	5.49E-12	5.33E-12	NaN	7.94E-12	1.21E-11	7.91E-12	5.85E-12	6.07E-12	1.07E-11	7.64E-12	1.28E-11	NaN	2.96E-11	7.88E-12	7.4E-12	1.0700E-12	2.08E-11	NaN
	10 $10^{-11}$ NaN	5.6E-12	4.39E-12	NaN	NaN	9.8E-12	NaN	5.46E-12	6.53E-12	4.27E-12	8.57E-12	NaN	NaN	NaN	NaN	6.83E-12	5.7E-13	NaN	NaN
	11 $10^{-12}$ NaN	NaN	NaN	NaN	NaN	NaN	NaN	NaN	NaN	NaN	NaN	NaN	NaN	NaN	NaN	NaN	7.25E-13	NaN	NaN
}\datatableFeas

\begin{groupplot}[max space between ticks=15, group style = {group size = 2 by 1, horizontal sep = 25pt}, width = 8cm, height = 4.5cm]
    \nextgroupplot[ title = {Number of Iterations},
                	ymin=0,
                	ymax=160,
                	xmin=0.8,
                	xmax=11.4,
                	ymajorgrids, yminorgrids,
                    xmajorgrids, xminorgrids,
            		minor xtick={2,4,...,10},
                	xtick={1,3,...,11},
                	xticklabels={$10^{-2}$,$10^{-4}$,$10^{-6}$,$10^{-8}$,$10^{-10}$,$10^{-12}$},
                	yticklabel style={font=\scriptsize},
                	xticklabel style={font=\scriptsize},
                	ylabel style={font=\scriptsize},
                	xlabel style={font=\scriptsize},
                	title style={font=\footnotesize,at={(0.5,0.95)}},
                	xlabel={Tolerance},
                	cycle list name=exotic,
                    legend style={legend cell align=left,align=right,draw=white!85!black,font=\tiny,legend columns=5, legend to name = grouplegend,}]
        \pgfplotstablegetcolsof{\datatableIters}
        \pgfmathparse{\pgfplotsretval-1}
          \foreach \i in {2,...,\pgfmathresult}{%
            \pgfplotstablegetcolumnnamebyindex{\i}\of{\datatableIters}\to{\colname};
            \addplot table[x index=0,y index=\i] {\datatableIters};
            \addlegendentryexpanded{\colname};
          }
    \nextgroupplot[ title = {Feasibility Error},
                	ymode=log, 
                	ymin=1e-13,
                	ymax=1.e-4,
                	xmin=0.8,
                	xmax=11.4,
                	ymajorgrids, 
                    xmajorgrids, xminorgrids,
					minor xtick={2,4,...,10},
					xtick={1,3,...,11},
					xticklabels={$10^{-2}$,$10^{-4}$,$10^{-6}$,$10^{-8}$,$10^{-10}$,$10^{-12}$},
                	yticklabel style={font=\scriptsize},
                	xticklabel style={font=\scriptsize},
                	ylabel style={font=\scriptsize},
                	xlabel style={font=\scriptsize},
                	title style={font=\footnotesize,at={(0.5,0.93)}},
                	xlabel={Tolerance},
                	cycle list name=exotic,]
        \pgfplotstablegetcolsof{\datatableFeas}
        \pgfmathparse{\pgfplotsretval-1}
          \foreach \i in {2,...,\pgfmathresult}{%
            \addplot table[x index=0,y index=\i] {\datatableFeas};
          }
\end{groupplot}

\end{tikzpicture}
	\end{subfigure}
	\vspace{-0.8cm}
	\caption{Statistics for different convergence tolerances for the standard benchmark cases -- \KNITRO{} 12. \label{fig:convergenceKnitro}}
\end{figure*}

\begin{figure}[ht!]
	\centering
	\begin{tikzpicture}[every mark/.append style={mark size=1.5pt}, every line/.append style={thick}]


\pgfplotstableread{
Id Tolerance case1951rte case2383wp case2736sp	case2737sop case2746wop case2746wp case2868rte case2869pegase case3012wp case3120sp case3375wp case6468rte case6470rte case6495rte case6515rte  case9241pegase caseACTIVSg2000 caseACTIVSg10k case13659pegase
1 $10^{-2}$ 	70	29	20	20	22	21	48	23	60	24	42	111	40	36	43	28	22	26	117
2 $10^{-3}$ 	70	29	20	20	22	21	48	23	60	24	42	111	40	36	43	28	22	26	117
3 $10^{-4}$ 	70	29	20	20	22	21	48	23	60	24	42	111	40	36	43	28	22	26	117
4 $10^{-5}$ 	70	30	20	21	23	21	48	23	60	24	42	111	71	36	43	28	22	27	214
5 $10^{-6}$ 	70	30	21	21	24	21	48	23	62	24	42	113	99	36	52	28	22	27	NaN 
6 $10^{-7}$ 	70	30	22	21	24	22	48	23	65	25	42	113	173	40	52	28	23	27	NaN 
7 $10^{-8}$ 	72	33	22	22	24	22	65	27	65	26	43	114	278	41	55	30	28	27	NaN 
8 $10^{-9}$ 	255	36	25	24	24	26	66	57	68	34	44	114	NaN 141	66	56	36	34	NaN 
9 $10^{-10}$	499	390	116	NaN 279	359	NaN 459	189	197	109	NaN NaN NaN NaN 369	34	86	NaN 
10 $10^{-11}$ NaN 192	163	NaN NaN 292	NaN 371	238	402	NaN NaN NaN NaN NaN 302	141	119	NaN 
11 $10^{-12}$ NaN NaN NaN NaN NaN NaN NaN NaN NaN NaN NaN NaN NaN NaN NaN NaN 203	NaN NaN 
}\datatableIters

\pgfplotstableread{
Id Tolerance case1951rte case2383wp case2736sp	case2737sop case2746wop case2746wp case2868rte case2869pegase case3012wp case3120sp case3375wp case6468rte case6470rte case6495rte case6515rte  case9241pegase caseACTIVSg2000 caseACTIVSg10k case13659pegase
1 $10^{-2}$ 	3.0200E-08	6.1500E-08	5.6200E-09	8.8300E-08	1.1600E-07	3.6300E-09	2.0700E-09	1.7100E-09	1.8500E-08	1.2700E-08	1.5300E-08	1.3500E-06	4.7900E-07	2.5900E-07	5.0000E-06	1.8500E-08	1.4600E-07	2.9800E-07	1.2900E-06
2 $10^{-3}$ 	3.0200E-08	6.1500E-08	5.6200E-09	8.8300E-08	1.1600E-07	3.6300E-09	2.0700E-09	1.7100E-09	1.8500E-08	1.2700E-08	1.5300E-08	1.3500E-06	4.7900E-07	2.5900E-07	5.0000E-06	1.8500E-08	1.4600E-07	2.9800E-07	1.2900E-06
3 $10^{-4}$ 	3.0200E-08	6.1500E-08	5.6200E-09	8.8300E-08	1.1600E-07	3.6300E-09	2.0700E-09	1.7100E-09	1.8500E-08	1.2700E-08	1.5300E-08	1.3500E-06	4.7900E-07	2.5900E-07	5.0000E-06	1.8500E-08	1.4600E-07	2.9800E-07	1.2900E-06
4 $10^{-5}$ 	3.0200E-08	3.4200E-09	5.6200E-09	6.2000E-09	5.2900E-09	3.6300E-09	2.0700E-09	1.7100E-09	1.8500E-08	1.2700E-08	1.5300E-08	1.3500E-06	3.1500E-08	2.5900E-07	5.0000E-06	1.8500E-08	1.4600E-07	9.2300E-09	1.9900E-09
5 $10^{-6}$ 	3.0200E-08	3.4200E-09	1.3100E-09	6.2000E-09	8.4500E-11	3.6300E-09	2.0700E-09	1.7100E-09	6.0600E-10	1.2700E-08	1.5300E-08	6.3500E-08	5.1100E-10	2.5900E-07	6.2700E-08	1.8500E-08	1.4600E-07	9.2300E-09	NaN
6 $10^{-7}$ 	3.0200E-08	3.4200E-09	1.1200E-10	6.2000E-09	8.4500E-11	3.6900E-11	2.0700E-09	1.7100E-09	9.3500E-12	1.4200E-09	1.5300E-08	6.3500E-08	1.1000E-11	1.7100E-08	6.2700E-08	1.8500E-08	2.0000E-08	9.2300E-09	NaN
7 $10^{-8}$ 	7.5500E-09	4.6100E-10	1.1200E-10	2.2800E-10	8.4500E-11	3.6900E-11	3.0200E-09	1.6800E-11	9.3500E-12	7.4200E-10	6.1700E-10	9.5200E-09	1.0400E-11	8.5400E-09	6.9000E-10	1.3200E-09	1.2500E-12	9.2300E-09	NaN
8 $10^{-9}$ 	1.2200E-10	5.8300E-11	1.5000E-11	5.7000E-11	8.3900E-11	1.2600E-11	1.5200E-10	2.0100E-10	9.5500E-12	1.1400E-11	9.3200E-12	1.3600E-11	NaN 7.0200E-10	4.0600E-10	5.0400E-10	6.4900E-13	5.7900E-12	NaN
9 $10^{-10}$	1.0700E-11	6.3100E-12	5.7900E-12	NaN 2.3800E-11	9.3200E-12	NaN 4.8800E-12	9.1700E-12	7.1400E-12	8.7100E-12	NaN NaN NaN NaN 3.8800E-12	6.2500E-13	2.5000E-11	NaN
10 $10^{-11}$ NaN 9.3100E-12	5.1000E-12	NaN NaN 7.9600E-12	NaN 7.6900E-12	8.9000E-12	6.7100E-12	NaN NaN NaN NaN NaN 6.9900E-12	8.9300E-13	6.1900E-12	NaN
11 $10^{-12}$ NaN NaN NaN NaN NaN NaN NaN NaN NaN NaN NaN NaN NaN NaN NaN NaN 8.4600E-13	NaN NaN
}\datatableFeas

\begin{groupplot}[max space between ticks=15, group style = {group size = 2 by 1, horizontal sep = 25pt}, width = 8cm, height = 4.5cm]
    \nextgroupplot[ title = {Number of Iterations},
                	ymin=0,
                	ymax=160,
                	xmin=0.6,
                	xmax=11.4,
                	ymajorgrids, yminorgrids,
                    xmajorgrids, xminorgrids,
            		minor xtick={2,4,...,10},
                	xtick={1,3,...,11},
                	xticklabels={$10^{-2}$,$10^{-4}$,$10^{-6}$,$10^{-8}$,$10^{-10}$,$10^{-12}$},
                	yticklabel style={font=\scriptsize},
                	xticklabel style={font=\scriptsize},
                	ylabel style={font=\scriptsize},
                	xlabel style={font=\scriptsize},
                	title style={font=\footnotesize,at={(0.5,0.95)}},
                	xlabel={Tolerance},
                	cycle list name=exotic,
                    legend style={legend cell align=left,align=right,draw=white!85!black,font=\tiny,legend columns=6, legend to name = grouplegend,}]
        \pgfplotstablegetcolsof{\datatableIters}
        \pgfmathparse{\pgfplotsretval-1}
          \foreach \i in {2,...,\pgfmathresult}{%
            \pgfplotstablegetcolumnnamebyindex{\i}\of{\datatableIters}\to{\colname};
            \addplot table[x index=0,y index=\i] {\datatableIters};
            \addlegendentryexpanded{\colname};
          }
    \nextgroupplot[ title = {Feasibility Error},
                	ymode=log, 
                	ymax=1.e-4,
                	xmin=0.6,
                	xmax=11.4,
                	ymajorgrids, 
                    xmajorgrids, xminorgrids,
					minor xtick={2,4,...,10},
					xtick={1,3,...,11},
					xticklabels={$10^{-2}$,$10^{-4}$,$10^{-6}$,$10^{-8}$,$10^{-10}$,$10^{-12}$},
                	yticklabel style={font=\scriptsize},
                	xticklabel style={font=\scriptsize},
                	ylabel style={font=\scriptsize},
                	xlabel style={font=\scriptsize},
                	title style={font=\footnotesize,at={(0.5,0.93)}},
                	xlabel={Tolerance},
                	cycle list name=exotic,]
        \pgfplotstablegetcolsof{\datatableFeas}
        \pgfmathparse{\pgfplotsretval-1}
          \foreach \i in {2,...,\pgfmathresult}{%
            \addplot table[x index=0,y index=\i] {\datatableFeas};
          }
\end{groupplot}

\end{tikzpicture}
	\vspace{-0.2cm}
	\caption{Statistics for different convergence tolerances for the standard benchmark cases -- KNITRO{} 11. }
\end{figure}

\begin{figure}[ht!]
	\centering
	\begin{tikzpicture}[every mark/.append style={mark size=1.5pt}, every line/.append style={thick}]


\pgfplotstableread{
Id Tolerance case1951rte case2383wp case2736sp	case2737sop case2746wop case2746wp case2868rte case2869pegase case3012wp case3120sp case3375wp case6468rte case6470rte case6495rte case6515rte  case9241pegase caseACTIVSg2000 caseACTIVSg10k case13659pegase
1 $10^{-2}$ 	27	72	21	36	42	21	35	31	27	32	36	31	47	50	51	39	23	31	66
2 $10^{-3}$ 	27	72	21	37	43	23	36	33	29	33	36	31	47	50	51	40	25	31	109
3 $10^{-4}$ 	28	75	23	39	45	23	54	33	30	34	38	31	69	50	51	40	27	34	237
4 $10^{-5}$ 	28	75	25	39	45	24	55	33	32	36	38	33	159	50	51	40	29	34	443
5 $10^{-6}$ 	29	75	26	41	45	26	56	33	34	38	39	33	250	52	54	40	29	36	NaN
6 $10^{-7}$ 	31	NaN	28	41	48	42	58	34	NaN	NaN	NaN	35	348	52	54	42	31	36	NaN
7 $10^{-8}$ 	31	NaN	NaN	NaN	NaN	NaN	60	36	NaN	NaN	NaN	39	477	54	56	45	33	NaN	NaN
8 $10^{-9}$ 	49	NaN	NaN	NaN	NaN	NaN	154	46	NaN	NaN	NaN	40	NaN	68	161	NaN	NaN	NaN	NaN
9 $10^{-10}$	NaN	NaN	NaN	NaN	NaN	NaN	NaN	NaN	NaN	NaN	NaN	NaN	NaN	NaN	NaN	NaN	NaN	NaN	NaN
10 $10^{-11}$ NaN	NaN	NaN	NaN	NaN	NaN	NaN	NaN	NaN	NaN	NaN	NaN	NaN	NaN	NaN	NaN	NaN	NaN	NaN
11 $10^{-12}$ NaN	NaN	NaN	NaN	NaN	NaN	NaN	NaN	NaN	NaN	NaN	NaN	NaN	NaN	NaN	NaN	NaN	NaN	NaN
}\datatableIters

\pgfplotstableread{
Id Tolerance case1951rte case2383wp case2736sp	case2737sop case2746wop case2746wp case2868rte case2869pegase case3012wp case3120sp case3375wp case6468rte case6470rte case6495rte case6515rte  case9241pegase caseACTIVSg2000 caseACTIVSg10k case13659pegase
1 $10^{-2}$ 	8.48000E-06	4.76000E-09	6.34000E-08	4.04000E-07	1.07000E-07	1.05000E-07	7.00000E-06	4.89000E-08	1.97000E-07	3.38000E-07	8.81000E-08	1.48000E-06	3.06000E-06	1.69000E-07	2.89000E-07	1.38000E-06	3.10000E-05	6.45000E-07	5.83000E-05
2 $10^{-3}$ 	8.48000E-06	4.76000E-09	6.34000E-08	1.98000E-08	1.40000E-08	4.55000E-09	2.54000E-05	7.27000E-09	7.84000E-09	5.42000E-08	8.81000E-08	1.48000E-06	3.06000E-06	1.69000E-07	2.89000E-07	1.30000E-07	2.70000E-06	6.45000E-07	4.67000E-06
3 $10^{-4}$ 	1.73000E-06	1.46000E-11	1.19000E-09	7.26000E-10	1.11000E-10	4.55000E-09	2.33000E-05	7.27000E-09	1.60000E-09	6.94000E-09	5.34000E-10	1.48000E-06	1.58000E-08	1.69000E-07	2.89000E-07	1.30000E-07	1.77000E-07	1.14000E-08	1.06000E-07
4 $10^{-5}$ 	1.73000E-06	1.46000E-11	5.80000E-11	7.26000E-10	1.11000E-10	5.00000E-10	7.25000E-07	7.27000E-09	1.06000E-11	1.91000E-10	5.34000E-10	2.04000E-07	8.33000E-10	1.69000E-07	2.89000E-07	1.30000E-07	1.39000E-08	1.14000E-08	1.64000E-10
5 $10^{-6}$ 	3.25000E-07	1.46000E-11	1.04000E-11	8.99000E-12	1.11000E-10	1.12000E-11	3.42000E-07	7.27000E-09	9.31000E-12	8.57000E-12	2.04000E-11	2.04000E-07	1.35000E-11	9.22000E-10	3.45000E-10	1.30000E-07	1.39000E-08	3.33000E-10	NaN
6 $10^{-7}$ 	6.47000E-10	NaN	4.84000E-12	8.99000E-12	1.12000E-11	1.58000E-11	3.66000E-09	7.77000E-10	NaN	NaN	NaN	7.10000E-08	1.42000E-11	9.22000E-10	3.45000E-10	4.03000E-09	9.16000E-10	3.33000E-10	NaN
7 $10^{-8}$ 	6.47000E-10	NaN	NaN	NaN	NaN	NaN	2.82000E-11	1.34000E-11	NaN	NaN	NaN	1.81000E-10	2.58000E-11	1.16000E-10	7.94000E-12	5.59000E-11	2.40000E-11	NaN	NaN
8 $10^{-9}$ 	1.11000E-11	NaN	NaN	NaN	NaN	NaN	2.20000E-11	6.08000E-12	NaN	NaN	NaN	1.04000E-11	NaN	2.45000E-10	8.34000E-10	NaN	NaN	NaN	NaN
9 $10^{-10}$	NaN	NaN	NaN	NaN	NaN	NaN	NaN	NaN	NaN	NaN	NaN	NaN	NaN	NaN	NaN	NaN	NaN	NaN	NaN
10 $10^{-11}$ NaN	NaN	NaN	NaN	NaN	NaN	NaN	NaN	NaN	NaN	NaN	NaN	NaN	NaN	NaN	NaN	NaN	NaN	NaN
11 $10^{-12}$ NaN	NaN	NaN	NaN	NaN	NaN	NaN	NaN	NaN	NaN	NaN	NaN	NaN	NaN	NaN	NaN	NaN	NaN	NaN
}\datatableFeas

\begin{groupplot}[max space between ticks=15, group style = {group size = 2 by 1, horizontal sep = 25pt}, width = 8cm, height = 4.5cm]
    \nextgroupplot[ title = {Number of Iterations},
                	ymin=0,
                	ymax=160,
                	xmin=0.6,
                	xmax=11.4,
                	ymajorgrids, yminorgrids,
                    xmajorgrids, xminorgrids,
                	xtick={1,2,3,4,5,6,7,8,9,10,11},
                	xticklabels={$10^{-2}$,$10^{-3}$,$10^{-4}$,$10^{-5}$,$10^{-6}$,$10^{-7}$,$10^{-8}$,$10^{-9}$,$10^{-10}$,$10^{-11}$,$10^{-12}$},
                	yticklabel style={font=\scriptsize},
                	xticklabel style={font=\scriptsize},
                	ylabel style={font=\scriptsize},
                	xlabel style={font=\scriptsize},
                	title style={font=\footnotesize,at={(0.5,0.95)}},
                	xlabel={Tolerance},
                	cycle list name=exotic,
                    legend style={legend cell align=left,align=right,draw=white!85!black,font=\tiny,legend columns=6, legend to name = grouplegend,}]
        \pgfplotstablegetcolsof{\datatableIters}
        \pgfmathparse{\pgfplotsretval-1}
          \foreach \i in {2,...,\pgfmathresult}{%
            \pgfplotstablegetcolumnnamebyindex{\i}\of{\datatableIters}\to{\colname};
            \addplot table[x index=0,y index=\i] {\datatableIters};
            \addlegendentryexpanded{\colname};
          }
    \nextgroupplot[ title = {Feasibility Error},
                	ymode=log, 
                	ymin=1e-13,
                	ymax=1.e-4,
                	xmin=0.6,
                	xmax=11.4,
                	ymajorgrids, 
                    xmajorgrids, xminorgrids,
                	xtick={1,2,3,4,5,6,7,8,9,10,11},
                	xticklabels={$10^{-2}$,$10^{-3}$,$10^{-4}$,$10^{-5}$,$10^{-6}$,$10^{-7}$,$10^{-8}$,$10^{-9}$,$10^{-10}$,$10^{-11}$,$10^{-12}$},
                	yticklabel style={font=\scriptsize},
                	xticklabel style={font=\scriptsize},
                	ylabel style={font=\scriptsize},
                	xlabel style={font=\scriptsize},
                	title style={font=\footnotesize,at={(0.5,0.93)}},
                	xlabel={Tolerance},
                	cycle list name=exotic,]
        \pgfplotstablegetcolsof{\datatableFeas}
        \pgfmathparse{\pgfplotsretval-1}
          \foreach \i in {2,...,\pgfmathresult}{%
            \addplot table[x index=0,y index=\i] {\datatableFeas};
          }
\end{groupplot}

\end{tikzpicture}
	\vspace{-0.2cm}
	\caption{Statistics for different convergence tolerances for the standard benchmark cases -- \IPOPT{}. \label{fig:convergenceIPOPT}}
\end{figure}

\begin{figure}[ht!]
	\centering
	\begin{tikzpicture}[every mark/.append style={mark size=1.5pt}, every line/.append style={thick}]


\pgfplotstableread{
Id Tolerance case1951rte case2383wp case2736sp	case2737sop case2746wop case2746wp case2868rte case2869pegase case3012wp case3120sp case3375wp case6468rte case6470rte case6495rte case6515rte  case9241pegase caseACTIVSg2000 caseACTIVSg10k case13659pegase
1 $10^{-2}$ 	14	25	21	20	20	20	11	20	26	31	26	21	24	49	34	29	38	73	23
2 $10^{-3}$ 	21	29	25	23	27	26	14	27	28	32	29	27	31	64	46	39	42	77	46
3 $10^{-4}$ 	26	31	27	25	28	28	26	29	29	33	30	39	44	67	51	41	43	80	50
4 $10^{-5}$ 	28	32	27	25	28	28	29	31	30	33	31	41	47	70	55	44	45	81	53
5 $10^{-6}$ 	28	32	27	25	29	28	31	32	30	33	31	43	47	72	55	46	47	81	54
6 $10^{-7}$ 	29	33	28	26	30	29	31	32	NaN 34	31	43	47	72	55	46	47	82	57
7 $10^{-8}$ 	29	34	29	27	31	30	31	32	NaN 35	31	43	47	72	55	46	48	82	NaN
8 $10^{-9}$ 	29	NaN 30	28	32	33	31	32	NaN NaN NaN 43	47	72	55	48	49	83	NaN
9 $10^{-10}$	30	NaN NaN 29	33	53	31	34	NaN NaN NaN NaN 48	NaN 56	48	50	84	NaN
10 $10^{-11}$ 32	NaN NaN 30	34	54	NaN NaN NaN NaN NaN NaN 49	NaN NaN NaN 51	85	NaN
11 $10^{-12}$ NaN NaN NaN NaN NaN NaN NaN NaN NaN NaN NaN NaN NaN NaN NaN NaN 52	NaN NaN
}\datatableIters

\pgfplotstableread{
Id Tolerance case1951rte case2383wp case2736sp	case2737sop case2746wop case2746wp case2868rte case2869pegase case3012wp case3120sp case3375wp case6468rte case6470rte case6495rte case6515rte  case9241pegase caseACTIVSg2000 caseACTIVSg10k case13659pegase
1 $10^{-2}$ 	5.77254E-06	3.38922E-06	1.44635E-06	2.20859E-06	2.92690E-06	3.10709E-06	1.10804E-05	2.47616E-05	3.16230E-07	1.52352E-07	5.71170E-07	3.17007E-05	3.02291E-05	1.08747E-04	2.62769E-05	1.52536E-05	4.31723E-06	1.72732E-07	1.02360E-04
2 $10^{-3}$ 	1.65794E-06	7.33075E-07	6.81588E-08	1.64565E-07	4.50359E-09	7.32797E-08	4.52371E-06	6.15401E-07	2.17674E-08	1.52753E-08	3.42567E-08	9.51904E-06	6.03711E-06	1.43247E-06	8.36781E-06	6.95725E-07	1.74620E-07	2.34485E-08	4.85903E-06
3 $10^{-4}$ 	2.16364E-07	8.05161E-08	5.46034E-11	3.17877E-10	3.54773E-10	1.89035E-11	2.16426E-07	5.15656E-08	1.34714E-09	1.16267E-09	2.05330E-09	9.81749E-08	3.46936E-08	1.38797E-07	2.13352E-07	1.13883E-07	5.44686E-08	1.99357E-10	1.64816E-05
4 $10^{-5}$ 	1.31557E-09	3.04432E-10	5.46034E-11	3.17877E-10	3.54773E-10	1.89035E-11	6.23528E-09	3.91808E-07	1.52173E-10	1.16267E-09	1.13773E-10	7.26147E-09	2.78234E-10	1.22227E-08	3.21313E-10	1.79013E-07	2.74599E-09	1.99762E-11	7.61363E-07
5 $10^{-6}$ 	1.31557E-09	3.04432E-10	5.46034E-11	3.17877E-10	3.12167E-11	1.89035E-11	8.97798E-11	8.42364E-11	1.52173E-10	1.16267E-09	1.13773E-10	3.23097E-10	2.78234E-10	4.46590E-10	3.21313E-10	1.56215E-11	1.32974E-10	1.99762E-11	4.43248E-07
6 $10^{-7}$ 	3.47411E-11	4.69852E-11	1.74189E-11	2.72766E-11	4.01116E-12	9.75612E-13	8.97798E-11	8.42364E-11	NaN 1.19564E-10	1.13773E-10	3.23097E-10	2.78234E-10	4.46590E-10	3.21313E-10	1.56215E-11	1.32974E-10	2.46267E-12	8.99834E-08
7 $10^{-8}$ 	3.47411E-11	1.12510E-11	4.22393E-12	4.67791E-12	1.43140E-13	6.03126E-14	8.97798E-11	8.42364E-11	NaN 9.69673E-12	1.13773E-10	3.23097E-10	2.78234E-10	4.46590E-10	3.21313E-10	1.56215E-11	4.70140E-11	2.46267E-12	NaN
8 $10^{-9}$ 	3.47411E-11	NaN 9.08036E-13	3.19983E-13	6.02536E-14	3.04828E-14	8.97798E-11	8.42364E-11	NaN NaN NaN 3.23097E-10	2.78234E-10	4.46590E-10	3.21313E-10	1.87032E-12	1.13280E-11	4.15685E-13	NaN
9 $10^{-10}$	1.40336E-11	NaN NaN 4.52178E-14	5.82666E-14	3.03895E-14	8.97798E-11	3.50666E-12	NaN NaN NaN NaN 1.60041E-12	NaN 4.31463E-11	1.87032E-12	2.68186E-12	6.05266E-14	NaN
10 $10^{-11}$ 8.13566E-12	NaN NaN 4.32883E-14	5.12818E-14	2.47274E-14	NaN NaN NaN NaN NaN NaN 1.30973E-12	NaN NaN NaN 5.25627E-13	5.00310E-15	NaN
11 $10^{-12}$ NaN NaN NaN NaN NaN NaN NaN NaN NaN NaN NaN NaN NaN NaN NaN NaN 9.01477E-14	NaN NaN
}\datatableFeas

\begin{groupplot}[max space between ticks=15, group style = {group size = 2 by 1, horizontal sep = 25pt}, width = 8cm, height = 4.5cm]
    \nextgroupplot[ title = {Number of Iterations},
                	ymax=160,
                	xmin=0.6,
                	xmax=11.4,
                	ymajorgrids, yminorgrids,
                    xmajorgrids, xminorgrids,
                	xtick={1,2,3,4,5,6,7,8,9,10,11},
                	xticklabels={$10^{-2}$,$10^{-3}$,$10^{-4}$,$10^{-5}$,$10^{-6}$,$10^{-7}$,$10^{-8}$,$10^{-9}$,$10^{-10}$,$10^{-11}$,$10^{-12}$},
                	yticklabel style={font=\scriptsize},
                	xticklabel style={font=\scriptsize},
                	ylabel style={font=\scriptsize},
                	xlabel style={font=\scriptsize},
                	title style={font=\footnotesize,at={(0.5,0.95)}},
                	xlabel={Tolerance},
                	cycle list name=exotic,
                    legend style={legend cell align=left,align=right,draw=white!85!black,font=\tiny,legend columns=6, legend to name = grouplegend,}]
        \pgfplotstablegetcolsof{\datatableIters}
        \pgfmathparse{\pgfplotsretval-1}
          \foreach \i in {2,...,\pgfmathresult}{%
            \pgfplotstablegetcolumnnamebyindex{\i}\of{\datatableIters}\to{\colname};
            \addplot table[x index=0,y index=\i] {\datatableIters};
            \addlegendentryexpanded{\colname};
          }
    \nextgroupplot[ title = {Feasibility Error},
                	ymode=log, 
                	ymax=2.e-4,
                	xmin=0.6,
                	xmax=11.4,
                	ymajorgrids, 
                    xmajorgrids, xminorgrids,
                	xtick={1,2,3,4,5,6,7,8,9,10,11},
                	xticklabels={$10^{-2}$,$10^{-3}$,$10^{-4}$,$10^{-5}$,$10^{-6}$,$10^{-7}$,$10^{-8}$,$10^{-9}$,$10^{-10}$,$10^{-11}$,$10^{-12}$},
                	yticklabel style={font=\scriptsize},
                	xticklabel style={font=\scriptsize},
                	ylabel style={font=\scriptsize},
                	xlabel style={font=\scriptsize},
                	title style={font=\footnotesize,at={(0.5,0.93)}},
                	xlabel={Tolerance},
                	cycle list name=exotic,]
        \pgfplotstablegetcolsof{\datatableFeas}
        \pgfmathparse{\pgfplotsretval-1}
          \foreach \i in {2,...,\pgfmathresult}{%
            \addplot table[x index=0,y index=\i] {\datatableFeas};
          }
\end{groupplot}

\end{tikzpicture}
	\vspace{-0.2cm}
	\caption{Statistics for different convergence tolerances for the standard benchmark cases -- \MIPS{}. \label{fig:convergenceMips}}
\end{figure}

\begin{figure}[ht!]
	\centering
	\begin{tikzpicture}[every mark/.append style={mark size=1.5pt}, every line/.append style={thick}]


\pgfplotstableread{
Id Tolerance case1951rte case2383wp case2736sp	case2737sop case2746wop case2746wp case2868rte case2869pegase case3012wp case3120sp case3375wp case6468rte case6470rte case6495rte case6515rte  case9241pegase caseACTIVSg2000 caseACTIVSg10k case13659pegase
1 $10^{-2}$ 	51	64	36	27	34	36	45	13	NaN 52	NaN NaN 24	47	54	31	36	40	24
2 $10^{-3}$ 	62	84	38	30	39	40	45	28	NaN 52	NaN NaN 24	47	74	44	38	52	24
3 $10^{-4}$ 	62	85	44	31	39	41	61	28	NaN 52	NaN NaN 45	47	74	44	40	52	64
4 $10^{-5}$ 	62	85	45	32	39	41	61	28	NaN 56	NaN NaN 45	58	74	44	40	52	72
5 $10^{-6}$ 	65	85	45	32	39	44	61	34	NaN 56	NaN NaN 67	58	80	44	40	52	85
6 $10^{-7}$ 	66	86	45	33	41	44	71	34	NaN 56	NaN NaN 67	58	80	49	41	57	159
7 $10^{-8}$ 	66	86	45	35	41	44	71	35	NaN 57	NaN NaN 67	62	80	50	43	57	159
8 $10^{-9}$ 	68	88	47	35	43	46	71	39	NaN 59	NaN NaN 70	63	83	51	45	57	159
9 $10^{-10}$	68	90	50	37	44	46	73	41	NaN 60	NaN NaN 70	65	84	55	47	60	159
10 $10^{-11}$ 94	90	62	38	44	50	80	60	NaN 61	NaN NaN 87	73	NaN	64	48	65	159
11 $10^{-12}$ 94	106	NaN NaN 49	58	80	NaN NaN NaN NaN NaN 87	73	NaN	64	49	65	159
}\datatableIters

\pgfplotstableread{
Id Tolerance case1951rte case2383wp case2736sp	case2737sop case2746wop case2746wp case2868rte case2869pegase case3012wp case3120sp case3375wp case6468rte case6470rte case6495rte case6515rte  case9241pegase caseACTIVSg2000 caseACTIVSg10k case13659pegase
1 $10^{-2}$ 	8.15E-05	3.84E-03	1.10E-03	3.64E-04	1.02E-02	1.18E-02	9.20E-04	5.14E-04	NaN	2.15E-05	NaN	NaN	9.07E-04	3.53E-04	2.97E-05	7.71E-04	3.29E-02	1.06E-02	5.97E-04
2 $10^{-3}$ 	1.15E-04	1.30E-07	1.59E-04	2.53E-05	3.03E-07	1.62E-04	9.20E-04	1.09E-06	NaN	2.15E-05	NaN	NaN	9.07E-04	3.53E-04	8.98E-09	4.60E-06	1.08E-03	2.00E-06	5.97E-04
3 $10^{-4}$ 	1.15E-04	2.47E-07	1.00E-06	6.39E-06	3.03E-07	3.42E-07	8.69E-09	1.09E-06	NaN	2.15E-05	NaN	NaN	1.26E-05	3.53E-04	8.98E-09	4.60E-06	3.18E-06	2.00E-06	4.83E-07
4 $10^{-5}$ 	1.15E-04	2.47E-07	2.11E-08	3.76E-07	3.03E-07	3.42E-07	8.69E-09	1.09E-06	NaN	1.24E-07	NaN	NaN	1.26E-05	1.64E-06	8.98E-09	4.60E-06	3.18E-06	2.00E-06	1.52E-04
5 $10^{-6}$ 	4.96E-06	2.47E-07	2.11E-08	3.76E-07	3.03E-07	9.69E-10	8.69E-09	7.71E-08	NaN	1.24E-07	NaN	NaN	2.19E-11	1.64E-06	2.62E-07	4.60E-06	3.18E-06	2.00E-06	5.54E-09
6 $10^{-7}$ 	5.13E-07	4.51E-09	2.11E-08	1.62E-08	4.82E-09	9.69E-10	1.59E-06	7.71E-08	NaN	1.24E-07	NaN	NaN	2.19E-11	1.64E-06	2.62E-07	1.49E-06	1.21E-06	4.47E-09	3.92E-12
7 $10^{-8}$ 	5.13E-07	4.51E-09	2.11E-08	1.10E-09	4.82E-09	9.69E-10	1.59E-06	4.87E-09	NaN	1.70E-08	NaN	NaN	2.19E-11	4.69E-09	2.62E-07	7.50E-08	9.90E-08	4.47E-09	3.92E-12
8 $10^{-9}$ 	6.30E-09	3.66E-10	3.78E-10	1.10E-09	2.42E-10	9.54E-12	1.59E-06	6.88E-10	NaN	1.36E-09	NaN	NaN	1.16E-08	9.47E-10	4.95E-10	2.14E-09	7.91E-09	4.47E-09	3.92E-12
9 $10^{-10}$	6.30E-09	7.06E-12	6.48E-12	7.08E-12	1.28E-11	9.54E-12	1.37E-09	1.20E-11	NaN	2.06E-10	NaN	NaN	1.16E-08	2.33E-10	4.60E-11	4.57E-10	5.23E-10	1.66E-11	3.92E-12
10 $10^{-11}$ 9.42E-12	7.06E-12	6.20E-12	8.83E-12	1.28E-11	7.98E-12	1.32E-11	4.15E-12	NaN	3.13E-11	NaN	NaN	9.42E-12	8.96E-12	NaN	6.47E-12	9.84E-11	1.24E-11	3.92E-12
11 $10^{-12}$ 9.42E-12	5.92E-12	NaN	NaN	5.66E-12	6.53E-12	1.32E-11	NaN	NaN	NaN	NaN	NaN	9.42E-12	8.96E-12	NaN	6.47E-12	8.29E-12	1.24E-11	3.92E-12
}\datatableFeas

\begin{groupplot}[max space between ticks=15, group style = {group size = 2 by 1, horizontal sep = 25pt}, width = 8cm, height = 4.5cm]
    \nextgroupplot[ title = {Number of Iterations},
                	ymax=160,
                	xmin=0.6,
                	xmax=11.4,
                	ymajorgrids, yminorgrids,
                    xmajorgrids, xminorgrids,
                	xtick={1,2,3,4,5,6,7,8,9,10,11},
                	xticklabels={$10^{-2}$,$10^{-3}$,$10^{-4}$,$10^{-5}$,$10^{-6}$,$10^{-7}$,$10^{-8}$,$10^{-9}$,$10^{-10}$,$10^{-11}$,$10^{-12}$},
                	yticklabel style={font=\scriptsize},
                	xticklabel style={font=\scriptsize},
                	ylabel style={font=\scriptsize},
                	xlabel style={font=\scriptsize},
                	title style={font=\footnotesize,at={(0.5,0.95)}},
                	xlabel={Tolerance},
                	cycle list name=exotic,
                    legend style={legend cell align=left,align=right,draw=white!85!black,font=\tiny,legend columns=6, legend to name = grouplegend,}]
        \pgfplotstablegetcolsof{\datatableIters}
        \pgfmathparse{\pgfplotsretval-1}
          \foreach \i in {2,...,\pgfmathresult}{%
            \pgfplotstablegetcolumnnamebyindex{\i}\of{\datatableIters}\to{\colname};
            \addplot table[x index=0,y index=\i] {\datatableIters};
            \addlegendentryexpanded{\colname};
          }
    \nextgroupplot[ title = {Feasibility Error},
                	ymode=log, 
					ymin=1e-13,
                	ymax=2.e-2,
                	xmin=0.6,
                	xmax=11.4,
                	ymajorgrids, 
                    xmajorgrids, xminorgrids,
                	xtick={1,2,3,4,5,6,7,8,9,10,11},
                	xticklabels={$10^{-2}$,$10^{-3}$,$10^{-4}$,$10^{-5}$,$10^{-6}$,$10^{-7}$,$10^{-8}$,$10^{-9}$,$10^{-10}$,$10^{-11}$,$10^{-12}$},
                	yticklabel style={font=\scriptsize},
                	xticklabel style={font=\scriptsize},
                	ylabel style={font=\scriptsize},
                	xlabel style={font=\scriptsize},
                	title style={font=\footnotesize,at={(0.5,0.93)}},
                	xlabel={Tolerance},
                	cycle list name=exotic,]
        \pgfplotstablegetcolsof{\datatableFeas}
        \pgfmathparse{\pgfplotsretval-1}
          \foreach \i in {2,...,\pgfmathresult}{%
            \addplot table[x index=0,y index=\i] {\datatableFeas};
          }
\end{groupplot}

\end{tikzpicture}
	\vspace{-0.2cm}
	\caption{Statistics for different convergence tolerances for the standard benchmark cases -- \FMINCON{}. \label{fig:convergenceFmincon}}
\end{figure}

\subsection{Initial point \label{sec:resultsInit}}

The performance of gradient-based optimization methods, in general, is sensitive to the initial point. In order to evaluate its influence, the OPF problems are solved from three different initial points currently provided by the \matpower{} option \texttt{opf.start}. The initial point for ``flat start" is heuristically chosen to be the average of the upper and lower bounds, or close to the bound if bounded only from one side. This option is the default in \matpower{} and it does not satisfy any nonlinear constraints. The option ``\matpower{} case data" uses the values of variables specified in the input \matpower{} case (usually an approximate solution, but this depends on the specific benchmark) and the third option ``power-flow solution" is the solution of the power flow equations for the given case initialized from the ``\matpower{} case data", which guarantees that the power balance constraints are satisfied. The Newton's method was used for solution of the power-flow equations, using tolerance $10^{-8}$.
The performance profiles for the number of iterations and overall time of the OPF benchmarks starting from different points are presented in Figures~\ref{fig:profileInitBeltistos} and \ref{fig:profileInitKnitro}. The default OPF formulation with polar voltage coordinates and power balance equations was considered. The results were obtained using \BELTISTOSOPF{} and \KNITRO{}, since these were the most successful optimizers, as demonstrated in Table~\ref{tab:profileGuess}.
The most successful initial points in our set of benchmark cases were the options ``\matpower{} case data" and ``power-flow solution".
In what follows we use the option ``power-flow solution" 
in order to avoid the ambiguity of the data specified in the \matpower{} case data.
The option ``\matpower{} case data" assumes that the case is well constructed and contains high-quality data, which might not always be the case.
%
The OPF problems are non-convex and different local minima may be reached from different starting points. We observed, however, that the relative differences between solutions obtained from different initial points or different optimizers were less than $10^{-5}$. 

\begin{figure*}[ht!]
	\centering
	\begin{subfigure}[b]{0.49\columnwidth}
		\begin{tikzpicture}

\pgfplotstableread{
ALPHA         SOLVER1         SOLVER2         SOLVER3
1.0         0.24000         0.56000         0.20000
1.33         0.64000         0.96000         0.84000
1.6600000000000001         0.68000         0.96000         0.92000
1.9900000000000002         0.68000         1.00000         1.00000
2.3200000000000003         0.68000         1.00000         1.00000
2.6500000000000004         0.68000         1.00000         1.00000
2.9800000000000004         0.72000         1.00000         1.00000
3.3100000000000005         0.72000         1.00000         1.00000
3.6400000000000006         0.72000         1.00000         1.00000
3.9700000000000006         0.72000         1.00000         1.00000
4.300000000000001         0.72000         1.00000         1.00000
4.630000000000001         0.76000         1.00000         1.00000
4.960000000000001         0.84000         1.00000         1.00000
5.290000000000001         0.96000         1.00000         1.00000
5.620000000000001         0.96000         1.00000         1.00000
5.950000000000001         0.96000         1.00000         1.00000
6.280000000000001         0.96000         1.00000         1.00000
6.610000000000001         0.96000         1.00000         1.00000
6.940000000000001         0.96000         1.00000         1.00000
7.270000000000001         0.96000         1.00000         1.00000
7.600000000000001         0.96000         1.00000         1.00000
7.9300000000000015         0.96000         1.00000         1.00000
8.260000000000002         0.96000         1.00000         1.00000
8.590000000000002         0.96000         1.00000         1.00000
8.920000000000002         0.96000         1.00000         1.00000
9.250000000000002         0.96000         1.00000         1.00000
9.580000000000002         0.96000         1.00000         1.00000
9.910000000000002         0.96000         1.00000         1.00000
10.240000000000002         0.96000         1.00000         1.00000
10.570000000000002         0.96000         1.00000         1.00000
10.900000000000002         0.96000         1.00000         1.00000
11.230000000000002         0.96000         1.00000         1.00000
11.560000000000002         0.96000         1.00000         1.00000
11.890000000000002         0.96000         1.00000         1.00000
12.220000000000002         0.96000         1.00000         1.00000
12.550000000000002         0.96000         1.00000         1.00000
12.880000000000003         0.96000         1.00000         1.00000
13.210000000000003         0.96000         1.00000         1.00000
13.540000000000003         0.96000         1.00000         1.00000
13.870000000000003         0.96000         1.00000         1.00000
14.200000000000003         1.00000         1.00000         1.00000
14.530000000000003         1.00000         1.00000         1.00000
14.860000000000003         1.00000         1.00000         1.00000
15.190000000000003         1.00000         1.00000         1.00000
}\datatable

  \begin{axis}[name=symb,
  	width=\columnwidth,
  	height=4cm,
	ymin=-0.02,
	ymax=1.02,
    enlarge x limits=0.05,
	axis lines*=left, 
	ymajorgrids, yminorgrids,
    xmajorgrids, xminorgrids,
	xticklabel style={rotate=0, xshift=-0.0cm, anchor=north, font=\scriptsize},
	ytick={0, 0.1, 0.2, 0.3, 0.4, 0.5, 0.6, 0.7, 0.8, 0.9, 1},
    xmax=15,
	yticklabel style={font=\scriptsize},
	ylabel style={font=\scriptsize},
	xlabel style={font=\scriptsize},
	ylabel={$p_m({\alpha})$},
	xlabel={$\alpha$},
    legend style={at={(1.0,0.55)},legend cell align=left,align=right,draw=white!85!black,font=\footnotesize,legend columns=1}
	]
	
    \addplot[color=mycolor1, ultra thick] table[x=ALPHA, y=SOLVER1] {\datatable}; 
    \addlegendentry{Flat start}

    \addplot[color=mycolor2, dashed, ultra thick] table[x=ALPHA, y=SOLVER2] {\datatable}; 
    \addlegendentry{\matpower{} case data}

    \addplot[color=mycolor3, dashdotted, ultra thick] table[x=ALPHA, y=SOLVER3] {\datatable}; 
    \addlegendentry{Power flow solution}
   
  \end{axis}
  
\end{tikzpicture}
		\caption{ Overall time \label{fig:profileInitBeltistos-t} }
	\end{subfigure}
	\begin{subfigure}[b]{0.49\columnwidth}
		\begin{tikzpicture}

\pgfplotstableread{
ALPHA         SOLVER1         SOLVER2         SOLVER3
1.0         0.32000         0.56000         0.24000
1.33         0.48000         0.92000         0.84000
1.6600000000000001         0.60000         0.92000         0.92000
1.9900000000000002         0.64000         1.00000         0.96000
2.3200000000000003         0.68000         1.00000         1.00000
2.6500000000000004         0.68000         1.00000         1.00000
2.9800000000000004         0.68000         1.00000         1.00000
3.3100000000000005         0.72000         1.00000         1.00000
3.6400000000000006         0.72000         1.00000         1.00000
3.9700000000000006         0.72000         1.00000         1.00000
4.300000000000001         0.72000         1.00000         1.00000
4.630000000000001         0.72000         1.00000         1.00000
4.960000000000001         0.76000         1.00000         1.00000
5.290000000000001         0.76000         1.00000         1.00000
5.620000000000001         0.80000         1.00000         1.00000
5.950000000000001         0.80000         1.00000         1.00000
6.280000000000001         0.84000         1.00000         1.00000
6.610000000000001         0.88000         1.00000         1.00000
6.940000000000001         0.96000         1.00000         1.00000
7.270000000000001         0.96000         1.00000         1.00000
7.600000000000001         0.96000         1.00000         1.00000
7.9300000000000015         0.96000         1.00000         1.00000
8.260000000000002         0.96000         1.00000         1.00000
8.590000000000002         0.96000         1.00000         1.00000
8.920000000000002         0.96000         1.00000         1.00000
9.250000000000002         0.96000         1.00000         1.00000
9.580000000000002         0.96000         1.00000         1.00000
9.910000000000002         0.96000         1.00000         1.00000
10.240000000000002         0.96000         1.00000         1.00000
10.570000000000002         0.96000         1.00000         1.00000
10.900000000000002         0.96000         1.00000         1.00000
11.230000000000002         0.96000         1.00000         1.00000
11.560000000000002         0.96000         1.00000         1.00000
11.90000000000001         0.96000         1.00000         1.00000
12.00000000000001         1.00000         1.00000         1.00000
13.00000000000001         1.00000         1.00000         1.00000
}\datatable

  \begin{axis}[name=symb,
  	width=\columnwidth,
  	height=4cm,
	ymin=-0.02,
	ymax=1.02,
    enlarge x limits=0.05,
	axis lines*=left, 
	ymajorgrids, yminorgrids,
    xmajorgrids, xminorgrids,
	xticklabel style={rotate=0, xshift=-0.0cm, anchor=north, font=\scriptsize},
	ytick={0, 0.1, 0.2, 0.3, 0.4, 0.5, 0.6, 0.7, 0.8, 0.9, 1},
    xmax=12.5,
	yticklabel style={font=\scriptsize},
	ylabel style={font=\scriptsize},
	xlabel style={font=\scriptsize},
	ylabel={$p_m({\alpha})$},
	xlabel={$\alpha$},
    legend style={at={(1.0,0.55)},legend cell align=left,align=right,draw=white!85!black,font=\footnotesize,legend columns=1}
	]
	
    \addplot[color=mycolor1, ultra thick] table[x=ALPHA, y=SOLVER1] {\datatable}; 
    \addlegendentry{Flat start}

    \addplot[color=mycolor2, dashed, ultra thick] table[x=ALPHA, y=SOLVER2] {\datatable}; 
    \addlegendentry{\matpower{} case data}

    \addplot[color=mycolor3, dashdotted, ultra thick] table[x=ALPHA, y=SOLVER3] {\datatable}; 
    \addlegendentry{Power flow solution}
   
  \end{axis}
  
\end{tikzpicture}
		\caption{ Number of iterations \label{fig:profileInitBeltistos-i} }
	\end{subfigure}
	\caption{Performance profiles for the three initial points  using \BELTISTOSOPF{} considering all benchmarks. \label{fig:profileInitBeltistos}}
\end{figure*}

\begin{figure*}[ht!]
	\centering
	\begin{subfigure}[b]{0.49\columnwidth}
		\begin{tikzpicture}

\pgfplotstableread{
ALPHA         SOLVER1         SOLVER2         SOLVER3         SOLVER4
1.0         0.32000         0.08000         0.60000         0.00000
1.2         0.68000         0.68000         0.76000         0.40000
1.4         0.72000         0.96000         0.88000         0.64000
1.5999999999999999         0.72000         0.96000         0.96000         0.72000
1.7999999999999998         0.80000         1.00000         0.96000         0.88000
1.9999999999999998         0.80000         1.00000         1.00000         0.88000
2.1999999999999997         0.84000         1.00000         1.00000         0.88000
2.3999999999999995         0.84000         1.00000         1.00000         0.88000
2.5999999999999996         0.84000         1.00000         1.00000         0.92000
2.8         0.84000         1.00000         1.00000         0.92000
2.9999999999999996         0.84000         1.00000         1.00000         0.92000
3.1999999999999993         0.84000         1.00000         1.00000         0.92000
3.3999999999999995         0.88000         1.00000         1.00000         0.96000
3.5999999999999996         0.92000         1.00000         1.00000         0.96000
3.7999999999999994         0.92000         1.00000         1.00000         0.96000
3.999999999999999         0.92000         1.00000         1.00000         0.96000
4.199999999999999         0.92000         1.00000         1.00000         0.96000
4.3999999999999995         0.92000         1.00000         1.00000         0.96000
4.6         0.92000         1.00000         1.00000         0.96000
4.799999999999999         0.92000         1.00000         1.00000         0.96000
4.999999999999999         0.92000         1.00000         1.00000         0.96000
5.199999999999999         0.92000         1.00000         1.00000         0.96000
5.399999999999999         0.92000         1.00000         1.00000         0.96000
5.599999999999999         0.92000         1.00000         1.00000         0.96000
5.799999999999999         0.92000         1.00000         1.00000         0.96000
5.999999999999999         0.92000         1.00000         1.00000         0.96000
6.199999999999999         0.92000         1.00000         1.00000         1.00000
6.399999999999999         0.92000         1.00000         1.00000         1.00000
6.599999999999999         0.92000         1.00000         1.00000         1.00000
6.799999999999999         0.92000         1.00000         1.00000         1.00000
6.999999999999998         0.92000         1.00000         1.00000         1.00000
7.199999999999998         0.96000         1.00000         1.00000         1.00000
7.399999999999999         0.96000         1.00000         1.00000         1.00000
7.599999999999999         0.96000         1.00000         1.00000         1.00000
7.799999999999999         0.96000         1.00000         1.00000         1.00000
7.999999999999998         0.96000         1.00000         1.00000         1.00000
8.2         0.96000         1.00000         1.00000         1.00000
8.399999999999999         0.96000         1.00000         1.00000         1.00000
8.599999999999998         0.96000         1.00000         1.00000         1.00000
8.799999999999997         0.96000         1.00000         1.00000         1.00000
8.999999999999998         0.96000         1.00000         1.00000         1.00000
}\datatable

  \begin{axis}[name=symb,
  	width=\columnwidth,
  	height=4cm,
	ymin=-0.02,
	ymax=1.02,
    enlarge x limits=0.05,
	axis lines*=left, 
	ymajorgrids, yminorgrids,
    xmajorgrids, xminorgrids,
	xticklabel style={rotate=0, xshift=-0.0cm, anchor=north, font=\scriptsize},
	ytick={0, 0.1, 0.2, 0.3, 0.4, 0.5, 0.6, 0.7, 0.8, 0.9, 1},
    xmax=9,
	yticklabel style={font=\scriptsize},
	ylabel style={font=\scriptsize},
	xlabel style={font=\scriptsize},
	ylabel={$p_m({\alpha})$},
	xlabel={$\alpha$},
    legend style={at={(1.0,0.55)},legend cell align=left,align=right,draw=white!85!black,font=\footnotesize,legend columns=1}
	]
	
    \addplot[color=mycolor1, ultra thick] table[x=ALPHA, y=SOLVER1] {\datatable}; 
    \addlegendentry{Flat start}

    \addplot[color=mycolor2, dashed, ultra thick] table[x=ALPHA, y=SOLVER2] {\datatable}; 
    \addlegendentry{\matpower{} case data}

    \addplot[color=mycolor3, dashdotted, ultra thick] table[x=ALPHA, y=SOLVER3] {\datatable}; 
    \addlegendentry{Power flow solution}
   
  \end{axis}
  
\end{tikzpicture}
		\caption{ Overall time \label{fig:profileInitKnitro-t} }
	\end{subfigure}
	\begin{subfigure}[b]{0.49\columnwidth}
		\begin{tikzpicture}

\pgfplotstableread{
ALPHA         SOLVER1         SOLVER2         SOLVER3         SOLVER4
1.0         0.48000         0.32000         0.24000         0.24000
1.2         0.64000         0.96000         0.72000         0.48000
1.4         0.64000         1.00000         0.76000         0.68000
1.5999999999999999         0.72000         1.00000         0.88000         0.84000
1.7999999999999998         0.72000         1.00000         0.92000         0.84000
1.9999999999999998         0.72000         1.00000         0.92000         0.96000
2.1999999999999997         0.76000         1.00000         1.00000         0.96000
2.3999999999999995         0.76000         1.00000         1.00000         0.96000
2.5999999999999996         0.84000         1.00000         1.00000         0.96000
2.8         0.84000         1.00000         1.00000         0.96000
2.9999999999999996         0.84000         1.00000         1.00000         0.96000
3.1999999999999993         0.84000         1.00000         1.00000         0.96000
3.3999999999999995         0.84000         1.00000         1.00000         1.00000
3.5999999999999996         0.84000         1.00000         1.00000         1.00000
3.7999999999999994         0.84000         1.00000         1.00000         1.00000
3.999999999999999         0.84000         1.00000         1.00000         1.00000
4.199999999999999         0.88000         1.00000         1.00000         1.00000
4.3999999999999995         0.88000         1.00000         1.00000         1.00000
4.6         0.88000         1.00000         1.00000         1.00000
4.799999999999999         0.88000         1.00000         1.00000         1.00000
4.999999999999999         0.92000         1.00000         1.00000         1.00000
5.199999999999999         0.92000         1.00000         1.00000         1.00000
5.399999999999999         0.92000         1.00000         1.00000         1.00000
5.599999999999999         0.92000         1.00000         1.00000         1.00000
5.799999999999999         0.92000         1.00000         1.00000         1.00000
5.999999999999999         0.92000         1.00000         1.00000         1.00000
6.199999999999999         0.92000         1.00000         1.00000         1.00000
6.399999999999999         0.92000         1.00000         1.00000         1.00000
6.599999999999999         0.92000         1.00000         1.00000         1.00000
6.799999999999999         0.92000         1.00000         1.00000         1.00000
6.999999999999998         0.92000         1.00000         1.00000         1.00000
7.199999999999998         0.92000         1.00000         1.00000         1.00000
7.399999999999999         0.92000         1.00000         1.00000         1.00000
7.599999999999999         0.92000         1.00000         1.00000         1.00000
7.799999999999999         0.92000         1.00000         1.00000         1.00000
7.999999999999998         0.92000         1.00000         1.00000         1.00000
8.2         0.92000         1.00000         1.00000         1.00000
8.399999999999999         0.92000         1.00000         1.00000         1.00000
8.599999999999998         0.96000         1.00000         1.00000         1.00000
8.799999999999997         0.96000         1.00000         1.00000         1.00000
8.999999999999998         0.96000         1.00000         1.00000         1.00000
}\datatable

  \begin{axis}[name=symb,
  	width=\columnwidth,
  	height=4cm,
	ymin=-0.02,
	ymax=1.02,
    enlarge x limits=0.05,
	axis lines*=left, 
	ymajorgrids, yminorgrids,
    xmajorgrids, xminorgrids,
	xticklabel style={rotate=0, xshift=-0.0cm, anchor=north, font=\scriptsize},
	ytick={0, 0.1, 0.2, 0.3, 0.4, 0.5, 0.6, 0.7, 0.8, 0.9, 1},
    xmax=9,
	yticklabel style={font=\scriptsize},
	ylabel style={font=\scriptsize},
	xlabel style={font=\scriptsize},
	ylabel={$p_m({\alpha})$},
	xlabel={$\alpha$},
    legend style={at={(1.0,0.55)},legend cell align=left,align=right,draw=white!85!black,font=\footnotesize,legend columns=1}
	]
	
    \addplot[color=mycolor1, ultra thick] table[x=ALPHA, y=SOLVER1] {\datatable}; 
    \addlegendentry{Flat start}

    \addplot[color=mycolor2, dashed, ultra thick] table[x=ALPHA, y=SOLVER2] {\datatable}; 
    \addlegendentry{\matpower{} case data}

    \addplot[color=mycolor3, dashdotted, ultra thick] table[x=ALPHA, y=SOLVER3] {\datatable}; 
    \addlegendentry{Power flow solution}
   
  \end{axis}
  
\end{tikzpicture}
		\caption{ Number of iterations \label{fig:profileInitKnitro-i} }
	\end{subfigure}
	\caption{Performance profiles for the three initial points  using \KNITRO{} considering all benchmarks. \label{fig:profileInitKnitro}}
\end{figure*}
%
%
\pgfplotstableread{
	Optimizer flat mpc pfsolve Total
	BELTISTOS-OPF      25 25 25 75/75
	KNITRO11             25 24 24  73/75
	IPOPT-PARDISO      23 22 24 69/75
	IPOPT-MA57         23 22 21 66/75
	MIPS-PARDISO       16 24 24 64/75
	MIPS-MATLAB'\textbackslash'       14 22 22 58/75
	FMINCON              17 19 20 56/75
}\initSuccess
%
\begin{table}[ht!]
\footnotesize
\centering
\caption{Number of solved benchmarks out of twenty-five test cases for different starting points. Gray background indicates use of the PARDISO linear solver. \label{tab:profileGuess}}
\pgfplotstabletypeset[
    every head row/.style={ before row=\toprule,after row=\midrule},
    every last row/.style={ after row=\bottomrule},
    precision=0, fixed zerofill, column type={c},
    columns/Optimizer/.style={string type,column type=l},
    every row no 1/.style={ before row={\rowcolor{tablecolor1}}},
    every row no 2/.style={ before row={\rowcolor{tablecolor1}}},
    every row no 4/.style={ before row={\rowcolor{tablecolor1}}},
    columns/flat/.style={column name={Flat start}},
    columns/mpc/.style={column name={\matpower{} case data}},
    columns/pfsolve/.style={column name={Power flow solution}},
    columns/Total/.style={string type,column type=r}
]\initSuccess
\end{table}

\begin{figure}[t!]
	\centering
	\begin{tikzpicture}




\pgfplotstableread{
id case flat mpc pfsolve
1 1951rte      160	24	25
2 2383wp       25	48	25
3 2736sp       20	13	15
4 2737sop      20	12	15
5 2746wop      20	9	18
6 2746wp       20	15	16
7 2868rte      164	26	28
8 2869pegase   29	22	19
9 3012wp       22	40	36
10 3120sp      26	29	29
11 3375wp      24	27	31
12 6468rte     192	34	29
13 6470rte     191	34	34
14 6495rte     137	50	45
15 6515rte     177	38	65
16 9241pegase  30	25	26
17 ACTIVSg2000 204	18	17
18 ACTIVSg10k  108	18	22
19 13659pegase 51	51	75
20 ACTIVSg25k  42	32	36
21 ACTIVSg70k  45	38	41
22 case21k     46	55	53
23 case42k     51	55	59
24 case99k     66	82	70
25 case193k    91	90	96
}\datatable

\begin{axis}[
width=\columnwidth, height=4cm, 
ymajorgrids, yminorgrids,
xmajorgrids, xminorgrids,
area style, 
smooth,
enlarge x limits=0.0,
xtick=data,
ymin=0,
yticklabel style={font=\footnotesize},
xticklabels from table={\datatable}{case},
xticklabel style={font=\footnotesize, rotate=40},
ylabel=Number of Iterations,
ylabel style={font=\footnotesize},
legend style={at={(0,0.999)},
	anchor=north west,legend columns=2,fill=none,draw=none, font=\footnotesize},
axis on top=false
]

\addplot[fill=mycolor1, color=mycolor1] table[x=id, y=flat] {\datatable} \closedcycle;
\addlegendentry{Flat start}
\addplot[fill=mycolor3,color=mycolor3] table[x=id, y=pfsolve] {\datatable} \closedcycle;
\addlegendentry{Power flow}
\addplot[fill=mycolor2,color=mycolor2] table[x=id, y=mpc] {\datatable} \closedcycle;
\addlegendentry{\matpower{} case data}

\addplot[color=mycolor1, dashed, thick] table[x=id, y=flat] {\datatable};
\addplot[color=mycolor2, dotted, thick] table[x=id, y=mpc] {\datatable};
\addplot[color=mycolor3, dashdotted, thick] table[x=id, y=pfsolve] {\datatable};

\end{axis}
\end{tikzpicture}
	\vspace{-5mm}
	\caption{Number of iterations for various initial points (\BELTISTOSOPF{}) \label{fig:OPFstartBeltistos}}
\end{figure}

\pgfplotstableread{
case flat mpc pfsolve  
case1951rte             17.380  4.200   5.170
case2383wp              29.230  9.030   8.080
case2736sp              30.310  4.210   3.630
case2737sop             4.830   3.380   4.160
case2746wop             5.870   3.140   4.170
case2746wp              4.760   3.020   3.930
case2868rte             15.210  5.970   6.160
case2869pegase          9.660   10.220  15.940
case3012wp              9.110   36.550  9.370
case3120sp              9.510   9.750   8.650
case3375wp              8.330   11.010  18.880
case6468rte             80.930  28.620  27.030
case6470rte             104.450 62.960  85.710
case6495rte             15.170  14.320  13.620
case6515rte             49.620  11.540  11.770
case9241pegase          235.390 118.080 52.520
case\_ACTIVSg2000       41.080  3.530   3.990
case\_ACTIVSg10k        125.450 19.320  121.120
case13659pegase         {} {} 180.440
case\_ACTIVSg25k        57.600  55.030  362.330
case\_ACTIVSg70k        332.800 {} 345.750
case21k                 161.490 323.110 330.290
case42k                 629.440 850.520 2024.700
case99k                 4192.140    4983.650    3852.760
case193k                {} {} {}
}\OPFtabletimeIPOPTdefault

\pgfplotstableread{
case flat mpc pfsolve  
case1951rte             238 43  53
case2383wp              177 51  43
case2736sp              309 30  22
case2737sop             34  21  27
case2746wop             45  18  28
case2746wp              37  17  26
case2868rte             136 41  50
case2869pegase          46  47  72
case3012wp              41  192 40
case3120sp              41  44  37
case3375wp              35  47  86
case6468rte             414 162 146
case6470rte             486 284 393
case6495rte             68  73  73
case6515rte             274 58  60
case9241pegase          335 178 65
case\_ACTIVSg2000       212 26  29
case\_ACTIVSg10k        287 49  345
case13659pegase         {} {} 244
case\_ACTIVSg25k        58  62  377
case\_ACTIVSg70k        102 {} 104
case21k                 70  145 114
case42k                 70  73  180
case99k                 91  95  94
case193k                {} {} {}
}\OPFtableitersIPOPTdefault

\pgfplotstableread{
	case flat mpc pfsolve  
case1951rte	14.9	2.99	3.3
case2383wp	3.92	6.58	4.56
case2736sp	2.98	3	2.95
case2737sop	3.23	2.53	2.99
case2746wop	3.07	1.99	3.38
case2746wp	2.94	2.78	3.7
case2868rte	19.32	3.95	4.77
case2869pegase	5.16	4.2	4
case3012wp	5.24	5.99	5.99
case3120sp	5.16	5.62	5.65
case3375wp	4.32	4.71	5.45
case6468rte	39.35	8	7.45
case6470rte	43.61	8.7	8.79
case6495rte	30.83	10.55	11.22
case6515rte	40.48	9.07	15.63
case9241pegase	13.16	11.66	12.68
case\_ACTIVSg2000	41.13	2.92	3.06
case\_ACTIVSg10k	49.1	10.6	10.28
case13659pegase	21.32	20.86	25.51
case\_ACTIVSg25k	44.48	37.17	43.04
case\_ACTIVSg70k	148.77	112.92	152.38
case21k	56.81	70.09	65.27
case42k	177.02	198.13	220.83
case99k	1387.91	1402.08	1231.9
case193k	4397.02	3832.95	4360.57
}\OPFtabletimeBELTISTOS

\pgfplotstableread{
case flat mpc pfsolve  
case1951rte	160	24	25
case2383wp	25	48	25
case2736sp	20	13	15
case2737sop	20	12	15
case2746wop	20	9	18
case2746wp	20	15	16
case2868rte	164	26	28
case2869pegase	29	22	19
case3012wp	22	40	36
case3120sp	26	29	29
case3375wp	24	27	31
case6468rte	192	34	29
case6470rte	191	34	34
case6495rte	137	50	45
case6515rte	177	38	65
case9241pegase	30	25	26
case\_ACTIVSg2000	204	18	17
case\_ACTIVSg10k	108	18	22
case13659pegase	51	51	75
case\_ACTIVSg25k	42	32	36
case\_ACTIVSg70k	45	38	41
case21k	46	55	53
case42k	51	55	59
case99k	66	82	70
case193k	91	90	96
}\OPFtableitersBELTISTOS


\pgfplotstableread{
case flat mpc pfsolve  
case1951rte             12.09   3.31    3.4
case2383wp              3.58    4.78    4.11
case2736sp              3.73    2.47    3.24
case2737sop             3.54    2.65    3.64
case2746wop             3.46    2.32    3.69
case2746wp              3.73    2.39    3.97
case2868rte             16.01   4.71    4.77
case2869pegase          4.68    3.91    3.89
case3012wp              5.18    5.23    5.57
case3120sp              4.78    4.66    4.86
case3375wp              4.27    4.47    5.4
case6468rte             24.75   6.29    6.65
case6470rte             33.64   12.18   12.51
case6495rte             35.49   9.1 10.17
case6515rte             23.15   10.59   9.86
case9241pegase          15.05   20.42   17.76
case\_ACTIVSg2000       44.24   3.44    4
case\_ACTIVSg10k        38.47   11.58   13.99
case13659pegase         83.52   71.93   86.76
case\_ACTIVSg25k        54.34   46.75   50.06
case\_ACTIVSg70k        167.8   144.37  152.56
case21k                 120.91  598.92  {}
case42k                 7918.16 {} {}
case99k                 {} {} {}
case193k                {} {} {}
}\OPFtabletimeIPOPTma

\pgfplotstableread{
case flat mpc pfsolve  
case1951rte             177 36  34
case2383wp              36  54  39
case2736sp              30  17  23
case2737sop             31  20  28
case2746wop             30  15  29
case2746wp              33  16  28
case2868rte             160 41  39
case2869pegase          36  31  29
case3012wp              43  49  48
case3120sp              40  40  38
case3375wp              34  36  43
case6468rte             169 38  38
case6470rte             206 70  68
case6495rte             224 56  60
case6515rte             142 62  57
case9241pegase          41  54  45
case\_ACTIVSg2000       212 26  29
case\_ACTIVSg10k        99  31  35
case13659pegase         225 200 246
case\_ACTIVSg25k        56  51  50
case\_ACTIVSg70k        72  62  65
case21k                 80  299 {}
case42k                 181 {} {}
case99k                 {} {} {}
case193k                {} {} {}
}\OPFtableitersIPOPTma

\pgfplotstableread{
case flat mpc pfsolve  
case1951rte             {} 3.51    3.71
case2383wp              5   4.71    5.14
case2736sp              4.68    4.28    4.44
case2737sop             4.51    4.55    4.57
case2746wop             4.99    4.39    4.94
case2746wp              5.44    4.68    4.85
case2868rte             {} 5.13    5.28
case2869pegase          7.22    24.83   6.11
case3012wp              7.81    5.16    5.5
case3120sp              7.45    21.61   6.36
case3375wp              8.77    5.72    5.91
case6468rte             {} 13.95   14.18
case6470rte             {} 13.43   14.16
case6495rte             {} 23.34   24.29
case6515rte             {} 18.42   19.36
case9241pegase          23.84   36.05   25.9
case\_ACTIVSg2000       4.9 3.84    4.38
case\_ACTIVSg10k        {} 20.26   45.32
case13659pegase         36.51   24.42   30.53
case\_ACTIVSg25k        {} 72.43   108.18
case\_ACTIVSg70k        {} {} {}
case21k                 291.87  238.12  241.2
case42k                 2264.12 1925.02 1957.95
case99k                 {} {} {}
case193k                {} {} {}
}\OPFtabletimeMIPSsc

\pgfplotstableread{
case flat mpc pfsolve  
case1951rte             {} 28  28
case2383wp              35  32  35
case2736sp              30  26  26
case2737sop             28  28  26
case2746wop             31  26  29
case2746wp              31  28  28
case2868rte             {} 32  31
case2869pegase          39  127 32
case3012wp              44  29  30
case3120sp              43  111 33
case3375wp              47  30  30
case6468rte             {} 47  44
case6470rte             {} 43  44
case6495rte             {} 74  76
case6515rte             {} 59  61
case9241pegase          45  68  46
case\_ACTIVSg2000       35  26  29
case\_ACTIVSg10k        {} 41  81
case13659pegase         65  42  54
case\_ACTIVSg25k        {} 56  76
case\_ACTIVSg70k        {} {} {}
case21k                 71  53  53
case42k                 83  64  64
case99k                 {} {} {}
case193k                {} {} {}
}\OPFtableitersMIPSsc


\pgfplotstableread{
case flat mpc pfsolve  
case1951rte             {} 4.54    4.75
case2383wp              6.49    6.02    6.6
case2736sp              6.04    5.47    5.98
case2737sop             5.7 5.83    5.68
case2746wop             6.63    5.92    6.4
case2746wp              6.53    6.1 6.26
case2868rte             {} 6.73    6.82
case2869pegase          9.04    23.44   7.66
case3012wp              9.86    6.52    7.13
case3120sp              9.4 28.23   7.6
case3375wp              11.23   7.2 7.37
case6468rte             {} 17.25   16.57
case6470rte             {} 16.84   17.31
case6495rte             {} 31.68   30.8
case6515rte             {} 23.51   24.71
case9241pegase          30.67   47.08   31.91
case\_ACTIVSg2000       6.41    5.04    5.7
case\_ACTIVSg10k        {} 28.07   60.06
case13659pegase         51.19   33.55   43.54
case\_ACTIVSg25k        {} 93.24   138.79
case\_ACTIVSg70k        {} {} {}
case21k                 135.4   97.4    98
case42k                 390.42  289.74  303.38
case99k                 1329.69 1082.13 1081.72
case193k                3703.2  2953.02 2982.83
}\OPFtabletimeMIPSscPardiso

\pgfplotstableread{
case flat mpc pfsolve  
case1951rte             {} 28  28
case2383wp              35  32  35
case2736sp              30  26  26
case2737sop             28  28  26
case2746wop             31  26  29
case2746wp              31  28  28
case2868rte             {} 32  31
case2869pegase          39  98  32
case3012wp              44  29  30
case3120sp              43  117 33
case3375wp              47  30  30
case6468rte             {} 47  44
case6470rte             {} 43  44
case6495rte             {} 77  76
case6515rte             {} 59  61
case9241pegase          45  68  46
case\_ACTIVSg2000       35  26  29
case\_ACTIVSg10k        {} 41  81
case13659pegase         64  42  54
case\_ACTIVSg25k        {} 56  76
case\_ACTIVSg70k        {} {} {}
case21k                 71  53  53
case42k                 83  64  64
case99k                 97  78  78
case193k                113 92  92
}\OPFtableitersMIPSscPardiso

\pgfplotstableread{
case flat mpc pfsolve
case1951rte      {} 37.8183330000 19.2719290000
case2383wp      30.3933380000 60.7413180000 64.5316740000
case2736sp      14.3935180000 {} 13.4904430000
case2737sop      10.7012570000 12.0072640000 9.0151000000
case2746wop      11.4215000000 19.9374390000 12.3532740000
case2746wp      13.0845840000 140.8455330000 12.6106050000
case2868rte      {} 49.1610210000 37.2584490000
case2869pegase      10.3962000000 22.0460560000 9.4018520000
case3012wp      33.2662160000 404.8784220000 {}
case3120sp      31.6527090000 {} 17.8540070000
case3375wp      35.7352800000 {} {}
case6468rte      {} {} {}
case6470rte      {} 92.6635830000 87.6310780000
case6495rte      {} 31.3131050000 34.6605350000
case6515rte      {} 108.8577350000 79.8955360000
case9241pegase      149.0412010000 304.8070510000 191.8276190000
case\_ACTIVSg2000     56.8867600000 10.6183720000 12.3843980000
case\_ACTIVSg10k      {} 19.3042190000 58.3209150000
case13659pegase      1259.1706270000 2601.3706430000 2960.2958620000
case\_ACTIVSg25k      481.6979450000 214.4435540000 272.0320370000
case\_ACTIVSg70k      1359.4226700000 836.1880860000 923.4247320000
case21k      669.9864010000 2742.7294370000 2228.1194000000
case42k      1834.8515810000 {} 6163.0426590000
case99k      26188.50 56236.69 62359.75
case193k     49073.79 {} {}
}\OPFtabletimeFMINCONa

\pgfplotstableread{
case flat mpc pfsolve
case1951rte      {} 77 62
case2383wp      111 184 85
case2736sp      46 {} 44
case2737sop      33 42 31
case2746wop      33 61 39
case2746wp      39 333 41
case2868rte      {} 58 61
case2869pegase      27 50 28
case3012wp      91 434 {}
case3120sp      90 {} 52
case3375wp      82 {} {}
case6468rte      {} {} {}
case6470rte      {} 57 45
case6495rte      {} 45 47
case6515rte      {} 108 74
case9241pegase      39 81 44
case\_ACTIVSg2000     67 35 40
case\_ACTIVSg10k      {} 18 52
case13659pegase      65 63 63
case\_ACTIVSg25k      116 48 48
case\_ACTIVSg70k      114 79 82
case21k      128 437 396
case42k      135 {} 386
case99k      141 297 354
case193k     153 {} {}
}\OPFtableitersFMINCONa

\pgfplotstableread{
case flat mpc pfsolve
case1951rte       {} 42.13 20.69
case2383wp       33.49 69.35 66.65
case2736sp       15.75 324.89 16.20
case2737sop       11.84 14.92 11.93
case2746wop       12.68 24.39 15.46
case2746wp       14.30 165.37 15.86
case2868rte       {} 54.29 41.41
case2869pegase       11.18 25.07 12.53
case3012wp       36.38 {} {}
case3120sp       36.06 {} 21.58
case3375wp       40.65 183.98 {}
case6468rte       {} {} {}
case6470rte       {} 92.65 97.92
case6495rte       {} 37.76 38.41
case6515rte       {} 123.73 80.07
case9241pegase       141.15 281.36 186.36
case\_ACTIVSg2000     52.67 12.81 14.06
case\_ACTIVSg10k      {} 26.88 61.70
case13659pegase       1848.93 1866.80 2269.08
case\_ACTIVSg25k      442.05 210.12 254.40
case\_ACTIVSg70k      1280.10 786.05 874.72
case21k       590.86 {} 2311.17
case42k       1624.20 8428.31 {}
case99k       8330.90 20554.62 24267.84
case193k      16680.00 70150.87 {}
}\OPFtabletimeFMINCONb

\pgfplotstableread{
case flat mpc pfsolve
case1951rte      {} 77 62
case2383wp      111 184 85
case2736sp      46 455 44
case2737sop      33 42 31
case2746wop      33 61 39
case2746wp      39 343 41
case2868rte      {} 58 61
case2869pegase      27 50 28
case3012wp      91 {} {}
case3120sp      90 {} 52
case3375wp      82 383 {}
case6468rte      {} {} {}
case6470rte      {} 57 45
case6495rte      {} 45 47
case6515rte      {} 108 74
case9241pegase      39 81 44
case\_ACTIVSg2000    67 35 40
case\_ACTIVSg10k     {} 18 52
case13659pegase      62 81 99
case\_ACTIVSg25k     116 48 48
case\_ACTIVSg70k     114 79 82
case21k      128 {} 383
case42k      135 435 {}
case99k      141 324 359
case193k     153 482 {}
}\OPFtableitersFMINCONb

\pgfplotstableread{
case flat mpc pfsolve  
case1951rte             {}   12.49   9.29
case2383wp              31.23   24.93   13.25
case2736sp              9.69    41.85   10.7
case2737sop             7.75    9.34    7.59
case2746wop             8.51    14.16   9.28
case2746wp              9.42    59.13   10.47
case2868rte             {}   70.72   22.09
case2869pegase          8.14    14.21   8.95
case3012wp              21.26   {}   {}
case3120sp              20.25   {}   14
case3375wp              22.39   {}   {}
case6468rte             {}   {}   {}
case6470rte             {}   108.83  110.68
case6495rte             {}   25  25.56
case6515rte             {}   35.77   35.61
case9241pegase          63.98   223.98  85.86
case\_ACTIVSg2000       26.67   9.76    9.44
case\_ACTIVSg10k        {}   18.88   38.61
case13659pegase         3091.9  2220.37 3975.64
case\_ACTIVSg25k        266.76  135.41  183.57
case\_ACTIVSg70k        869.9   635.57  644.3
case21k                 466.2   {}   1548.28
case42k                 1511.98 4721.76 {}
case99k                 {}   17013.29    17325.98
case193k                17674.56    {}   {}
}\OPFtabletimeFMINCONc

\pgfplotstableread{
case flat mpc pfsolve  
case1951rte             {}   59  48
case2383wp              92  124 62
case2736sp              44  114 46
case2737sop             34  43  31
case2746wop             36  62  39
case2746wp              39  273 45
case2868rte             {}   91  62
case2869pegase          33  50  34
case3012wp              93  {}   {}
case3120sp              90  {}   56
case3375wp              83  {}   {}
case6468rte             {}   {}   {}
case6470rte             {}   80  67
case6495rte             {}   57  60
case6515rte             {}   73  70
case9241pegase          39  67  46
case\_ACTIVSg2000       74  43  41
case\_ACTIVSg10k        {}   25  52
case13659pegase         137 103 77
case\_ACTIVSg25k        116 53  65
case\_ACTIVSg70k        114 86  83
case21k                 143 {}   415
case42k                 154 475 {}
case99k                 {}   315 348
case193k                173 {}   {}
}\OPFtableitersFMINCONc
\pgfplotstableread{
case flat mpc pfsolve
case1951rte            {} 10.22 11.44
case2383wp             5.09 4.92 4.84
case2736sp             4.05 {} 4.36
case2737sop            4.01 3.97 4.40
case2746wop            3.94 45.41 4.56
case2746wp             4.11 13.07 4.65
case2868rte            {} {} 23.06
case2869pegase         4.53 4.74 4.95
case3012wp             5.15 7.02 7.86
case3120sp             5.15 5.35 5.06
case3375wp             5.48 6.73 7.38
case6468rte            94.71 42.25 45.22
case6470rte            {} 10.17 9.91
case6495rte            {} 9.33 9.50
case6515rte            114.61 12.58 13.54
case9241pegase         101.92 16.87 14.42
case\_ACTIVSg2000      4.16 4.37 4.70
case\_ACTIVSg10k       {} 12.47 12.37
case13659pegase        166.12 93.37 160.42
case\_ACTIVSg25k       49.66 46.64 48.47
case\_ACTIVSg70k       163.39 162.04 183.33
case21k                51.71 206.98 224.30
case42k                164.11 1147.66 {}
case99k                785.26 2398.32 7541.11
case193k               1038.66 2778.65 2914.72
}\OPFtabletimeKNITROa

\pgfplotstableread{
case flat mpc pfsolve
case1951rte             {} 97 107
case2383wp              34 32 29
case2736sp              20 {} 20
case2737sop             20 20 20
case2746wop             20 440 22
case2746wp              20 117 20
case2868rte             {} {} 169
case2869pegase          22 25 23
case3012wp              28 48 55
case3120sp              27 28 24
case3375wp              28 39 42
case6468rte             360 185 194
case6470rte             {} 33 33
case6495rte             {} 36 36
case6515rte             469 46 46
case9241pegase          188 32 28
case\_ACTIVSg2000       18 21 20
case\_ACTIVSg10k        {} 23 25
case13659pegase         363 203 355
case\_ACTIVSg25k        46 45 45
case\_ACTIVSg70k        52 52 52
case21k                 43 181 163
case42k                 50 325 {}
case99k                 60 150 437
case193k                59 145 156
}\OPFtableitersKNITROa
\pgfplotstableread{
case flat mpc pfsolve  
case1951rte             6.59    4.99    5.24
case2383wp              3.87    3.67    3.82
case2736sp              3   40.3    3.23
case2737sop             2.97    3.01    3.51
case2746wop             3.15    19.58   3.44
case2746wp              3.09    18.49   3.44
case2868rte             36.14   6.96    5.91
case2869pegase          3.39    3.85    3.79
case3012wp              3.96    5.69    7.02
case3120sp              4.03    4.01    3.92
case3375wp              4.19    5.41    5.51
case6468rte             27.52   24.86   27.75
case6470rte             42.58   7.61    8.11
case6495rte             14.05   8.43    8.78
case6515rte             18.05   9.26    8.86
case9241pegase          106.27  23.62   13.42
case\_ACTIVSg2000       3.09    3.15    3.4
case\_ACTIVSg10k        151.69  11.46   11.73
case13659pegase         80.51   60.03   67.05
case\_ACTIVSg25k        50.13   46.78   48.36
case\_ACTIVSg70k        186.11  201.98  189.45
case21k                 55.73   173.18  189.18
case42k                 178.61  {} 1239.42
case99k                 905.67  2310.66 {}
case193k                1483.28 3636.18 3095.06
}\OPFtabletimeKNITROb

\pgfplotstableread{
case flat mpc pfsolve  
case1951rte             75  54  54
case2383wp              35  32  31
case2736sp              20  449 19
case2737sop             20  20  22
case2746wop             20  226 22
case2746wp              20  178 21
case2868rte             326 65  47
case2869pegase          22  27  23
case3012wp              28  49  60
case3120sp              27  28  24
case3375wp              28  40  40
case6468rte             131 134 152
case6470rte             192 31  32
case6495rte             65  40  38
case6515rte             72  41  40
case9241pegase          113 55  28
case\_ACTIVSg2000       19  20  20
case\_ACTIVSg10k        415 24  26
case13659pegase         165 110 131
case\_ACTIVSg25k        46  45  45
case\_ACTIVSg70k        53  53  52
case21k                 44  161 172
case42k                 50  {} 396
case99k                 60  168 {}
case193k                59  178 157
}\OPFtableitersKNITROb

\begin{minipage}{0.45\textwidth}
    \begin{table}[H]
        \footnotesize
    	\centering
    	\caption{Overall time (s) (\BELTISTOS{}). \label{tab:OPFstartIPOPT-t}}
    	\pgfplotstabletypeset[
    	every head row/.style={ before row={\toprule},after row=\midrule},
    	every last row/.style={ after row=\bottomrule},
    	precision=2, fixed zerofill, column type={r},
    	every row/.style={
    		column type=r,
    		dec sep align,
    		fixed,
    		fixed zerofill,
    	},
    	columns={case, flat, mpc, pfsolve},
    	every odd row/.style={before row={\rowcolor{tablecolor1}}},
        every even row/.style={before row={\rowcolor{tablecolor2}}},
    	columns/case/.style={string type, column type={l}, column name={Benchmark}},
    	columns/flat/.style={column name={Flat}},
    	columns/mpc/.style={column name={MPC}},
    	columns/pfsolve/.style={column name={PF}},
    	empty cells with={---} 
    	]\OPFtabletimeBELTISTOS
    \end{table}
\end{minipage}
\qquad
\begin{minipage}{0.45\textwidth} 
    \begin{table}[H]
        \footnotesize
    	\centering
    	\caption{ Number of iterations (\BELTISTOS{}). \label{tab:OPFstartIPOPT-i}}
    	\pgfplotstabletypeset[
    	every head row/.style={ before row={\toprule},after row=\midrule},
    	every last row/.style={ after row=\bottomrule},
    	precision=0, fixed zerofill, column type={r},
    	every row/.style={
    		column type=r,
    		dec sep align,
    		fixed,
    		fixed zerofill,
    	},
    	columns={case, flat, mpc, pfsolve},
    	every odd row/.style={before row={\rowcolor{tablecolor1}}},
        every even row/.style={before row={\rowcolor{tablecolor2}}},
    	columns/case/.style={string type, column type={l}, column name={Benchmark}},
    	columns/flat/.style={column name={Flat}},
    	columns/mpc/.style={column name={MPC}},
    	columns/pfsolve/.style={column name={PF}},
    	empty cells with={---} 
    	]\OPFtableitersBELTISTOS
    \end{table}
\end{minipage}  

\begin{minipage}{0.45\textwidth}
    \begin{table}[H]
        \footnotesize
    	\centering
    	\caption{Overall time (s) (\IPOPT{}-\PARDISO{}). \label{tab:OPFstartIPOPTdefault-t}}
    	\pgfplotstabletypeset[
    	every head row/.style={ before row={\toprule},after row=\midrule},
    	every last row/.style={ after row=\bottomrule},
    	precision=2, fixed zerofill, column type={r},
    	every row/.style={
    		column type=r,
    		dec sep align,
    		fixed,
    		fixed zerofill,
    	},
    	columns={case, flat, mpc, pfsolve},
    	every odd row/.style={before row={\rowcolor{tablecolor1}}},
        every even row/.style={before row={\rowcolor{tablecolor2}}},
    	columns/case/.style={string type, column type={l}, column name={Benchmark}},
    	columns/flat/.style={column name={Flat}},
    	columns/mpc/.style={column name={MPC}},
    	columns/pfsolve/.style={column name={PF}},
    	empty cells with={---} 
    	]\OPFtabletimeIPOPTdefault
    \end{table}
\end{minipage}
\qquad
\begin{minipage}{0.45\textwidth} 
    \begin{table}[H]
        \footnotesize
    	\centering
    	\caption{Number of iterations (\IPOPT{}-\PARDISO{}). \label{tab:OPFstartIPOPTdefault-i}}
    	\pgfplotstabletypeset[
    	every head row/.style={ before row={\toprule},after row=\midrule},
    	every last row/.style={ after row=\bottomrule},
    	precision=0, fixed zerofill, column type={r},
    	every row/.style={
    		column type=r,
    		dec sep align,
    		fixed,
    		fixed zerofill,
    	},
    	columns={case, flat, mpc, pfsolve},
    	every odd row/.style={before row={\rowcolor{tablecolor1}}},
        every even row/.style={before row={\rowcolor{tablecolor2}}},
    	columns/case/.style={string type, column type={l}, column name={Benchmark}},
    	columns/flat/.style={column name={Flat}},
    	columns/mpc/.style={column name={MPC}},
    	columns/pfsolve/.style={column name={PF}},
    	empty cells with={---} 
    	]\OPFtableitersIPOPTdefault
    \end{table}
\end{minipage}   

\begin{minipage}{0.45\textwidth}
    \begin{table}[H]
        \footnotesize
    	\centering
    	\caption{Overall time (s) (\IPOPT{}-MA57). \label{tab:OPFstartIPOPTma57-t}}
    	\pgfplotstabletypeset[
    	every head row/.style={ before row={\toprule},after row=\midrule},
    	every last row/.style={ after row=\bottomrule},
    	precision=2, fixed zerofill, column type={r},
    	every row/.style={
    		column type=r,
    		dec sep align,
    		fixed,
    		fixed zerofill,
    	},
    	columns={case, flat, mpc, pfsolve},
    	every odd row/.style={before row={\rowcolor{tablecolor1}}},
        every even row/.style={before row={\rowcolor{tablecolor2}}},
    	columns/case/.style={string type, column type={l}, column name={Benchmark}},
    	columns/flat/.style={column name={Flat}},
    	columns/mpc/.style={column name={MPC}},
    	columns/pfsolve/.style={column name={PF}},
    	empty cells with={---} 
    	]\OPFtabletimeIPOPTma
    \end{table}
\end{minipage}
\qquad
\begin{minipage}{0.45\textwidth} 
    \begin{table}[H]
        \footnotesize
    	\centering
    	\caption{Number of iterations (\IPOPT{}-MA57). \label{tab:OPFstartIPOPTma57-i}}
    	\pgfplotstabletypeset[
    	every head row/.style={ before row={\toprule},after row=\midrule},
    	every last row/.style={ after row=\bottomrule},
    	precision=0, fixed zerofill, column type={r},
    	every row/.style={
    		column type=r,
    		dec sep align,
    		fixed,
    		fixed zerofill,
    	},
    	columns={case, flat, mpc, pfsolve},
    	every odd row/.style={before row={\rowcolor{tablecolor1}}},
        every even row/.style={before row={\rowcolor{tablecolor2}}},
    	columns/case/.style={string type, column type={l}, column name={Benchmark}},
    	columns/flat/.style={column name={Flat}},
    	columns/mpc/.style={column name={MPC}},
    	columns/pfsolve/.style={column name={PF}},
    	empty cells with={---} 
    	]\OPFtableitersIPOPTma
    \end{table}
\end{minipage}  

\begin{minipage}{0.45\textwidth}    
    \begin{table}[H]
        \footnotesize
    	\centering
    	\caption{Overall time (s) (\MIPS{}-MATLAB'\textbackslash'). \label{tab:OPFstartMIPSsc-t}}
    	\pgfplotstabletypeset[
    	every head row/.style={ before row={\toprule},after row=\midrule},
    	every last row/.style={ after row=\bottomrule},
    	precision=2, fixed, fixed zerofill, column type={r},
    	every row/.style={
    		column type=r,
    		dec sep align,
    		fixed,
    		fixed zerofill,
    	},
    	columns={case, flat, mpc, pfsolve},
    	every odd row/.style={before row={\rowcolor{tablecolor1}}},
        every even row/.style={before row={\rowcolor{tablecolor2}}},
    	columns/case/.style={string type, column type={l}, column name={Benchmark}},
    	columns/flat/.style={column name={Flat}},
    	columns/mpc/.style={column name={MPC}},
    	columns/pfsolve/.style={column name={PF}},
    	empty cells with={---} 
    	]\OPFtabletimeMIPSsc
    \end{table}
\end{minipage}
\qquad
\begin{minipage}{0.45\textwidth}    
    \begin{table}[H]
        \footnotesize
    	\centering
    	\caption{Number of iterations (\MIPS{}-MATLAB'\textbackslash'). \label{tab:OPFstartMIPSsc-i}}
    	\pgfplotstabletypeset[
    	every head row/.style={ before row={\toprule},after row=\midrule},
    	every last row/.style={ after row=\bottomrule},
    	precision=0, fixed zerofill, column type={r},
    	every row/.style={
    		column type=r,
    		dec sep align,
    		fixed,
    		fixed zerofill,
    	},
    	columns={case, flat, mpc, pfsolve},
    	every odd row/.style={before row={\rowcolor{tablecolor1}}},
        every even row/.style={before row={\rowcolor{tablecolor2}}},
    	columns/case/.style={string type, column type={l}, column name={Benchmark}},
    	columns/flat/.style={column name={Flat}},
    	columns/mpc/.style={column name={MPC}},
    	columns/pfsolve/.style={column name={PF}},
    	empty cells with={---} 
    	]\OPFtableitersMIPSsc
    \end{table}
\end{minipage} 

\begin{minipage}{0.45\textwidth}    
    \begin{table}[H]
        \footnotesize
    	\centering
    	\caption{Overall time (s) (\MIPS{}-\PARDISO{}). \label{tab:OPFstartMIPSscPardiso-t}}
    	\pgfplotstabletypeset[
    	every head row/.style={ before row={\toprule},after row=\midrule},
    	every last row/.style={ after row=\bottomrule},
    	precision=2, fixed zerofill, column type={r},
    	every row/.style={
    		column type=r,
    		dec sep align,
    		fixed,
    		fixed zerofill,
    	},
    	columns={case, flat, mpc, pfsolve},
    	every odd row/.style={before row={\rowcolor{tablecolor1}}},
        every even row/.style={before row={\rowcolor{tablecolor2}}},
    	columns/case/.style={string type, column type={l}, column name={Benchmark}},
    	columns/flat/.style={column name={Flat}},
    	columns/mpc/.style={column name={MPC}},
    	columns/pfsolve/.style={column name={PF}},
    	empty cells with={---} 
    	]\OPFtabletimeMIPSscPardiso
    \end{table}
\end{minipage}
\qquad
\begin{minipage}{0.45\textwidth}    
    \begin{table}[H]
        \footnotesize
    	\centering
    	\caption{Number of iterations (\MIPS{}-\PARDISO{}). \label{tab:OPFstartMIPSscPardiso-i}}
    	\pgfplotstabletypeset[
    	every head row/.style={ before row={\toprule},after row=\midrule},
    	every last row/.style={ after row=\bottomrule},
    	precision=0, fixed zerofill, column type={r},
    	every row/.style={
    		column type=r,
    		dec sep align,
    		fixed,
    		fixed zerofill,
    	},
    	columns={case, flat, mpc, pfsolve},
    	every odd row/.style={before row={\rowcolor{tablecolor1}}},
        every even row/.style={before row={\rowcolor{tablecolor2}}},
    	columns/case/.style={string type, column type={l}, column name={Benchmark}},
    	columns/flat/.style={column name={Flat}},
    	columns/mpc/.style={column name={MPC}},
    	columns/pfsolve/.style={column name={PF}},
    	empty cells with={---} 
    	]\OPFtableitersMIPSscPardiso
    \end{table}
\end{minipage} 


\begin{minipage}{0.45\textwidth}
    \begin{table}[H]
        \footnotesize
    	\centering
    	\caption{Overall time (s) (\FMINCON{} 2018b). \label{tab:OPFstartFMINCONa-t}}
    	\pgfplotstabletypeset[
    	every head row/.style={ before row={\toprule},after row=\midrule},
    	every last row/.style={ after row=\bottomrule},
    	precision=2, fixed zerofill, column type={r},
    	every row/.style={
    		column type=r,
    		dec sep align,
    		fixed,
    		fixed zerofill,
    	},
    	columns={case, flat, mpc, pfsolve},
    	every odd row/.style={before row={\rowcolor{tablecolor1}}},
        every even row/.style={before row={\rowcolor{tablecolor2}}},
    	columns/case/.style={string type, column type={l}, column name={Benchmark}},
    	columns/flat/.style={column name={Flat}},
    	columns/mpc/.style={column name={MPC}},
    	columns/pfsolve/.style={column name={PF}},
    	empty cells with={---} 
    	]\OPFtabletimeFMINCONc
    \end{table}
\end{minipage}
\qquad
\begin{minipage}{0.45\textwidth}
    \begin{table}[H]
        \footnotesize
    	\centering
    	\caption{Number of iterations (\FMINCON{} 2018b). \label{tab:OPFstartFMINCONa-i}}
    	\pgfplotstabletypeset[
    	every head row/.style={ before row={\toprule},after row=\midrule},
    	every last row/.style={ after row=\bottomrule},
    	precision=0, fixed zerofill, column type={r},
    	every row/.style={
    		column type=r,
    		dec sep align,
    		fixed,
    		fixed zerofill,
    	},
    	columns={case, flat, mpc, pfsolve},
    	every odd row/.style={before row={\rowcolor{tablecolor1}}},
        every even row/.style={before row={\rowcolor{tablecolor2}}},
    	columns/case/.style={string type, column type={l}, column name={Benchmark}},
    	columns/flat/.style={column name={Flat}},
    	columns/mpc/.style={column name={MPC}},
    	columns/pfsolve/.style={column name={PF}},
    	empty cells with={---} 
    	]\OPFtableitersFMINCONc
    \end{table}
\end{minipage}

\begin{minipage}{0.45\textwidth}
    \begin{table}[H]
        \footnotesize
    	\centering
    	\caption{Overall time (s) (\FMINCON{} 2017b). \label{tab:OPFstartFMINCONb-t}}
    	\pgfplotstabletypeset[
    	every head row/.style={ before row={\toprule},after row=\midrule},
    	every last row/.style={ after row=\bottomrule},
    	precision=2, fixed zerofill, column type={r},
    	every row/.style={
    		column type=r,
    		dec sep align,
    		fixed,
    		fixed zerofill,
    	},
    	columns={case, flat, mpc, pfsolve},
    	every odd row/.style={before row={\rowcolor{tablecolor1}}},
        every even row/.style={before row={\rowcolor{tablecolor2}}},
    	columns/case/.style={string type, column type={l}, column name={Benchmark}},
    	columns/flat/.style={column name={Flat}},
    	columns/mpc/.style={column name={MPC}},
    	columns/pfsolve/.style={column name={PF}},
    	empty cells with={---} 
    	]\OPFtabletimeFMINCONb
    \end{table}
\end{minipage}
\qquad
\begin{minipage}{0.45\textwidth}
    \begin{table}[H]
        \footnotesize
    	\centering
    	\caption{Number of iterations (\FMINCON{} 2017b). \label{tab:OPFstartFMINCONb-i}}
    	\pgfplotstabletypeset[
    	every head row/.style={ before row={\toprule},after row=\midrule},
    	every last row/.style={ after row=\bottomrule},
    	precision=0, fixed zerofill, column type={r},
    	every row/.style={
    		column type=r,
    		dec sep align,
    		fixed,
    		fixed zerofill,
    	},
    	columns={case, flat, mpc, pfsolve},
    	every odd row/.style={before row={\rowcolor{tablecolor1}}},
        every even row/.style={before row={\rowcolor{tablecolor2}}},
    	columns/case/.style={string type, column type={l}, column name={Benchmark}},
    	columns/flat/.style={column name={Flat}},
    	columns/mpc/.style={column name={MPC}},
    	columns/pfsolve/.style={column name={PF}},
    	empty cells with={---} 
    	]\OPFtableitersFMINCONb
    \end{table}
\end{minipage}

\begin{minipage}{0.45\textwidth}
    \begin{table}[H]
        \footnotesize
    	\centering
    	\caption{Overall time (s) (\KNITRO{} 11). \label{tab:OPFstartKNITROb-t}}
    	\pgfplotstabletypeset[
    	every head row/.style={ before row={\toprule},after row=\midrule},
    	every last row/.style={ after row=\bottomrule},
    	precision=2, fixed zerofill, column type={r},
    	every row/.style={
    		column type=r,
    		dec sep align,
    		fixed,
    		fixed zerofill,
    	},
    	columns={case, flat, mpc, pfsolve},
    	every odd row/.style={before row={\rowcolor{tablecolor1}}},
        every even row/.style={before row={\rowcolor{tablecolor2}}},
    	columns/case/.style={string type, column type={l}, column name={Benchmark}},
    	columns/flat/.style={column name={Flat}},
    	columns/mpc/.style={column name={MPC}},
    	columns/pfsolve/.style={column name={PF}},
    	empty cells with={---} 
    	]\OPFtabletimeKNITROb
    \end{table}
\end{minipage}
\qquad
\begin{minipage}{0.45\textwidth}
    \begin{table}[H]
        \footnotesize
    	\centering
    	\caption{Number of iterations (\KNITRO{} 11). \label{tab:OPFstartKNITROb-i}}
    	\pgfplotstabletypeset[
    	every head row/.style={ before row={\toprule},after row=\midrule},
    	every last row/.style={ after row=\bottomrule},
    	precision=0, fixed zerofill, column type={r},
    	every row/.style={
    		column type=r,
    		dec sep align,
    		fixed,
    		fixed zerofill,
    	},
    	columns={case, flat, mpc, pfsolve},
    	every odd row/.style={before row={\rowcolor{tablecolor1}}},
        every even row/.style={before row={\rowcolor{tablecolor2}}},
    	columns/case/.style={string type, column type={l}, column name={Benchmark}},
    	columns/flat/.style={column name={Flat}},
    	columns/mpc/.style={column name={MPC}},
    	columns/pfsolve/.style={column name={PF}},
    	empty cells with={---} 
    	]\OPFtableitersKNITROb
    \end{table}
\end{minipage}

\subsection{OPF variants}
The  bus voltages in the standard  AC  OPF  problem  can be represented either in polar, or Cartesian coordinates. Another  formulation  of  the  AC  OPF  problem  uses current  balance  constraints  in  place  of  the  power  balance. Different representations of the complex voltage variables and formulation of the nodal balance equations lead to a different number of constraints and sparsity structure of the problem, which, in turn, influences the numerical behavior of the optimizer. The corresponding \matpower{} options are \texttt{opf.v\_cartesian}, specifying whether to use polar or Cartesian voltage coordinates, and option \texttt{opf. current\_balance}, which selects using either a current or power balance formulation for AC OPF.

The results presented in Table~\ref{tab:OPFvariants} suggest that the choice of formulation can significantly influence whether the benchmark case can be successfully solved. \MIPS{}, and \IPOPT{}-\PARDISO{} optimizers are more robust with polar voltage coordinates and nodal power balance  while the opposite is true for \FMINCON{} which is more robust with Cartesian voltage coordinates. Less robust optimizers are also not able to solve the large-scale cases due to the extensive time requirements of the linear solver or insufficient precision of the solution. 

\pgfplotstableread{
	Optimizer polarPower cartPower polarCurrent cartCurrent Total
	BELTISTOS-OPF          25 25 25 25 100/100 
	KNITRO12             25 25 25 24 99/100
	KNITRO11             24 25 25 25 99/100
	IPOPT-PARDISO      24 24 22 20 90/100
	FMINCON            20 24 20 24 88/100
	IPOPT-MA57         21 21 23 21 86/100
	MIPS-PARDISO       24 23 21 15 83/100
	MIPS-MATLAB'\textbackslash'       22 19 21 17 79/100
}\variantsSuccess

\begin{table}[ht!]
	\centering
	\caption{Number of solved benchmarks for different OPF formulations. \label{tab:OPFvariants}}
	\pgfplotstabletypeset[
	every head row/.style={ before row=\toprule,after row=\midrule},
	every last row/.style={ after row=\bottomrule},
	precision=0, fixed zerofill, column type={c},
	columns/Optimizer/.style={string type,column type=l},
	every odd row/.style={before row={\rowcolor{tablecolor1}}},
	every even row/.style={before row={\rowcolor{tablecolor2}}},
	columns/polarPower/.style={column name={Polar Power}},
	columns/polarCurrent/.style={column name={Polar Current}},
	columns/cartPower/.style={column name={Cartesian Power}},
	columns/cartCurrent/.style={column name={Cartesian Current}},
	columns/Total/.style={column name={Total}, string type,column type=r}
	]\variantsSuccess
\end{table}
\begin{figure*}[ht!]
	\centering
	\begin{subfigure}[b]{0.49\columnwidth}
		\begin{tikzpicture}

\pgfplotstableread{
ALPHA         SOLVER1         SOLVER2         SOLVER3         SOLVER4
1.0         0.32000         0.00000         0.60000         0.08000
1.1         0.56000         0.12000         0.72000         0.20000
1.2000000000000002         0.80000         0.28000         0.84000         0.24000
1.3000000000000003         0.84000         0.48000         0.92000         0.40000
1.4000000000000004         0.84000         0.60000         1.00000         0.64000
1.5000000000000004         0.96000         0.68000         1.00000         0.68000
1.6000000000000005         0.96000         0.72000         1.00000         0.76000
1.7000000000000006         1.00000         0.88000         1.00000         0.76000
1.8000000000000007         1.00000         0.92000         1.00000         0.80000
1.9000000000000008         1.00000         0.96000         1.00000         0.88000
2.000000000000001         1.00000         0.96000         1.00000         0.88000
2.100000000000001         1.00000         0.96000         1.00000         0.92000
2.200000000000001         1.00000         0.96000         1.00000         0.92000
2.300000000000001         1.00000         1.00000         1.00000         0.92000
2.4000000000000012         1.00000         1.00000         1.00000         0.92000
2.5000000000000013         1.00000         1.00000         1.00000         0.92000
2.6000000000000014         1.00000         1.00000         1.00000         0.92000
2.7000000000000015         1.00000         1.00000         1.00000         0.92000
2.8000000000000016         1.00000         1.00000         1.00000         0.92000
2.9000000000000017         1.00000         1.00000         1.00000         0.92000
3.0000000000000018         1.00000         1.00000         1.00000         0.92000
3.100000000000002         1.00000         1.00000         1.00000         0.92000
3.200000000000002         1.00000         1.00000         1.00000         0.96000
3.300000000000002         1.00000         1.00000         1.00000         0.96000
3.400000000000002         1.00000         1.00000         1.00000         0.96000
3.500000000000002         1.00000         1.00000         1.00000         0.96000
3.6000000000000023         1.00000         1.00000         1.00000         0.96000
3.7000000000000024         1.00000         1.00000         1.00000         0.96000
3.8000000000000025         1.00000         1.00000         1.00000         0.96000
3.9000000000000026         1.00000         1.00000         1.00000         1.00000
4.000000000000003         1.00000         1.00000         1.00000         1.00000
4.100000000000003         1.00000         1.00000         1.00000         1.00000
4.200000000000003         1.00000         1.00000         1.00000         1.00000
4.3000000000000025         1.00000         1.00000         1.00000         1.00000
4.400000000000003         1.00000         1.00000         1.00000         1.00000
4.5000000000000036         1.00000         1.00000         1.00000         1.00000
4.600000000000003         1.00000         1.00000         1.00000         1.00000
4.700000000000003         1.00000         1.00000         1.00000         1.00000
4.800000000000003         1.00000         1.00000         1.00000         1.00000
4.900000000000004         1.00000         1.00000         1.00000         1.00000
5.0000000000000036         1.00000         1.00000         1.00000         1.00000
}\datatable

  \begin{axis}[name=symb,
  	width=\columnwidth,
  	height=4cm,
	ymin=-0.02,
	ymax=1.02,
    enlarge x limits=0.05,
	axis lines*=left, 
	ymajorgrids, yminorgrids,
    xmajorgrids, xminorgrids,
	xticklabel style={rotate=0, xshift=-0.0cm, anchor=north, font=\scriptsize},
	ytick={0,  0.2,  0.4, 0.6, 0.8, 1},
    xmax=5,
	yticklabel style={font=\scriptsize},
	ylabel style={font=\scriptsize},
	xlabel style={font=\scriptsize},
	ylabel={$p_m({\alpha})$},
	xlabel={$\alpha$},
	legend style={at={(1.0,0.7)},legend cell align=left,align=right,draw=white!85!black,font=\scriptsize,legend columns=1}
	]
	
    \addplot[color=mycolor1, ultra thick] table[x=ALPHA, y=SOLVER1] {\datatable}; 
    \addlegendentry{Polar-Power}

    \addplot[color=mycolor2, dashdotted, ultra thick] table[x=ALPHA, y=SOLVER2] {\datatable}; 
    \addlegendentry{Rect-Power}
    
    \addplot[color=mycolor3, dashed, ultra thick] table[x=ALPHA, y=SOLVER3] {\datatable}; 
    \addlegendentry{Polar-Current}

	\addplot[color=mycolor4, dotted, ultra thick] table[x=ALPHA, y=SOLVER4] {\datatable}; 
	\addlegendentry{Rect-Current}

  \end{axis}
  
\end{tikzpicture}
		\caption{Overall time ($\alpha_{max}=4$)\label{fig:profileTvariantsBELTISTOS} }
	\end{subfigure}
	\begin{subfigure}[b]{0.5\columnwidth}
		\begin{tikzpicture}

\pgfplotstableread{
ALPHA         SOLVER1         SOLVER2         SOLVER3         SOLVER4
1.0         0.32000         0.20000         0.32000         0.20000
1.05         0.40000         0.32000         0.36000         0.24000
1.1         0.68000         0.48000         0.60000         0.36000
1.1500000000000001         0.68000         0.56000         0.60000         0.44000
1.2000000000000002         0.84000         0.64000         0.72000         0.52000
1.2500000000000002         0.84000         0.64000         0.80000         0.60000
1.3000000000000003         0.84000         0.68000         0.84000         0.72000
1.3500000000000003         0.92000         0.68000         0.88000         0.72000
1.4000000000000004         0.96000         0.76000         0.88000         0.72000
1.4500000000000004         0.96000         0.80000         0.92000         0.80000
1.5000000000000004         0.96000         0.80000         0.96000         0.80000
1.5500000000000005         0.96000         0.88000         0.96000         0.80000
1.6000000000000005         0.96000         0.88000         1.00000         0.80000
1.6500000000000006         0.96000         0.88000         1.00000         0.80000
1.7000000000000006         0.96000         0.88000         1.00000         0.84000
1.7500000000000007         0.96000         1.00000         1.00000         0.88000
1.8000000000000007         1.00000         1.00000         1.00000         0.92000
1.8500000000000008         1.00000         1.00000         1.00000         0.92000
1.9000000000000008         1.00000         1.00000         1.00000         0.92000
1.9500000000000008         1.00000         1.00000         1.00000         0.92000
2.000000000000001         1.00000         1.00000         1.00000         0.92000
2.0500000000000007         1.00000         1.00000         1.00000         0.92000
2.100000000000001         1.00000         1.00000         1.00000         0.92000
2.1500000000000012         1.00000         1.00000         1.00000         0.92000
2.200000000000001         1.00000         1.00000         1.00000         0.92000
2.250000000000001         1.00000         1.00000         1.00000         0.92000
2.300000000000001         1.00000         1.00000         1.00000         0.96000
2.3500000000000014         1.00000         1.00000         1.00000         0.96000
2.4000000000000012         1.00000         1.00000         1.00000         0.96000
2.450000000000001         1.00000         1.00000         1.00000         0.96000
2.5000000000000013         1.00000         1.00000         1.00000         0.96000
2.5500000000000016         1.00000         1.00000         1.00000         0.96000
2.6000000000000014         1.00000         1.00000         1.00000         0.96000
2.6500000000000012         1.00000         1.00000         1.00000         0.96000
2.7000000000000015         1.00000         1.00000         1.00000         0.96000
2.7500000000000018         1.00000         1.00000         1.00000         1.00000
2.8000000000000016         1.00000         1.00000         1.00000         1.00000
2.8500000000000014         1.00000         1.00000         1.00000         1.00000
2.9000000000000017         1.00000         1.00000         1.00000         1.00000
2.950000000000002         1.00000         1.00000         1.00000         1.00000
3.0000000000000018         1.00000         1.00000         1.00000         1.00000
}\datatable

  \begin{axis}[name=symb,
  	width=\columnwidth,
  	height=4cm,
	ymin=-0.02,
	ymax=1.02,
    enlarge x limits=0.05,
	axis lines*=left, 
	ymajorgrids, yminorgrids,
    xmajorgrids, xminorgrids,
	xticklabel style={rotate=0, xshift=-0.0cm, anchor=north, font=\scriptsize},
	ytick={0,  0.2,  0.4, 0.6, 0.8, 1},
    xmax=3,
	yticklabel style={font=\scriptsize},
	ylabel style={font=\scriptsize},
	xlabel style={font=\scriptsize},
	ylabel={$p_m({\alpha})$},
	xlabel={$\alpha$},
	legend style={at={(1.0,0.7)},legend cell align=left,align=right,draw=white!85!black,font=\scriptsize,legend columns=1}
	]
	
    \addplot[color=mycolor1, ultra thick] table[x=ALPHA, y=SOLVER1] {\datatable}; 
    \addlegendentry{Polar-Power}

    \addplot[color=mycolor2, dashdotted, ultra thick] table[x=ALPHA, y=SOLVER2] {\datatable}; 
    \addlegendentry{Rect-Power}
    
    \addplot[color=mycolor3, dashed, ultra thick] table[x=ALPHA, y=SOLVER3] {\datatable}; 
    \addlegendentry{Polar-Current}

	\addplot[color=mycolor4, dotted, ultra thick] table[x=ALPHA, y=SOLVER4] {\datatable}; 
	\addlegendentry{Rect-Current}

  \end{axis}
  
\end{tikzpicture}
		\caption{ Number of iterations ($\alpha_{max}=3$)\label{fig:profileIvariantsBELTISTOS} }
	\end{subfigure}
	\caption{\BELTISTOSOPF{} performance profiles for various OPF formulations.\label{fig:profileVariantsBELTISTOS}}
\end{figure*}

\begin{figure*}[ht!]
	\centering
	\begin{subfigure}[b]{0.49\columnwidth}
		\begin{tikzpicture}

\pgfplotstableread{
ALPHA         SOLVER1         SOLVER2         SOLVER3         SOLVER4
1.0         0.32000         0.08000         0.60000         0.00000
1.2         0.68000         0.68000         0.76000         0.40000
1.4         0.72000         0.96000         0.88000         0.64000
1.5999999999999999         0.72000         0.96000         0.96000         0.72000
1.7999999999999998         0.80000         1.00000         0.96000         0.88000
1.9999999999999998         0.80000         1.00000         1.00000         0.88000
2.1999999999999997         0.84000         1.00000         1.00000         0.88000
2.3999999999999995         0.84000         1.00000         1.00000         0.88000
2.5999999999999996         0.84000         1.00000         1.00000         0.92000
2.8         0.84000         1.00000         1.00000         0.92000
2.9999999999999996         0.84000         1.00000         1.00000         0.92000
3.1999999999999993         0.84000         1.00000         1.00000         0.92000
3.3999999999999995         0.88000         1.00000         1.00000         0.96000
3.5999999999999996         0.92000         1.00000         1.00000         0.96000
3.7999999999999994         0.92000         1.00000         1.00000         0.96000
3.999999999999999         0.92000         1.00000         1.00000         0.96000
4.199999999999999         0.92000         1.00000         1.00000         0.96000
4.3999999999999995         0.92000         1.00000         1.00000         0.96000
4.6         0.92000         1.00000         1.00000         0.96000
4.799999999999999         0.92000         1.00000         1.00000         0.96000
4.999999999999999         0.92000         1.00000         1.00000         0.96000
5.199999999999999         0.92000         1.00000         1.00000         0.96000
5.399999999999999         0.92000         1.00000         1.00000         0.96000
5.599999999999999         0.92000         1.00000         1.00000         0.96000
5.799999999999999         0.92000         1.00000         1.00000         0.96000
5.999999999999999         0.92000         1.00000         1.00000         0.96000
6.199999999999999         0.92000         1.00000         1.00000         1.00000
6.399999999999999         0.92000         1.00000         1.00000         1.00000
6.599999999999999         0.92000         1.00000         1.00000         1.00000
6.799999999999999         0.92000         1.00000         1.00000         1.00000
6.999999999999998         0.92000         1.00000         1.00000         1.00000
7.199999999999998         0.96000         1.00000         1.00000         1.00000
7.399999999999999         0.96000         1.00000         1.00000         1.00000
7.599999999999999         0.96000         1.00000         1.00000         1.00000
7.799999999999999         0.96000         1.00000         1.00000         1.00000
7.999999999999998         0.96000         1.00000         1.00000         1.00000
}\datatable

  \begin{axis}[name=symb,
  	width=\columnwidth,
  	height=4cm,
	ymin=-0.02,
	ymax=1.02,
    enlarge x limits=0.05,
	axis lines*=left, 
	ymajorgrids, yminorgrids,
    xmajorgrids, xminorgrids,
	xticklabel style={rotate=0, xshift=-0.0cm, anchor=north, font=\scriptsize},
	ytick={0,  0.2,  0.4, 0.6, 0.8, 1},
    xmax=8,
	yticklabel style={font=\scriptsize},
	ylabel style={font=\scriptsize},
	xlabel style={font=\scriptsize},
	ylabel={$p_m({\alpha})$},
	xlabel={$\alpha$},
	legend style={at={(1.0,0.70)},legend cell align=left,align=right,draw=white!85!black,font=\scriptsize,legend columns=1}
	]
	
    \addplot[color=mycolor1, ultra thick] table[x=ALPHA, y=SOLVER1] {\datatable}; 
    \addlegendentry{Polar-Power}

    \addplot[color=mycolor2, dashdotted, ultra thick] table[x=ALPHA, y=SOLVER2] {\datatable}; 
    \addlegendentry{Rect-Power}
    
    \addplot[color=mycolor3, dashed, ultra thick] table[x=ALPHA, y=SOLVER3] {\datatable}; 
    \addlegendentry{Polar-Current}

	\addplot[color=mycolor4, dotted, ultra thick] table[x=ALPHA, y=SOLVER4] {\datatable}; 
	\addlegendentry{Rect-Current}

  \end{axis}
  
\end{tikzpicture}
		\caption{Overall time ($\alpha_{max}=8$)\label{fig:profileTvariantsKNITRO} }
	\end{subfigure}
	\begin{subfigure}[b]{0.5\columnwidth}
		\begin{tikzpicture}

\pgfplotstableread{
ALPHA         SOLVER1         SOLVER2         SOLVER3         SOLVER4
1.0         0.48000         0.32000         0.24000         0.24000
1.2         0.64000         0.96000         0.72000         0.48000
1.4         0.64000         1.00000         0.76000         0.68000
1.5999999999999999         0.72000         1.00000         0.88000         0.84000
1.7999999999999998         0.72000         1.00000         0.92000         0.84000
1.9999999999999998         0.72000         1.00000         0.92000         0.96000
2.1999999999999997         0.76000         1.00000         1.00000         0.96000
2.3999999999999995         0.76000         1.00000         1.00000         0.96000
2.5999999999999996         0.84000         1.00000         1.00000         0.96000
2.8         0.84000         1.00000         1.00000         0.96000
2.9999999999999996         0.84000         1.00000         1.00000         0.96000
3.1999999999999993         0.84000         1.00000         1.00000         0.96000
3.3999999999999995         0.84000         1.00000         1.00000         1.00000
3.5999999999999996         0.84000         1.00000         1.00000         1.00000
3.7999999999999994         0.84000         1.00000         1.00000         1.00000
3.999999999999999         0.84000         1.00000         1.00000         1.00000
4.199999999999999         0.88000         1.00000         1.00000         1.00000
4.3999999999999995         0.88000         1.00000         1.00000         1.00000
4.6         0.88000         1.00000         1.00000         1.00000
4.799999999999999         0.88000         1.00000         1.00000         1.00000
4.999999999999999         0.92000         1.00000         1.00000         1.00000
5.199999999999999         0.92000         1.00000         1.00000         1.00000
5.399999999999999         0.92000         1.00000         1.00000         1.00000
5.599999999999999         0.92000         1.00000         1.00000         1.00000
5.799999999999999         0.92000         1.00000         1.00000         1.00000
5.999999999999999         0.92000         1.00000         1.00000         1.00000
6.199999999999999         0.92000         1.00000         1.00000         1.00000
6.399999999999999         0.92000         1.00000         1.00000         1.00000
6.599999999999999         0.92000         1.00000         1.00000         1.00000
6.799999999999999         0.92000         1.00000         1.00000         1.00000
6.999999999999998         0.92000         1.00000         1.00000         1.00000
7.199999999999998         0.92000         1.00000         1.00000         1.00000
7.399999999999999         0.92000         1.00000         1.00000         1.00000
7.599999999999999         0.92000         1.00000         1.00000         1.00000
7.799999999999999         0.92000         1.00000         1.00000         1.00000
7.999999999999998         0.92000         1.00000         1.00000         1.00000
8.2         0.92000         1.00000         1.00000         1.00000
8.399999999999999         0.92000         1.00000         1.00000         1.00000
8.599999999999998         0.96000         1.00000         1.00000         1.00000
8.799999999999997         0.96000         1.00000         1.00000         1.00000
8.999999999999998         0.96000         1.00000         1.00000         1.00000
9.199999999999998         0.96000         1.00000         1.00000         1.00000
9.399999999999999         0.96000         1.00000         1.00000         1.00000
9.599999999999998         0.96000         1.00000         1.00000         1.00000
9.799999999999997         0.96000         1.00000         1.00000         1.00000
9.999999999999998         0.96000         1.00000         1.00000         1.00000
}\datatable

  \begin{axis}[name=symb,
  	width=\columnwidth,
  	height=4cm,
	ymin=-0.02,
	ymax=1.02,
    enlarge x limits=0.05,
	axis lines*=left, 
	ymajorgrids, yminorgrids,
    xmajorgrids, xminorgrids,
	xticklabel style={rotate=0, xshift=-0.0cm, anchor=north, font=\scriptsize},
	ytick={0,  0.2,  0.4, 0.6, 0.8, 1},
    xmax=10,
	yticklabel style={font=\scriptsize},
	ylabel style={font=\scriptsize},
	xlabel style={font=\scriptsize},
	ylabel={$p_m({\alpha})$},
	xlabel={$\alpha$},
	legend style={at={(1.0,0.70)},legend cell align=left,align=right,draw=white!85!black,font=\scriptsize,legend columns=1}
	]
	
    \addplot[color=mycolor1, ultra thick] table[x=ALPHA, y=SOLVER1] {\datatable}; 
    \addlegendentry{Polar-Power}

    \addplot[color=mycolor2, dashdotted, ultra thick] table[x=ALPHA, y=SOLVER2] {\datatable}; 
    \addlegendentry{Rect-Power}
    
    \addplot[color=mycolor3, dashed, ultra thick] table[x=ALPHA, y=SOLVER3] {\datatable}; 
    \addlegendentry{Polar-Current}

	\addplot[color=mycolor4, dotted, ultra thick] table[x=ALPHA, y=SOLVER4] {\datatable}; 
	\addlegendentry{Rect-Current}

  \end{axis}
  
\end{tikzpicture}
		\caption{ Number of iterations ($\alpha_{max}=9$)\label{fig:profileIvariantsKNITRO} }
	\end{subfigure}
	\caption{\KNITRO{}11 performance profiles for various OPF formulations.\label{fig:profileVariantsKNITRO}}
\end{figure*}

\begin{figure*}[ht!]
	\centering
	\begin{subfigure}[b]{0.49\columnwidth}
		\begin{tikzpicture}

\pgfplotstableread{
ALPHA         SOLVER1         SOLVER2         SOLVER3         SOLVER4
1.0         0.56000         0.12000         0.32000         0.00000
1.2         0.84000         0.72000         0.84000         0.32000
1.4         0.92000         0.92000         0.88000         0.60000
1.5999999999999999         0.92000         0.96000         0.88000         0.72000
1.7999999999999998         1.00000         0.96000         0.96000         0.88000
1.9999999999999998         1.00000         0.96000         0.96000         0.92000
2.1999999999999997         1.00000         0.96000         0.96000         0.96000
2.3999999999999995         1.00000         0.96000         0.96000         0.96000
2.5999999999999996         1.00000         0.96000         0.96000         0.96000
2.8         1.00000         0.96000         0.96000         0.96000
2.9999999999999996         1.00000         0.96000         0.96000         0.96000
3.1999999999999993         1.00000         0.96000         0.96000         0.96000
3.3999999999999995         1.00000         0.96000         0.96000         0.96000
3.5999999999999996         1.00000         0.96000         0.96000         0.96000
3.7999999999999994         1.00000         0.96000         0.96000         0.96000
3.999999999999999         1.00000         0.96000         0.96000         0.96000
4.199999999999999         1.00000         0.96000         0.96000         0.96000
4.3999999999999995         1.00000         0.96000         0.96000         0.96000
4.6         1.00000         0.96000         0.96000         0.96000
4.799999999999999         1.00000         0.96000         0.96000         0.96000
4.999999999999999         1.00000         0.96000         0.96000         0.96000
5.199999999999999         1.00000         0.96000         0.96000         0.96000
5.399999999999999         1.00000         0.96000         0.96000         0.96000
5.599999999999999         1.00000         0.96000         0.96000         0.96000
5.799999999999999         1.00000         0.96000         0.96000         0.96000
5.999999999999999         1.00000         0.96000         0.96000         0.96000
6.199999999999999         1.00000         0.96000         0.96000         0.96000
6.399999999999999         1.00000         0.96000         0.96000         0.96000
6.599999999999999         1.00000         0.96000         0.96000         0.96000
6.799999999999999         1.00000         0.96000         0.96000         0.96000
6.999999999999998         1.00000         0.96000         0.96000         0.96000
7.199999999999998         1.00000         0.96000         0.96000         0.96000
7.399999999999999         1.00000         0.96000         0.96000         0.96000
7.599999999999999         1.00000         0.96000         0.96000         0.96000
7.799999999999999         1.00000         0.96000         0.96000         0.96000
7.999999999999998         1.00000         0.96000         0.96000         0.96000
8.2         1.00000         0.96000         0.96000         0.96000
8.399999999999999         1.00000         0.96000         0.96000         0.96000
8.599999999999998         1.00000         0.96000         0.96000         0.96000
8.799999999999997         1.00000         0.96000         0.96000         0.96000
8.999999999999998         1.00000         0.96000         1.00000         0.96000
9.999999999999998         1.00000         0.96000         1.00000         0.96000
}\datatable

  \begin{axis}[name=symb,
  	width=\columnwidth,
  	height=4cm,
	ymin=-0.02,
	ymax=1.02,
    enlarge x limits=0.05,
	axis lines*=left, 
	ymajorgrids, yminorgrids,
    xmajorgrids, xminorgrids,
	xticklabel style={rotate=0, xshift=-0.0cm, anchor=north, font=\scriptsize},
	ytick={0,  0.2,  0.4, 0.6, 0.8, 1},
    xmax=9.5,
	yticklabel style={font=\scriptsize},
	ylabel style={font=\scriptsize},
	xlabel style={font=\scriptsize},
	ylabel={$p_m({\alpha})$},
	xlabel={$\alpha$},
	legend style={at={(1.0,0.70)},legend cell align=left,align=right,draw=white!85!black,font=\scriptsize,legend columns=1}
	]
	
    \addplot[color=mycolor1, ultra thick] table[x=ALPHA, y=SOLVER1] {\datatable}; 
    \addlegendentry{Polar-Power}

    \addplot[color=mycolor2, dashdotted, ultra thick] table[x=ALPHA, y=SOLVER2] {\datatable}; 
    \addlegendentry{Rect-Power}
    
    \addplot[color=mycolor3, dashed, ultra thick] table[x=ALPHA, y=SOLVER3] {\datatable}; 
    \addlegendentry{Polar-Current}

	\addplot[color=mycolor4, dotted, ultra thick] table[x=ALPHA, y=SOLVER4] {\datatable}; 
	\addlegendentry{Rect-Current}

  \end{axis}
  
\end{tikzpicture}
		\caption{Overall time ($\alpha_{max}=9$)\label{fig:profileTvariantsKNITRO12} }
	\end{subfigure}
	\begin{subfigure}[b]{0.5\columnwidth}
		\begin{tikzpicture}

\pgfplotstableread{
ALPHA         SOLVER1         SOLVER2         SOLVER3         SOLVER4
1.0         0.80000         0.04000         0.44000         0.00000
1.2         0.96000         0.80000         0.76000         0.56000
1.4         0.96000         0.96000         0.76000         0.68000
1.5999999999999999         0.96000         0.96000         0.88000         0.80000
1.7999999999999998         0.96000         0.96000         0.92000         0.84000
1.9999999999999998         0.96000         0.96000         0.92000         0.88000
2.1999999999999997         1.00000         0.96000         0.96000         0.96000
2.3999999999999995         1.00000         0.96000         0.96000         0.96000
2.5999999999999996         1.00000         0.96000         0.96000         0.96000
2.8         1.00000         0.96000         0.96000         0.96000
2.9999999999999996         1.00000         0.96000         0.96000         0.96000
3.1999999999999993         1.00000         0.96000         0.96000         0.96000
3.3999999999999995         1.00000         0.96000         0.96000         0.96000
3.5999999999999996         1.00000         0.96000         0.96000         0.96000
3.7999999999999994         1.00000         0.96000         1.00000         0.96000
3.999999999999999         1.00000         0.96000         1.00000         0.96000
4.199999999999999         1.00000         0.96000         1.00000         0.96000
4.3999999999999995         1.00000         0.96000         1.00000         0.96000
4.6         1.00000         0.96000         1.00000         0.96000
4.799999999999999         1.00000         0.96000         1.00000         0.96000
4.999999999999999         1.00000         0.96000         1.00000         0.96000
}\datatable

  \begin{axis}[name=symb,
  	width=\columnwidth,
  	height=4cm,
	ymin=-0.02,
	ymax=1.02,
    enlarge x limits=0.05,
	axis lines*=left, 
	ymajorgrids, yminorgrids,
    xmajorgrids, xminorgrids,
	xticklabel style={rotate=0, xshift=-0.0cm, anchor=north, font=\scriptsize},
	ytick={0,  0.2,  0.4, 0.6, 0.8, 1},
    xmax=5,
	yticklabel style={font=\scriptsize},
	ylabel style={font=\scriptsize},
	xlabel style={font=\scriptsize},
	ylabel={$p_m({\alpha})$},
	xlabel={$\alpha$},
	legend style={at={(1.0,0.70)},legend cell align=left,align=right,draw=white!85!black,font=\scriptsize,legend columns=1}
	]
	
    \addplot[color=mycolor1, ultra thick] table[x=ALPHA, y=SOLVER1] {\datatable}; 
    \addlegendentry{Polar-Power}

    \addplot[color=mycolor2, dashdotted, ultra thick] table[x=ALPHA, y=SOLVER2] {\datatable}; 
    \addlegendentry{Rect-Power}
    
    \addplot[color=mycolor3, dashed, ultra thick] table[x=ALPHA, y=SOLVER3] {\datatable}; 
    \addlegendentry{Polar-Current}

	\addplot[color=mycolor4, dotted, ultra thick] table[x=ALPHA, y=SOLVER4] {\datatable}; 
	\addlegendentry{Rect-Current}

  \end{axis}
  
\end{tikzpicture}
		\caption{ Number of iterations ($\alpha_{max}=5$)\label{fig:profileIvariantsKNITRO12} }
	\end{subfigure}
	\caption{\KNITRO{}12 performance profiles for various OPF formulations.\label{fig:profileVariantsKNITRO12}}
\end{figure*}

\begin{figure*}[ht!]
	\centering
	\begin{subfigure}[b]{0.49\columnwidth}
		\input{profilesNEW/OPFvariantsFMINCON-all-t.tex}
		\caption{Overall time ($\alpha_{max}=27$)\label{fig:profileTvariantsFMINCON} }
	\end{subfigure}
	\begin{subfigure}[b]{0.5\columnwidth}
		\begin{tikzpicture}

\pgfplotstableread{
ALPHA         SOLVER1         SOLVER2         SOLVER3         SOLVER4
1.0         0.28000         0.52000         0.12000         0.08000
1.2         0.52000         0.68000         0.28000         0.52000
1.4         0.64000         0.80000         0.52000         0.72000
1.5999999999999999         0.68000         0.84000         0.60000         0.76000
1.7999999999999998         0.68000         0.88000         0.60000         0.80000
1.9999999999999998         0.72000         0.92000         0.60000         0.80000
2.1999999999999997         0.72000         0.92000         0.60000         0.92000
2.3999999999999995         0.72000         0.92000         0.64000         0.92000
2.5999999999999996         0.72000         0.92000         0.64000         0.92000
2.8         0.72000         0.92000         0.68000         0.92000
2.9999999999999996         0.72000         0.96000         0.68000         0.92000
3.1999999999999993         0.72000         0.96000         0.72000         0.92000
3.3999999999999995         0.72000         0.96000         0.72000         0.92000
3.5999999999999996         0.76000         0.96000         0.72000         0.92000
3.7999999999999994         0.76000         0.96000         0.72000         0.92000
3.999999999999999         0.76000         0.96000         0.72000         0.92000
4.199999999999999         0.76000         0.96000         0.72000         0.92000
4.3999999999999995         0.80000         0.96000         0.76000         0.92000
4.6         0.80000         0.96000         0.76000         0.92000
4.799999999999999         0.80000         0.96000         0.76000         0.92000
4.999999999999999         0.80000         0.96000         0.76000         0.92000
5.199999999999999         0.80000         0.96000         0.76000         0.92000
5.399999999999999         0.80000         0.96000         0.76000         0.92000
5.599999999999999         0.80000         0.96000         0.76000         0.92000
5.799999999999999         0.80000         0.96000         0.76000         0.92000
5.999999999999999         0.80000         0.96000         0.76000         0.92000
6.199999999999999         0.80000         0.96000         0.76000         0.92000
6.399999999999999         0.80000         0.96000         0.76000         0.92000
6.599999999999999         0.80000         0.96000         0.76000         0.92000
6.799999999999999         0.80000         0.96000         0.76000         0.92000
6.999999999999998         0.80000         0.96000         0.76000         0.92000
7.199999999999998         0.80000         0.96000         0.76000         0.92000
7.399999999999999         0.80000         0.96000         0.76000         0.92000
7.599999999999999         0.80000         0.96000         0.76000         0.96000
7.799999999999999         0.80000         0.96000         0.76000         0.96000
7.999999999999998         0.80000         0.96000         0.76000         0.96000
8.2         0.80000         0.96000         0.76000         0.96000
8.399999999999999         0.80000         0.96000         0.76000         0.96000
8.599999999999998         0.80000         0.96000         0.76000         0.96000
8.799999999999997         0.80000         0.96000         0.76000         0.96000
8.999999999999998         0.80000         0.96000         0.76000         0.96000
9.199999999999998         0.80000         0.96000         0.76000         0.96000
9.399999999999999         0.80000         0.96000         0.76000         0.96000
9.599999999999998         0.80000         0.96000         0.76000         0.96000
9.799999999999997         0.80000         0.96000         0.76000         0.96000
9.999999999999998         0.80000         0.96000         0.76000         0.96000
10.199999999999998         0.80000         0.96000         0.76000         0.96000
10.399999999999999         0.80000         0.96000         0.76000         0.96000
10.599999999999998         0.80000         0.96000         0.76000         0.96000
10.799999999999997         0.80000         0.96000         0.76000         0.96000
10.999999999999998         0.80000         0.96000         0.76000         0.96000
11.199999999999998         0.80000         0.96000         0.76000         0.96000
11.399999999999999         0.80000         0.96000         0.76000         0.96000
11.599999999999998         0.80000         0.96000         0.76000         0.96000
11.799999999999997         0.80000         0.96000         0.76000         0.96000
11.999999999999998         0.80000         0.96000         0.76000         0.96000
12.199999999999998         0.80000         0.96000         0.76000         0.96000
12.399999999999997         0.80000         0.96000         0.76000         0.96000
12.599999999999998         0.80000         0.96000         0.76000         0.96000
12.799999999999997         0.80000         0.96000         0.76000         0.96000
12.999999999999996         0.80000         0.96000         0.76000         0.96000
13.199999999999998         0.80000         0.96000         0.76000         0.96000
13.399999999999997         0.80000         0.96000         0.76000         0.96000
13.599999999999998         0.80000         0.96000         0.76000         0.96000
13.799999999999997         0.80000         0.96000         0.76000         0.96000
13.999999999999996         0.80000         0.96000         0.76000         0.96000
14.199999999999998         0.80000         0.96000         0.76000         0.96000
14.399999999999997         0.80000         0.96000         0.76000         0.96000
14.599999999999998         0.80000         0.96000         0.76000         0.96000
14.799999999999997         0.80000         0.96000         0.76000         0.96000
14.999999999999996         0.80000         0.96000         0.76000         0.96000
15.199999999999998         0.80000         0.96000         0.80000         0.96000
15.399999999999997         0.80000         0.96000         0.80000         0.96000
15.599999999999996         0.80000         0.96000         0.80000         0.96000
15.799999999999997         0.80000         0.96000         0.80000         0.96000
15.999999999999996         0.80000         0.96000         0.80000         0.96000
}\datatable

  \begin{axis}[name=symb,
  	width=\columnwidth,
  	height=4cm,
	ymin=-0.02,
	ymax=1.02,
    enlarge x limits=0.05,
	axis lines*=left, 
	ymajorgrids, yminorgrids,
    xmajorgrids, xminorgrids,
	xticklabel style={rotate=0, xshift=-0.0cm, anchor=north, font=\scriptsize},
	ytick={0,  0.2,  0.4, 0.6, 0.8, 1},
    xmax=16,
	yticklabel style={font=\scriptsize},
	ylabel style={font=\scriptsize},
	xlabel style={font=\scriptsize},
	ylabel={$p_m({\alpha})$},
	xlabel={$\alpha$},
	legend style={at={(1.0,0.7)},legend cell align=left,align=right,draw=white!85!black,font=\scriptsize,legend columns=1}
	]
	
    \addplot[color=mycolor1, ultra thick] table[x=ALPHA, y=SOLVER1] {\datatable}; 
    \addlegendentry{Polar-Power}

    \addplot[color=mycolor2, dashdotted, ultra thick] table[x=ALPHA, y=SOLVER2] {\datatable}; 
    \addlegendentry{Rect-Power}
    
    \addplot[color=mycolor3, dashed, ultra thick] table[x=ALPHA, y=SOLVER3] {\datatable}; 
    \addlegendentry{Polar-Current}

	\addplot[color=mycolor4, dotted, ultra thick] table[x=ALPHA, y=SOLVER4] {\datatable}; 
	\addlegendentry{Rect-Current}

  \end{axis}
  
\end{tikzpicture}
		\caption{ Number of iterations ($\alpha_{max}=16$)\label{fig:profileIvariantsFMINCON} }
	\end{subfigure}
	\caption{\FMINCON{} performance profiles for various OPF formulations. \label{fig:profileVariantsFMINCON}} 
\end{figure*}

\begin{figure*}[ht!]
\centering
\begin{subfigure}[b]{0.49\columnwidth}
	\begin{tikzpicture}

\pgfplotstableread{
ALPHA         SOLVER1         SOLVER2         SOLVER3         SOLVER4
1.0         0.48000         0.44000         0.04000         0.00000
1.5         0.64000         0.76000         0.40000         0.08000
2.0         0.80000         0.88000         0.52000         0.28000
2.5         0.80000         0.92000         0.64000         0.40000
3.0         0.84000         0.92000         0.76000         0.40000
3.5         0.84000         0.92000         0.76000         0.48000
4.0         0.84000         0.92000         0.80000         0.56000
4.5         0.88000         0.96000         0.80000         0.56000
5.0         0.88000         0.96000         0.84000         0.56000
5.5         0.88000         0.96000         0.84000         0.56000
6.0         0.88000         0.96000         0.84000         0.56000
6.5         0.88000         0.96000         0.84000         0.60000
7.0         0.92000         0.96000         0.84000         0.60000
7.5         0.96000         0.96000         0.84000         0.60000
8.0         0.96000         0.96000         0.84000         0.60000
8.5         0.96000         0.96000         0.84000         0.64000
9.0         0.96000         0.96000         0.84000         0.64000
9.5         0.96000         0.96000         0.84000         0.64000
10.0         0.96000         0.96000         0.88000         0.64000
10.5         0.96000         0.96000         0.88000         0.64000
11.0         0.96000         0.96000         0.88000         0.64000
11.5         0.96000         0.96000         0.88000         0.64000
12.0         0.96000         0.96000         0.88000         0.64000
12.5         0.96000         0.96000         0.88000         0.68000
13.0         0.96000         0.96000         0.88000         0.68000
13.5         0.96000         0.96000         0.88000         0.72000
14.0         0.96000         0.96000         0.88000         0.72000
14.5         0.96000         0.96000         0.88000         0.72000
15.0         0.96000         0.96000         0.88000         0.76000
15.5         0.96000         0.96000         0.88000         0.76000
16.0         0.96000         0.96000         0.88000         0.76000
16.5         0.96000         0.96000         0.88000         0.76000
17.0         0.96000         0.96000         0.88000         0.76000
17.5         0.96000         0.96000         0.88000         0.76000
18.0         0.96000         0.96000         0.88000         0.76000
18.5         0.96000         0.96000         0.88000         0.76000
19.0         0.96000         0.96000         0.88000         0.76000
19.5         0.96000         0.96000         0.88000         0.76000
20.0         0.96000         0.96000         0.88000         0.76000
20.5         0.96000         0.96000         0.88000         0.76000
21.0         0.96000         0.96000         0.88000         0.76000
21.5         0.96000         0.96000         0.88000         0.76000
22.0         0.96000         0.96000         0.88000         0.76000
22.5         0.96000         0.96000         0.88000         0.76000
23.0         0.96000         0.96000         0.88000         0.76000
23.5         0.96000         0.96000         0.88000         0.76000
24.0         0.96000         0.96000         0.88000         0.76000
24.5         0.96000         0.96000         0.88000         0.76000
25.0         0.96000         0.96000         0.88000         0.76000
25.5         0.96000         0.96000         0.88000         0.76000
26.0         0.96000         0.96000         0.88000         0.76000
26.5         0.96000         0.96000         0.88000         0.76000
27.0         0.96000         0.96000         0.88000         0.76000
27.5         0.96000         0.96000         0.88000         0.76000
28.0         0.96000         0.96000         0.88000         0.76000
28.5         0.96000         0.96000         0.88000         0.76000
29.0         0.96000         0.96000         0.88000         0.76000
29.5         0.96000         0.96000         0.88000         0.76000
30.0         0.96000         0.96000         0.88000         0.76000
30.5         0.96000         0.96000         0.88000         0.76000
31.0         0.96000         0.96000         0.88000         0.76000
31.5         0.96000         0.96000         0.88000         0.76000
32.0         0.96000         0.96000         0.88000         0.76000
32.5         0.96000         0.96000         0.88000         0.76000
33.0         0.96000         0.96000         0.88000         0.76000
33.5         0.96000         0.96000         0.88000         0.76000
34.0         0.96000         0.96000         0.88000         0.76000
34.5         0.96000         0.96000         0.88000         0.76000
35.0         0.96000         0.96000         0.88000         0.76000
35.5         0.96000         0.96000         0.88000         0.76000
36.0         0.96000         0.96000         0.88000         0.76000
36.5         0.96000         0.96000         0.88000         0.76000
37.0         0.96000         0.96000         0.88000         0.76000
37.5         0.96000         0.96000         0.88000         0.76000
38.0         0.96000         0.96000         0.88000         0.76000
38.5         0.96000         0.96000         0.88000         0.76000
39.0         0.96000         0.96000         0.88000         0.76000
39.5         0.96000         0.96000         0.88000         0.76000
40.0         0.96000         0.96000         0.88000         0.80000
40.5         0.96000         0.96000         0.88000         0.80000
41.0         0.96000         0.96000         0.88000         0.80000
41.5         0.96000         0.96000         0.88000         0.80000
42.0         0.96000         0.96000         0.88000         0.80000
}\datatable

  \begin{axis}[name=symb,
  	width=\columnwidth,
  	height=5cm,
	ymin=-0.02,
	ymax=1.02,
    enlarge x limits=0.05,
	axis lines*=left, 
	ymajorgrids, yminorgrids,
    xmajorgrids, xminorgrids,
	xticklabel style={rotate=0, xshift=-0.0cm, anchor=north, font=\scriptsize},
	ytick={0, 0.1, 0.2, 0.3, 0.4, 0.5, 0.6, 0.7, 0.8, 0.9, 1},
    xmax=18,
	yticklabel style={font=\scriptsize},
	ylabel style={font=\scriptsize},
	xlabel style={font=\scriptsize},
	ylabel={$p_m({\alpha})$},
	xlabel={$\alpha$},
	legend style={at={(1.0,0.60)},legend cell align=left,align=right,draw=white!85!black,font=\footnotesize,legend columns=1}
	]
	
    \addplot[color=mycolor1, ultra thick] table[x=ALPHA, y=SOLVER1] {\datatable}; 
    \addlegendentry{Polar-Power}

    \addplot[color=mycolor2, dashdotted, ultra thick] table[x=ALPHA, y=SOLVER2] {\datatable}; 
    \addlegendentry{Rect-Power}
    
    \addplot[color=mycolor3, dashed, ultra thick] table[x=ALPHA, y=SOLVER3] {\datatable}; 
    \addlegendentry{Polar-Current}

	\addplot[color=mycolor4, dotted, ultra thick] table[x=ALPHA, y=SOLVER4] {\datatable}; 
	\addlegendentry{Rect-Current}

  \end{axis}
  
\end{tikzpicture}
	\caption{Overall time ($\alpha_{max}=40$)\label{fig:profileTvariantsIPOPTpardiso} }
\end{subfigure}
\begin{subfigure}[b]{0.5\columnwidth}
	\begin{tikzpicture}

\pgfplotstableread{
ALPHA         SOLVER1         SOLVER2         SOLVER3         SOLVER4
1.0         0.28000         0.60000         0.08000         0.08000
1.2         0.48000         0.80000         0.40000         0.24000
1.4         0.56000         0.88000         0.52000         0.28000
1.5999999999999999         0.60000         0.88000         0.52000         0.36000
1.7999999999999998         0.68000         0.92000         0.52000         0.44000
1.9999999999999998         0.68000         0.92000         0.56000         0.52000
2.1999999999999997         0.72000         0.92000         0.64000         0.52000
2.3999999999999995         0.80000         0.96000         0.64000         0.52000
2.5999999999999996         0.80000         0.96000         0.64000         0.60000
2.8         0.80000         0.96000         0.64000         0.64000
2.9999999999999996         0.80000         0.96000         0.64000         0.64000
3.1999999999999993         0.80000         0.96000         0.64000         0.64000
3.3999999999999995         0.80000         0.96000         0.72000         0.64000
3.5999999999999996         0.84000         0.96000         0.72000         0.64000
3.7999999999999994         0.84000         0.96000         0.76000         0.68000
3.999999999999999         0.84000         0.96000         0.76000         0.68000
4.199999999999999         0.84000         0.96000         0.80000         0.68000
4.3999999999999995         0.84000         0.96000         0.84000         0.68000
4.6         0.84000         0.96000         0.84000         0.68000
4.799999999999999         0.84000         0.96000         0.84000         0.68000
4.999999999999999         0.84000         0.96000         0.84000         0.68000
5.199999999999999         0.84000         0.96000         0.88000         0.68000
5.399999999999999         0.88000         0.96000         0.88000         0.76000
5.599999999999999         0.88000         0.96000         0.88000         0.76000
5.799999999999999         0.88000         0.96000         0.88000         0.76000
5.999999999999999         0.88000         0.96000         0.88000         0.76000
6.199999999999999         0.88000         0.96000         0.88000         0.76000
6.399999999999999         0.88000         0.96000         0.88000         0.76000
6.599999999999999         0.88000         0.96000         0.88000         0.80000
6.799999999999999         0.88000         0.96000         0.88000         0.80000
6.999999999999998         0.92000         0.96000         0.88000         0.80000
7.199999999999998         0.92000         0.96000         0.88000         0.80000
7.399999999999999         0.92000         0.96000         0.88000         0.80000
7.599999999999999         0.92000         0.96000         0.88000         0.80000
7.799999999999999         0.92000         0.96000         0.88000         0.80000
7.999999999999998         0.92000         0.96000         0.88000         0.80000
8.2         0.92000         0.96000         0.88000         0.80000
8.399999999999999         0.92000         0.96000         0.88000         0.80000
8.599999999999998         0.92000         0.96000         0.88000         0.80000
8.799999999999997         0.92000         0.96000         0.88000         0.80000
8.999999999999998         0.92000         0.96000         0.88000         0.80000
9.199999999999998         0.92000         0.96000         0.88000         0.80000
9.399999999999999         0.92000         0.96000         0.88000         0.80000
9.599999999999998         0.96000         0.96000         0.88000         0.80000
9.799999999999997         0.96000         0.96000         0.88000         0.80000
9.999999999999998         0.96000         0.96000         0.88000         0.80000	
}\datatable

  \begin{axis}[name=symb,
  	width=\columnwidth,
  	height=5cm,
	ymin=-0.02,
	ymax=1.02,
    enlarge x limits=0.05,
	axis lines*=left, 
	ymajorgrids, yminorgrids,
    xmajorgrids, xminorgrids,
	xticklabel style={rotate=0, xshift=-0.0cm, anchor=north, font=\scriptsize},
	ytick={0, 0.1, 0.2, 0.3, 0.4, 0.5, 0.6, 0.7, 0.8, 0.9, 1},
    xmax=10,
	yticklabel style={font=\scriptsize},
	ylabel style={font=\scriptsize},
	xlabel style={font=\scriptsize},
	ylabel={$p_m({\alpha})$},
	xlabel={$\alpha$},
	legend style={at={(1.0,0.60)},legend cell align=left,align=right,draw=white!85!black,font=\footnotesize,legend columns=1}
	]
	
    \addplot[color=mycolor1, ultra thick] table[x=ALPHA, y=SOLVER1] {\datatable}; 
    \addlegendentry{Polar-Power}

    \addplot[color=mycolor2, dashdotted, ultra thick] table[x=ALPHA, y=SOLVER2] {\datatable}; 
    \addlegendentry{Rect-Power}
    
    \addplot[color=mycolor3, dashed, ultra thick] table[x=ALPHA, y=SOLVER3] {\datatable}; 
    \addlegendentry{Polar-Current}

	\addplot[color=mycolor4, dotted, ultra thick] table[x=ALPHA, y=SOLVER4] {\datatable}; 
	\addlegendentry{Rect-Current}

  \end{axis}
  
\end{tikzpicture}
	\caption{ Number of iterations ($\alpha_{max}=10$)\label{fig:profileIvariantsIPOPTpardiso} }
\end{subfigure}
\caption{\IPOPT{}-\PARDISO{} performance profiles for various OPF formulations starting from the PF solution.\label{fig:profileVariantsIPOPTpardiso}}

\centering
\begin{subfigure}[b]{0.49\columnwidth}
	\begin{tikzpicture}

\pgfplotstableread{
ALPHA         SOLVER1         SOLVER2         SOLVER3         SOLVER4
1.0         0.76000         0.12000         0.12000         0.00000
1.2         0.92000         0.60000         0.56000         0.16000
1.4         0.96000         0.68000         0.60000         0.32000
1.5999999999999999         0.96000         0.72000         0.68000         0.36000
1.7999999999999998         0.96000         0.72000         0.76000         0.44000
1.9999999999999998         0.96000         0.72000         0.76000         0.48000
2.1999999999999997         0.96000         0.72000         0.80000         0.48000
2.3999999999999995         0.96000         0.80000         0.80000         0.48000
2.5999999999999996         0.96000         0.80000         0.84000         0.48000
2.8         0.96000         0.80000         0.84000         0.52000
2.9999999999999996         0.96000         0.80000         0.84000         0.52000
3.1999999999999993         0.96000         0.80000         0.84000         0.56000
3.3999999999999995         0.96000         0.84000         0.84000         0.56000
3.5999999999999996         0.96000         0.84000         0.84000         0.56000
3.7999999999999994         0.96000         0.84000         0.84000         0.56000
3.999999999999999         0.96000         0.84000         0.84000         0.56000
4.199999999999999         0.96000         0.84000         0.84000         0.56000
4.3999999999999995         0.96000         0.84000         0.84000         0.56000
4.6         0.96000         0.84000         0.84000         0.56000
4.799999999999999         0.96000         0.88000         0.84000         0.56000
4.999999999999999         0.96000         0.88000         0.84000         0.56000
5.199999999999999         0.96000         0.88000         0.84000         0.56000
5.399999999999999         0.96000         0.88000         0.84000         0.56000
5.599999999999999         0.96000         0.88000         0.84000         0.56000
5.799999999999999         0.96000         0.88000         0.84000         0.56000
5.999999999999999         0.96000         0.92000         0.84000         0.56000
6.199999999999999         0.96000         0.92000         0.84000         0.56000
6.399999999999999         0.96000         0.92000         0.84000         0.60000
6.599999999999999         0.96000         0.92000         0.84000         0.60000
6.799999999999999         0.96000         0.92000         0.84000         0.60000
6.999999999999998         0.96000         0.92000         0.84000         0.60000
7.199999999999998         0.96000         0.92000         0.84000         0.60000
7.399999999999999         0.96000         0.92000         0.84000         0.60000
7.599999999999999         0.96000         0.92000         0.84000         0.60000
7.799999999999999         0.96000         0.92000         0.84000         0.60000
7.999999999999998         0.96000         0.92000         0.84000         0.60000
}\datatable

  \begin{axis}[name=symb,
  	width=\columnwidth,
  	height=5cm,
	ymin=-0.02,
	ymax=1.02,
    enlarge x limits=0.05,
	axis lines*=left, 
	ymajorgrids, yminorgrids,
    xmajorgrids, xminorgrids,
	xticklabel style={rotate=0, xshift=-0.0cm, anchor=north, font=\scriptsize},
	ytick={0, 0.1, 0.2, 0.3, 0.4, 0.5, 0.6, 0.7, 0.8, 0.9, 1},
    xmax=8,
	yticklabel style={font=\scriptsize},
	ylabel style={font=\scriptsize},
	xlabel style={font=\scriptsize},
	ylabel={$p_m({\alpha})$},
	xlabel={$\alpha$},
	legend style={at={(1.0,0.50)},legend cell align=left,align=right,draw=white!85!black,font=\footnotesize,legend columns=1}
	]
	
    \addplot[color=mycolor1, ultra thick] table[x=ALPHA, y=SOLVER1] {\datatable}; 
    \addlegendentry{Polar-Power}

    \addplot[color=mycolor2, dashdotted, ultra thick] table[x=ALPHA, y=SOLVER2] {\datatable}; 
    \addlegendentry{Rect-Power}
    
    \addplot[color=mycolor3, dashed, ultra thick] table[x=ALPHA, y=SOLVER3] {\datatable}; 
    \addlegendentry{Polar-Current}

	\addplot[color=mycolor4, dotted, ultra thick] table[x=ALPHA, y=SOLVER4] {\datatable}; 
	\addlegendentry{Rect-Current}

  \end{axis}
  
\end{tikzpicture}
	\caption{Overall time ($\alpha_{max}=7$)\label{fig:profileTvariantsMIPSpardiso} }
\end{subfigure}
\begin{subfigure}[b]{0.5\columnwidth}
	\begin{tikzpicture}

\pgfplotstableread{
ALPHA         SOLVER1         SOLVER2         SOLVER3         SOLVER4
1.0         0.84000         0.12000         0.20000         0.00000
1.15         0.92000         0.36000         0.56000         0.08000
1.2999999999999998         0.96000         0.60000         0.60000         0.20000
1.4499999999999997         0.96000         0.68000         0.68000         0.36000
1.5999999999999996         0.96000         0.72000         0.72000         0.36000
1.7499999999999996         0.96000         0.72000         0.72000         0.44000
1.8999999999999995         0.96000         0.72000         0.76000         0.48000
2.0499999999999994         0.96000         0.76000         0.80000         0.48000
2.1999999999999993         0.96000         0.80000         0.84000         0.48000
2.349999999999999         0.96000         0.80000         0.84000         0.48000
2.499999999999999         0.96000         0.80000         0.84000         0.52000
2.649999999999999         0.96000         0.80000         0.84000         0.52000
2.799999999999999         0.96000         0.84000         0.84000         0.56000
2.949999999999999         0.96000         0.84000         0.84000         0.56000
3.0999999999999988         0.96000         0.84000         0.84000         0.56000
3.2499999999999987         0.96000         0.84000         0.84000         0.56000
3.3999999999999986         0.96000         0.84000         0.84000         0.56000
3.5499999999999985         0.96000         0.84000         0.84000         0.56000
3.6999999999999984         0.96000         0.84000         0.84000         0.56000
3.8499999999999983         0.96000         0.84000         0.84000         0.56000
3.9999999999999982         0.96000         0.88000         0.84000         0.56000
4.149999999999999         0.96000         0.88000         0.84000         0.56000
4.299999999999998         0.96000         0.88000         0.84000         0.56000
4.4499999999999975         0.96000         0.88000         0.84000         0.56000
4.599999999999998         0.96000         0.88000         0.84000         0.56000
4.749999999999998         0.96000         0.92000         0.84000         0.56000
4.899999999999998         0.96000         0.92000         0.84000         0.56000
5.049999999999997         0.96000         0.92000         0.84000         0.56000
5.1999999999999975         0.96000         0.92000         0.84000         0.56000
5.349999999999998         0.96000         0.92000         0.84000         0.56000
5.499999999999997         0.96000         0.92000         0.84000         0.56000
5.649999999999997         0.96000         0.92000         0.84000         0.56000
5.799999999999997         0.96000         0.92000         0.84000         0.56000
5.9499999999999975         0.96000         0.92000         0.84000         0.60000
6.099999999999997         0.96000         0.92000         0.84000         0.60000
}\datatable

  \begin{axis}[name=symb,
  	width=\columnwidth,
  	height=5cm,
	ymin=-0.02,
	ymax=1.02,
    enlarge x limits=0.05,
	axis lines*=left, 
	ymajorgrids, yminorgrids,
    xmajorgrids, xminorgrids,
	xticklabel style={rotate=0, xshift=-0.0cm, anchor=north, font=\scriptsize},
	ytick={0, 0.1, 0.2, 0.3, 0.4, 0.5, 0.6, 0.7, 0.8, 0.9, 1},
    xmax=6,
	yticklabel style={font=\scriptsize},
	ylabel style={font=\scriptsize},
	xlabel style={font=\scriptsize},
	ylabel={$p_m({\alpha})$},
	xlabel={$\alpha$},
	legend style={at={(1.0,0.50)},legend cell align=left,align=right,draw=white!85!black,font=\footnotesize,legend columns=1}
	]
	
    \addplot[color=mycolor1, ultra thick] table[x=ALPHA, y=SOLVER1] {\datatable}; 
    \addlegendentry{Polar-Power}

    \addplot[color=mycolor2, dashdotted, ultra thick] table[x=ALPHA, y=SOLVER2] {\datatable}; 
    \addlegendentry{Rect-Power}
    
    \addplot[color=mycolor3, dashed, ultra thick] table[x=ALPHA, y=SOLVER3] {\datatable}; 
    \addlegendentry{Polar-Current}

	\addplot[color=mycolor4, dotted, ultra thick] table[x=ALPHA, y=SOLVER4] {\datatable}; 
	\addlegendentry{Rect-Current}

  \end{axis}
  
\end{tikzpicture}
	\caption{ Number of iterations ($\alpha_{max}=6$)\label{fig:profileIvariantsMIPSpardiso} }
\end{subfigure}
\caption{\MIPS{}-\PARDISO{} performance profiles for various OPF formulations starting from the PF solution.\label{fig:profileVariantsMIPSpardiso}}
\end{figure*}

\clearpage
\newpage
\input{OPFvariantsBELTISTOS.tex}
\pgfplotstableread{
case polarPower cartPower polarCurrent cartCurrent 
case1951rte             5.17    3.85    6.16    9.34
case2383wp              8.08    9.98    6.17    11.18
case2736sp              3.63    4.15    4.21    6.14
case2737sop             4.16    4.08    4.6 4.87
case2746wop             4.17    4.21    4.48    16.36
case2746wp              3.93    4.38    4.13    4.14
case2868rte             6.16    7.32    11.31   9.7
case2869pegase          15.94   10.31   10.98   22.18
case3012wp              9.37    11.37   35.96   35.49
case3120sp              8.65    12.99   22.4    {}
case3375wp              18.88   11.42   25.87   38.24
case6468rte             27.03   10.74   30.63   18.31
case6470rte             85.71   19.65   22.06   159.11
case6495rte             13.62   13.78   14.41   45.24
case6515rte             11.77   13.38   58.29   171.48
case9241pegase          52.52   221.52  126.28  {}
case\_ACTIVSg2000       3.99    4.81    4.35    8.28
case\_ACTIVSg10k        121.12  16.68   45.53   105
case13659pegase         180.44  369.26  189.62  277.08
case\_ACTIVSg25k        362.33  52.14   96.87   {}
case\_ACTIVSg70k        345.75  214.18  528.73  8500.76
case21k                 330.29  227.3   {} 2839.95
case42k                 2024.7  1055.39 10250.8 14050.54
case99k                 3852.76 6133.73 {} {}
case193k                {} {} {} {}
}\OPFtabletimeIPOPTdefault

\pgfplotstableread{
case polarPower cartPower polarCurrent cartCurrent 
case1951rte             53  33  45  48
case2383wp              43  41  45  45
case2736sp              22  28  25  33
case2737sop             27  26  31  29
case2746wop             28  27  29  48
case2746wp              26  28  27  27
case2868rte             50  50  83  41
case2869pegase          72  33  45  63
case3012wp              40  38  166 93
case3120sp              37  44  125 {}
case3375wp              86  37  81  89
case6468rte             146 41  152 47
case6470rte             393 85  73  193
case6495rte             73  64  60  100
case6515rte             60  60  243 220
case9241pegase          65  149 122 {}
case\_ACTIVSg2000       29  33  30  38
case\_ACTIVSg10k        345 36  117 68
case13659pegase         244 353 226 216
case\_ACTIVSg25k        377 54  67  {}
case\_ACTIVSg70k        104 70  78  449
case21k                 114 68  {} 358
case42k                 180 76  388 406
case99k                 94  94  {} {}
case193k                {} {} {} {}
}\OPFtableitersIPOPTdefault


\begin{table}[H]
    \footnotesize
	\centering
	\caption{Overall time (s) (\IPOPT{}-\PARDISO{}). \label{tab:OPFvariantssummaryIPOPTdefault-t}}
	\pgfplotstabletypeset[
	every head row/.style={ before row={\toprule},after row=\midrule},
	every last row/.style={ after row=\bottomrule},
	precision=2, fixed zerofill, column type={r},
	columns/Statistic/.style={string type,column type=l},
	every row/.style={
		column type=r,
		dec sep align,
		fixed,
		fixed zerofill,
	},
	columns={case, polarPower, polarCurrent, cartPower, cartCurrent},
	every odd row/.style={before row={\rowcolor{tablecolor1}}},
    every even row/.style={before row={\rowcolor{tablecolor2}}},
	columns/case/.style={string type, column type={l}, column name={Benchmark}},
	columns/polarPower/.style={column name={Polar-Power}},
	columns/polarCurrent/.style={column name={Polar-Current}},
	columns/cartPower/.style={column name={Cartesian-Power}},
	columns/cartCurrent/.style={column name={Cartesian-Current}},
	empty cells with={---} 
	]\OPFtabletimeIPOPTdefault
\end{table}
\begin{table}[H]
    \footnotesize
	\centering
	\caption{Number of iterations (\IPOPT{}-\PARDISO{}). \label{tab:OPFvariantssummaryIPOPTdefault-i}}
	\pgfplotstabletypeset[
	every head row/.style={ before row={\toprule},after row=\midrule},
	every last row/.style={ after row=\bottomrule},
	precision=0, fixed zerofill, column type={r},
	columns/Statistic/.style={string type,column type=l},
	every row/.style={
		column type=r,
		dec sep align,
		fixed,
		fixed zerofill,
	},
	columns={case, polarPower, polarCurrent, cartPower, cartCurrent},
	every odd row/.style={before row={\rowcolor{tablecolor1}}},
    every even row/.style={before row={\rowcolor{tablecolor2}}},
	columns/case/.style={string type, column type={l}, column name={Benchmark}},
	columns/polarPower/.style={column name={Polar-Power}},
	columns/polarCurrent/.style={column name={Polar-Current}},
	columns/cartPower/.style={column name={Cartesian-Power}},
	columns/cartCurrent/.style={column name={Cartesian-Current}},
	empty cells with={---} 
	]\OPFtableitersIPOPTdefault
\end{table}
\pgfplotstableread{
case    polarPower  cartPower   polarCurrent    cartCurrent
case1951rte             3.4 3.44    3.5 4.04
case2383wp              4.11    4.36    3.9 5.01
case2736sp              3.24    3.7 3.1 4.21
case2737sop             3.64    3.55    3.46    3.81
case2746wop             3.69    3.81    3.56    4.06
case2746wp              3.97    3.75    3.42    3.7
case2868rte             4.77    6.59    4.63    5.14
case2869pegase          3.89    4.59    4.12    4.97
case3012wp              5.57    5.55    5.38    4.91
case3120sp              4.86    5.65    4.86    5.6
case3375wp              5.4 5   5.04    5.7
case6468rte             6.65    8.44    6.18    7.96
case6470rte             12.51   19.51   11.06   19.73
case6495rte             10.17   12.51   10.05   13.18
case6515rte             9.86    12.27   11.64   19.43
case9241pegase          17.76   16.24   30.79   16.8
case\_ACTIVSg2000       4   4.45    4.3 4.67
case\_ACTIVSg10k        13.99   14.39   13.26   15.3
case13659pegase         86.76   144.02  66.78   137.24
case\_ACTIVSg25k        50.06   191.89  37.01   117.16
case\_ACTIVSg70k        152.56  339.94  140.52  265.47
case21k                 {} {} 825.15  {}
case42k                 {} {} 8745.17 {}
case99k                 {} {} {} {}
case193k                {} {} {} {}
}\OPFtabletimeIPOPThsl

\pgfplotstableread{
case    polarPower  cartPower   polarCurrent    cartCurrent
case1951rte             34  33  37  37
case2383wp              39  41  38  46
case2736sp              23  28  23  33
case2737sop             28  26  26  29
case2746wop             29  27  28  32
case2746wp              28  28  26  27
case2868rte             39  50  39  41
case2869pegase          29  33  34  37
case3012wp              48  38  48  40
case3120sp              38  44  40  46
case3375wp              43  37  40  45
case6468rte             38  41  39  41
case6470rte             68  85  67  92
case6495rte             60  64  64  71
case6515rte             57  60  72  89
case9241pegase          45  41  76  47
case\_ACTIVSg2000       29  33  30  36
case\_ACTIVSg10k        35  36  36  38
case13659pegase         246 343 282 374
case\_ACTIVSg25k        50  155 51  76
case\_ACTIVSg70k        65  92  64  74
case21k                 {} {} 388 {}
case42k                 {} {} 343 {}
case99k                 {} {} {} {}
case193k                {} {} {} {}
}\OPFtableitersIPOPThsl


\begin{table}[H]
    \footnotesize
	\centering
	\caption{Overall time (s) (IPOPT-MA57). \label{tab:OPFvariantssummaryIPOPThsl-t}}
	\pgfplotstabletypeset[
	every head row/.style={ before row={\toprule},after row=\midrule},
	every last row/.style={ after row=\bottomrule},
	precision=2, fixed zerofill, column type={r},
	columns/Statistic/.style={string type,column type=l},
	every row/.style={
		column type=r,
		dec sep align,
		fixed,
		fixed zerofill,
	},
	columns={case, polarPower, polarCurrent, cartPower, cartCurrent},
	every odd row/.style={before row={\rowcolor{tablecolor1}}},
    every even row/.style={before row={\rowcolor{tablecolor2}}},
	columns/case/.style={string type, column type={l}, column name={Benchmark}},
	columns/polarPower/.style={column name={Polar-Power}},
	columns/polarCurrent/.style={column name={Polar-Current}},
	columns/cartPower/.style={column name={Cartesian-Power}},
	columns/cartCurrent/.style={column name={Cartesian-Current}},
	empty cells with={---} 
	]\OPFtabletimeIPOPThsl
\end{table}
\begin{table}[H]
    \footnotesize
	\centering
	\caption{Number of iterations (IPOPT-MA57). \label{tab:OPFvariantssummaryIPOPThsl-i}}
	\pgfplotstabletypeset[
	every head row/.style={ before row={\toprule},after row=\midrule},
	every last row/.style={ after row=\bottomrule},
	precision=0, fixed zerofill, column type={r},
	columns/Statistic/.style={string type,column type=l},
	every row/.style={
		column type=r,
		dec sep align,
		fixed,
		fixed zerofill,
	},
	columns={case, polarPower, polarCurrent, cartPower, cartCurrent},
	every odd row/.style={before row={\rowcolor{tablecolor1}}},
    every even row/.style={before row={\rowcolor{tablecolor2}}},
	columns/case/.style={string type, column type={l}, column name={Benchmark}},
	columns/polarPower/.style={column name={Polar-Power}},
	columns/polarCurrent/.style={column name={Polar-Current}},
	columns/cartPower/.style={column name={Cartesian-Power}},
	columns/cartCurrent/.style={column name={Cartesian-Current}},
	empty cells with={---} 
	]\OPFtableitersIPOPThsl
\end{table}

\pgfplotstableread{
case polarPower cartPower polarCurrent cartCurrent 
case1951rte             3.71    8.47    4.1 14.53
case2383wp              5.14    5.51    5.17    6.87
case2736sp              4.44    5.23    4.55    7.7
case2737sop             4.57    5.55    4.45    5.44
case2746wop             4.94    4.83    5.03    5.23
case2746wp              4.85    5.85    5.1 6.25
case2868rte             5.28    6.76    {} 96.18
case2869pegase          6.11    9.26    6.63    {}
case3012wp              5.5 6.38    5.41    6.67
case3120sp              6.36    6.94    6.14    6.31
case3375wp              5.91    7.17    6.42    7.46
case6468rte             14.18   {} 21.91   {}
case6470rte             14.16   67.94   33.07   {}
case6495rte             24.29   {} 42.42   58.21
case6515rte             19.36   55.22   19.44   {}
case9241pegase          25.9    {} 26.51   {}
case\_ACTIVSg2000       4.38    4.71    8.02    7.42
case\_ACTIVSg10k        45.32   34.48   62.39   48.05
case13659pegase         30.53   {} {} {}
case\_ACTIVSg25k        108.18  92.13   104.14  150.27
case\_ACTIVSg70k        {} 404.02  423.96  728.76
case21k                 241.2   253.8   230.33  250.33
case42k                 1957.95 2023.77 1966.92 2090.61
case99k                 {} {} {} {}
case193k                {} {} {} {}
}\OPFtabletimeMIPSsc

\pgfplotstableread{
case polarPower cartPower polarCurrent cartCurrent 
case1951rte             28  62  32  106
case2383wp              35  40  35  46
case2736sp              26  32  27  46
case2737sop             26  30  26  30
case2746wop             29  29  29  31
case2746wp              28  37  30  37
case2868rte             31  41  {} 428
case2869pegase          32  47  36  {}
case3012wp              30  36  29  37
case3120sp              33  38  34  36
case3375wp              30  37  33  40
case6468rte             44  {} 66  {}
case6470rte             44  178 94  {}
case6495rte             76  {} 131 165
case6515rte             61  160 59  {}
case9241pegase          46  {} 51  {}
case\_ACTIVSg2000       29  31  52  48
case\_ACTIVSg10k        81  66  108 83
case13659pegase         54  {} {} {}
case\_ACTIVSg25k        76  72  78  103
case\_ACTIVSg70k        {} 96  96  152
case21k                 53  61  53  61
case42k                 64  73  64  74
case99k                 {} {} {} {}
case193k                {} {} {} {}
}\OPFtableitersMIPSsc

\begin{table}[H]
    \footnotesize
	\centering
	\caption{Overall time (s) (\MIPS{}). \label{tab:OPFvariantssummaryMIPSsc-t}}
	\pgfplotstabletypeset[
	every head row/.style={ before row={\toprule},after row=\midrule},
	every last row/.style={ after row=\bottomrule},
	precision=2, fixed zerofill, column type={r},
	columns/Statistic/.style={string type,column type=l},
	every row/.style={
		column type=r,
		dec sep align,
		fixed,
		fixed zerofill,
	},
	columns={case, polarPower, polarCurrent, cartPower, cartCurrent},
	every odd row/.style={before row={\rowcolor{tablecolor1}}},
    every even row/.style={before row={\rowcolor{tablecolor2}}},
	columns/case/.style={string type, column type={l}, column name={Benchmark}},
	columns/polarPower/.style={column name={Polar-Power}},
	columns/polarCurrent/.style={column name={Polar-Current}},
	columns/cartPower/.style={column name={Cartesian-Power}},
	columns/cartCurrent/.style={column name={Cartesian-Current}},
	empty cells with={---} 
	]\OPFtabletimeMIPSsc
\end{table}
\begin{table}[H]
    \footnotesize
	\centering
	\caption{Number of iterations (\MIPS{}). \label{tab:OPFvariantssummaryMIPSsc-i}}
	\pgfplotstabletypeset[
	every head row/.style={ before row={\toprule},after row=\midrule},
	every last row/.style={ after row=\bottomrule},
	precision=0, fixed zerofill, column type={r},
	columns/Statistic/.style={string type,column type=l},
	every row/.style={
		column type=r,
		dec sep align,
		fixed,
		fixed zerofill,
	},
	columns={case, polarPower, polarCurrent, cartPower, cartCurrent},
	every odd row/.style={before row={\rowcolor{tablecolor1}}},
    every even row/.style={before row={\rowcolor{tablecolor2}}},
	columns/case/.style={string type, column type={l}, column name={Benchmark}},
	columns/polarPower/.style={column name={Polar-Power}},
	columns/polarCurrent/.style={column name={Polar-Current}},
	columns/cartPower/.style={column name={Cartesian-Power}},
	columns/cartCurrent/.style={column name={Cartesian-Current}},
	empty cells with={---} 
	]\OPFtableitersMIPSsc
\end{table}

\pgfplotstableread{
case polarPower cartPower polarCurrent cartCurrent 
case1951rte             4.75    10.46   5.69    29.97
case2383wp              6.6 7.38    6.55    8.63
case2736sp              5.98    6.82    6.05    9.62
case2737sop             5.68    6.28    5.91    6.32
case2746wop             6.4 6.36    6.84    6.87
case2746wp              6.26    7.67    6.66    7.81
case2868rte             6.82    8.66    13.79   {}
case2869pegase          7.66    12.18   8.39    {}
case3012wp              7.13    8.44    7.12    8.38
case3120sp              7.6 8.79    7.7 8.16
case3375wp              7.37    8.82    8.11    9.85
case6468rte             16.57   {} 19.93   {}
case6470rte             17.31   81.82   42.32   {}
case6495rte             30.8    72.56   49.91   94.68
case6515rte             24.71   79.67   29.13   68.41
case9241pegase          31.91   185.65  34.71   {}
case\_ACTIVSg2000       5.7 6.29    10.19   9.28
case\_ACTIVSg10k        60.06   47.41   69.08   62.79
case13659pegase         43.54   {} {} {}
case\_ACTIVSg25k        138.79  120.95  138.73  192.83
case\_ACTIVSg70k        {} 500.27  486.22  877.68
case21k                 98  110.55  {} {}
case42k                 303.38  333.91  449.81  {}
case99k                 1081.72 1174.49 {} {}
case193k                2982.83 3235.83 {} {}
}\OPFtabletimeMIPSscPardiso

\pgfplotstableread{
case polarPower cartPower polarCurrent cartCurrent 
case1951rte             28  61  32  164
case2383wp              35  40  35  46
case2736sp              26  32  27  46
case2737sop             26  30  26  30
case2746wop             29  29  29  31
case2746wp              28  37  30  37
case2868rte             31  41  63  {}
case2869pegase          32  47  36  {}
case3012wp              30  36  29  37
case3120sp              33  38  34  36
case3375wp              30  37  33  40
case6468rte             44  {} 54  {}
case6470rte             44  172 95  {}
case6495rte             76  153 116 205
case6515rte             61  165 67  152
case9241pegase          46  218 51  {}
case\_ACTIVSg2000       29  31  52  48
case\_ACTIVSg10k        81  67  94  84
case13659pegase         54  {} {} {}
case\_ACTIVSg25k        76  72  78  103
case\_ACTIVSg70k        {} 96  92  152
case21k                 53  61  {} {}
case42k                 64  73  92  {}
case99k                 78  85  {} {}
case193k                92  101 {} {}
}\OPFtableitersMIPSscPardiso

\begin{table}[H]
    \footnotesize
	\centering
	\caption{Overall time (s) (\MIPS{}-\PARDISO{}). \label{tab:OPFvariantssummaryMIPSscPardiso-t}}
	\pgfplotstabletypeset[
	every head row/.style={ before row={\toprule},after row=\midrule},
	every last row/.style={ after row=\bottomrule},
	precision=2, fixed zerofill, column type={r},
	columns/Statistic/.style={string type,column type=l},
	every row/.style={
		column type=r,
		dec sep align,
		fixed,
		fixed zerofill,
	},
	columns={case, polarPower, polarCurrent, cartPower, cartCurrent},
	every odd row/.style={before row={\rowcolor{tablecolor1}}},
    every even row/.style={before row={\rowcolor{tablecolor2}}},
	columns/case/.style={string type, column type={l}, column name={Benchmark}},
	columns/polarPower/.style={column name={Polar-Power}},
	columns/polarCurrent/.style={column name={Polar-Current}},
	columns/cartPower/.style={column name={Cartesian-Power}},
	columns/cartCurrent/.style={column name={Cartesian-Current}},
	empty cells with={---} 
	]\OPFtabletimeMIPSscPardiso
\end{table}
\begin{table}[H]
    \footnotesize
	\centering
	\caption{Number of iterations (\MIPS{}-\PARDISO{}). \label{tab:OPFvariantssummaryMIPSscPardiso-i}}
	\pgfplotstabletypeset[
	every head row/.style={ before row={\toprule},after row=\midrule},
	every last row/.style={ after row=\bottomrule},
	precision=0, fixed zerofill, column type={r},
	columns/Statistic/.style={string type,column type=l},
	every row/.style={
		column type=r,
		dec sep align,
		fixed,
		fixed zerofill,
	},
	columns={case, polarPower, polarCurrent, cartPower, cartCurrent},
	every odd row/.style={before row={\rowcolor{tablecolor1}}},
    every even row/.style={before row={\rowcolor{tablecolor2}}},
	columns/case/.style={string type, column type={l}, column name={Benchmark}},
	columns/polarPower/.style={column name={Polar-Power}},
	columns/polarCurrent/.style={column name={Polar-Current}},
	columns/cartPower/.style={column name={Cartesian-Power}},
	columns/cartCurrent/.style={column name={Cartesian-Current}},
	empty cells with={---} 
	]\OPFtableitersMIPSscPardiso
\end{table}

\pgfplotstableread{
case polarPower cartPower polarCurrent cartCurrent 
case1951rte          41.3033080000 10.7119680000 172.3579880000 26.9319370000
case2383wp           60.4933000000 92.9733820000 186.7791310000 25.5599920000
case2736sp           {} 15.8332620000 53.7199290000 19.3071760000
case2737sop          12.8088140000 11.9183750000 10.8141350000 11.7843110000
case2746wop          20.9215910000 9.8113540000 386.6310380000 14.6873200000
case2746wp           138.3175610000 13.7478810000 {} 13.3572820000
case2868rte          47.8652820000 35.2816170000 52.7816860000 20.5851340000
case2869pegase       21.9097440000 24.4455720000 28.9021400000 92.3998560000
case3012wp           390.1926990000 22.2998140000 {} 22.4550720000
case3120sp           {} 40.3665000000 {} 51.0952170000
case3375wp           {} 30.8690640000 {} 57.1750970000
case6468rte          {} 46.8318290000 64.3115860000 128.4643690000
case6470rte          89.3804650000 36.2921460000 72.5013020000 442.6277050000
case6495rte          31.7369140000 42.5896770000 174.0086030000 474.8213670000
case6515rte          105.2720130000 81.7043820000 191.5430240000 149.0664570000
case9241pegase       292.1678670000 226.6760130000 386.1305830000 1137.0500230000
case\_ACTIVSg2000    11.4091030000 12.0267480000 11.2158370000 15.0732060000
case\_ACTIVSg10k     25.4822690000 39.2320650000 43.1478830000 235.0943280000
case13659pegase      2510.0285900000 6468.6343610000 5061.6742940000 1850.4865000000
case\_ACTIVSg25k     214.4695870000 523.4196650000 217.2233180000 381.7470900000
case\_ACTIVSg70k     798.0774980000 15290.3365080000 1439.9979110000 2961.6604940000
case21k               2460.1359420000 367.3697480000 5140.3383320000 406.9842160000
case42k              {} 4602.00 {} 5093.77
case99k              53656.32 42371.89 {} {}
case193k             {} {}  {} {}
}\OPFtabletimeFMINCONa

\pgfplotstableread{
case polarPower cartPower polarCurrent cartCurrent 
case1951rte          77 39 90 45
case2383wp           184 60 465 70
case2736sp           {} 42 162 50
case2737sop          42 32 40 32
case2746wop          61 25 342 38
case2746wp           333 35 {} 35
case2868rte          58 31 51 29
case2869pegase       50 52 61 96
case3012wp           434 51 {} 53
case3120sp           {} 72 {} 106
case3375wp           {} 65 {} 19
case6468rte          {} 51 50 71
case6470rte          57 35 73 54
case6495rte          45 55 141 68
case6515rte          108 65 94 58
case9241pegase       81 75 92 94
case\_ACTIVSg2000    35 35 36 40
case\_ACTIVSg10k     18 27 35 50
case13659pegase      63 96 98 62
case\_ACTIVSg25k     48 100 51 65
case\_ACTIVSg70k     79 129 73 160
case21k              437 59 410 68
case42k              {} 64 {} 64
case99k              297 83 {} {}  
case193k             {} {}  {} {}  
}\OPFtableitersFMINCONa

\pgfplotstableread{
case polarPower cartPower polarCurrent cartCurrent 
case1951rte              42.13 14.09 46.21 32.63
case2383wp               69.35 102.23 163.32 30.05
case2736sp               324.89 19.03 60.62 22.62
case2737sop              14.92 15.33 13.87 14.60
case2746wop              24.39 13.01 511.82 17.44
case2746wp               165.37 17.27 {} 16.61
case2868rte              54.29 40.68 284.08 24.92
case2869pegase           25.07 30.08 32.82 103.72
case3012wp               {} 25.37 230.33 27.22
case3120sp               {} 45.17 {} 56.38
case3375wp               183.98 38.25 {} 60.56
case6468rte              {} 51.51 69.92 280.19
case6470rte              92.65 42.04 86.24 294.48
case6495rte              37.76 48.09 538.31 701.07
case6515rte              123.73 93.03 183.28 157.91
case9241pegase           281.36 247.55 311.00 1126.93
case\_ACTIVSg2000        12.81 15.30 14.73 19.94
case\_ACTIVSg10k         26.88 52.87 56.26 274.61
case13659pegase          1866.80 6085.81 1140.71 4190.40
case\_ACTIVSg25k         210.12 549.13 224.92 397.45
case\_ACTIVSg70k         786.05 15831.40 1464.73 10398.00
case21k                  {} 434.96 {} 486.57
case42k                  8428.31 1322.59 {} 1461.98
case99k                  20554.62 16603.51 {} 20352.13
case193k                 70150.87 14281.20 {} 45494.82
}\OPFtabletimeFMINCONb

\pgfplotstableread{
case polarPower cartPower polarCurrent cartCurrent 
case1951rte           77 39 78 45
case2383wp            184 60 462 70
case2736sp            455 42 159 50
case2737sop           42 32 40 32
case2746wop           61 25 407 38
case2746wp            343 35 {} 35
case2868rte           58 31 66 29
case2869pegase        50 52 61 87
case3012wp            {} 51 495 53
case3120sp            {} 72 {} 106
case3375wp            383 65 {} 19
case6468rte           {} 51 50 71
case6470rte           57 35 73 45
case6495rte           45 55 155 74
case6515rte           108 65 88 58
case9241pegase        81 75 94 81
case\_ACTIVSg2000     35 35 36 40
case\_ACTIVSg10k      18 27 35 50
case13659pegase       81 75 79 126
case\_ACTIVSg25k      48 100 51 65
case\_ACTIVSg70k      79 88 78 91
case21k               {} 59 {} 68
case42k               435 64 {} 64
case99k               324 83 {} 82
case193k              482 106 {} 108
}\OPFtableitersFMINCONb

\pgfplotstableread{
	case polarPower cartPower polarCurrent cartCurrent 
	case1951rte             9.29    8.29    72.99   86.23
	case2383wp              13.25   53.41   16.34   19.24
	case2736sp              10.7    12.01   10.81   12.47
	case2737sop             7.59    9.18    8.75    9.77
	case2746wop             9.28    8.5 9.93    9.64
	case2746wp              10.47   11.75   9.64    11.52
	case2868rte             22.09   23.72   310.19  593.45
	case2869pegase          8.95    13.56   8.7 13.61
	case3012wp              {}   16.63   {}   18.98
	case3120sp              14  17.38   12.45   23.29
	case3375wp              {}   22.84   109.52  22.09
	case6468rte             {}   34.39   81.76   52.77
	case6470rte             110.68  69.07   424.97  1599.78
	case6495rte             25.56   29.01   164.49  115.46
	case6515rte             35.61   61.98   118.64  184.59
	case9241pegase          85.86   73.96   204.36  139.45
	case\_ACTIVSg2000       9.44    11.11   9.98    12.17
	case\_ACTIVSg10k        38.61   58.12   46.11   105.04
	case13659pegase         3975.64 {}   4894.19 {}
	case\_ACTIVSg25k        183.57  388.48  226.34  291.82
	case\_ACTIVSg70k        644.3   1619.39 1269.91 2657.26
	case21k                 1548.28 640.1   {}   419.82
	case42k                 {}   1421.21 {}   1693.47
	case99k                 17325.98    6061.47 {}   11369.02
	case193k                {}   16911.01    {}   14333.75
}\OPFtabletimeFMINCONc

\pgfplotstableread{
	case polarPower cartPower polarCurrent cartCurrent 
	case1951rte             48  39  119 66
	case2383wp              62  89  84  77
	case2736sp              46  44  47  47
	case2737sop             31  32  38  35
	case2746wop             39  27  43  33
	case2746wp              45  38  43  41
	case2868rte             62  32  84  68
	case2869pegase          34  40  30  37
	case3012wp              {}   56  {}   64
	case3120sp              56  57  47  68
	case3375wp              {}   73  379 25
	case6468rte             {}   44  58  51
	case6470rte             67  49  67  372
	case6495rte             60  57  128 68
	case6515rte             70  54  67  55
	case9241pegase          46  47  73  58
	case\_ACTIVSg2000       41  40  39  44
	case\_ACTIVSg10k        52  60  59  111
	case13659pegase         77  {}   324 {}
	case\_ACTIVSg25k        65  112 76  82
	case\_ACTIVSg70k        83  160 111 171
	case21k                 415 127 {}   96
	case42k                 {}   104 {}   106
	case99k                 348 101 {}   103
	case193k                {}   135 {}   137
}\OPFtableitersFMINCONc
\begin{table}[H]
	\footnotesize
	\centering
	\caption{Time to solution in seconds (\FMINCON{} 2018b). \label{tab:OPFvariantssummaryFMINCONb-t}}
	\pgfplotstabletypeset[
	every head row/.style={ before row={\toprule},after row=\midrule},
	every last row/.style={ after row=\bottomrule},
	precision=2, fixed zerofill, column type={r},
	columns/Statistic/.style={string type,column type=l},
	every row/.style={
		column type=r,
		dec sep align,
		fixed,
		fixed zerofill,
	},
	columns={case, polarPower, polarCurrent, cartPower, cartCurrent},
	every odd row/.style={before row={\rowcolor{tablecolor1}}},
	every even row/.style={before row={\rowcolor{tablecolor2}}},
	columns/case/.style={string type, column type={l}, column name={Benchmark}},
	columns/polarPower/.style={column name={Polar-Power}},
	columns/polarCurrent/.style={column name={Polar-Current}},
	columns/cartPower/.style={column name={Cartesian-Power}},
	columns/cartCurrent/.style={column name={Cartesian-Current}},
	empty cells with={---} 
	]\OPFtabletimeFMINCONc
\end{table}
\begin{table}[H]
	\footnotesize
	\centering
	\caption{Iterations until convergence (\FMINCON{} 2018b). \label{tab:OPFvariantssummaryFMINCONb-i}}
	\pgfplotstabletypeset[
	every head row/.style={ before row={\toprule},after row=\midrule},
	every last row/.style={ after row=\bottomrule},
	precision=0, fixed zerofill, column type={r},
	columns/Statistic/.style={string type,column type=l},
	every row/.style={
		column type=r,
		dec sep align,
		fixed,
		fixed zerofill,
	},
	columns={case, polarPower, polarCurrent, cartPower, cartCurrent},
	every odd row/.style={before row={\rowcolor{tablecolor1}}},
	every even row/.style={before row={\rowcolor{tablecolor2}}},
	columns/case/.style={string type, column type={l}, column name={Benchmark}},
	columns/polarPower/.style={column name={Polar-Power}},
	columns/polarCurrent/.style={column name={Polar-Current}},
	columns/cartPower/.style={column name={Cartesian-Power}},
	columns/cartCurrent/.style={column name={Cartesian-Current}},
	empty cells with={---} 
	]\OPFtableitersFMINCONc
\end{table}

\pgfplotstableread{
case polarPower cartPower polarCurrent cartCurrent 
case1951rte           10.22 4.40 5.44 8.62
case2383wp            4.92 4.77 4.84 4.91
case2736sp            {} 4.39 4.85 4.53
case2737sop           3.97 4.55 4.21 4.21
case2746wop           45.41 5.27 3.94 4.13
case2746wp            13.07 4.51 4.61 4.67
case2868rte           {} 5.48 9.32 10.41
case2869pegase        4.74 5.63 5.69 6.46
case3012wp            7.02 5.61 5.26 5.36
case3120sp            5.35 5.12 5.51 5.82
case3375wp            6.73 5.79 5.89 6.07
case6468rte           42.25 8.84 9.35 14.71
case6470rte           10.17 10.58 12.52 14.80
case6495rte           9.33 11.81 13.31 21.10
case6515rte           12.58 10.01 13.80 17.39
case9241pegase        16.87 18.16 21.66 119.23
case\_ACTIVSg2000     4.37 4.41 4.30 4.77
case\_ACTIVSg10k      12.47 14.02 11.98 15.82
case13659pegase       93.37 133.52 {} {}
case\_ACTIVSg25k      46.64 55.67 47.32 51.77
case\_ACTIVSg70k      162.04 200.82 150.89 216.53
case21k               206.98 66.96 57.54 68.14
case42k               1147.66 157.81 160.24 156.70
case99k               2398.32 531.30 491.46 517.16
case193k              2778.65 1085.04 949.67 987.93
}\OPFtabletimeKNITROa

\pgfplotstableread{
case polarPower cartPower polarCurrent cartCurrent 
case1951rte           97 26 42 68
case2383wp            32 29 32 29
case2736sp            {} 21 25 22
case2737sop           20 21 23 21
case2746wop           440 19 18 20
case2746wp            117 22 25 22
case2868rte           {} 24 65 66
case2869pegase        25 28 34 33
case3012wp            48 25 28 24
case3120sp            28 28 30 30
case3375wp            39 26 27 26
case6468rte           185 31 39 51
case6470rte           33 35 45 46
case6495rte           36 40 54 68
case6515rte           46 35 52 59
case9241pegase        32 37 47 240
case\_ACTIVSg2000     21 22 20 23
case\_ACTIVSg10k      23 28 26 32
case13659pegase       203 234 {} {}
case\_ACTIVSg25k      45 47 48 48
case\_ACTIVSg70k      52 55 53 57
case21k               181 42 43 42
case42k               325 48 48 47
case99k               150 54 48 54
case193k              145 54 57 53
}\OPFtableitersKNITROa

\pgfplotstableread{
case polarPower cartPower polarCurrent cartCurrent 
case1951rte             5.24    3.18    3.87    5.44
case2383wp              3.82    3.94    3.67    4
case2736sp              3.23    3.56    3.32    3.63
case2737sop             3.51    3.74    3.39    3.73
case2746wop             3.44    3.66    3.43    3.81
case2746wp              3.44    3.76    3.33    3.61
case2868rte             5.91    5.83    7.73    9.43
case2869pegase          3.79    4.57    3.95    4.8
case3012wp              7.02    4.4 4.17    4.31
case3120sp              3.92    4.57    3.9 4.52
case3375wp              5.51    4.9 4.38    5.08
case6468rte             27.75   7.93    7.85    11.22
case6470rte             8.11    9.04    11.35   14.35
case6495rte             8.78    10.59   9.68    14.08
case6515rte             8.86    9.35    13.05   21.64
case9241pegase          13.42   16.74   19.36   45.23
case\_ACTIVSg2000       3.4 3.67    3.79    4.25
case\_ACTIVSg10k        11.73   15.17   11.51   14.89
case13659pegase         67.05   113.02  133.34  410.86
case\_ACTIVSg25k        48.36   56.19   46.98   56.53
case\_ACTIVSg70k        189.45  240.83  187.59  272.29
case21k                 189.18  69.91   57.04   75
case42k                 1239.42 205.36  176.63  187.31
case99k                 {} 728.7   704.54  719.09
case193k                3095.06 1985.92 1540.56 1943.29
}\OPFtabletimeKNITROb

\pgfplotstableread{
case polarPower cartPower polarCurrent cartCurrent 
case1951rte             54  25  37  49
case2383wp              31  30  29  30
case2736sp              19  21  21  22
case2737sop             22  23  22  23
case2746wop             22  21  22  23
case2746wp              21  22  21  21
case2868rte             47  31  65  47
case2869pegase          23  27  25  28
case3012wp              60  24  29  24
case3120sp              24  28  23  28
case3375wp              40  27  26  27
case6468rte             152 31  36  44
case6470rte             32  33  45  49
case6495rte             38  44  43  58
case6515rte             40  36  56  70
case9241pegase          28  33  47  91
case\_ACTIVSg2000       20  22  24  25
case\_ACTIVSg10k        26  31  27  32
case13659pegase         131 143 269 251
case\_ACTIVSg25k        45  47  45  49
case\_ACTIVSg70k        52  56  53  63
case21k                 172 42  43  42
case42k                 396 51  48  47
case99k                 {} 55  57  54
case193k                157 63  66  63
}\OPFtableitersKNITROb

\pgfplotstableread{
case polarPower cartPower polarCurrent cartCurrent 
case1951rte           3.08  3.28    4.93    5.19
case2383wp            3.9   4.26    3.79    4.29
case2736sp            3.28  3.45    3.31    3.56
case2737sop           3.39  3.66    3.44    3.65
case2746wop           3.4   3.82    3.31    3.77
case2746wp            3.74  3.63    3.49    3.71
case2868rte           5.34  4.65    5.49    6.56
case2869pegase        3.95  4.49    4.01    5.36
case3012wp            3.93  4.7 4.04    4.42
case3120sp            4.35  4.81    4.17    4.68
case3375wp            4.12  4.77    4.14    5.05
case6468rte           6.6   8.05    7.37    10.12
case6470rte           7.63  8.86    10.35   14.84
case6495rte           8.52  10.2    9.32    13.97
case6515rte           8.36  10.16   14.25   13.99
case9241pegase        28.36 15.94   18.6    33.14
case\_ACTIVSg2000     3.29  3.6 3.46    3.95
case\_ACTIVSg10k      10.78 13.82   10.84   13.98
case13659pegase       134.32    100.02  896.73  {}
case\_ACTIVSg25k      45.79 55.72   46.66   55.06
case\_ACTIVSg70k      187.65    218.07  182.97  235.23
case21k               47.05 64.65   49.06   66.25
case42k               161.59    175.47  158.84  202.64
case99k               753.07    892.22  620.7   729.39
case193k              1889.51   1899.8 1111.63 1778.61
}\OPFtabletimeKNITROc

\pgfplotstableread{
case polarPower cartPower polarCurrent cartCurrent 
case1951rte           21    22  45  45
case2383wp            29    32  29  32
case2736sp            17    18  18  19
case2737sop           18    20  18  20
case2746wop           19    22  18  21
case2746wp            19    20  20  20
case2868rte           27    28  39  41
case2869pegase        22    25  24  29
case3012wp            23    26  23  24
case3120sp            25    29  25  28
case3375wp            23    25  22  26
case6468rte           27    31  32  39
case6470rte           29    31  43  55
case6495rte           37    43  42  56
case6515rte           31    38  54  52
case9241pegase        64    32  47  70
case\_ACTIVSg2000     17    20  19  22
case\_ACTIVSg10k      24    29  25  30
case13659pegase       132   183 494 {}
case\_ACTIVSg25k      43    46  43  48
case\_ACTIVSg70k      53    55  55  56
case21k               38    42  38  41
case42k               44    50  42  47
case99k               50    55  49  54
case193k              57    65 57  62
}\OPFtableitersKNITROc

\begin{table}[H]
    \footnotesize
	\centering
	\caption{Overall time (s) (\KNITRO{} 11). \label{tab:OPFvariantssummaryKNITROb-t}}
	\pgfplotstabletypeset[
	every head row/.style={ before row={\toprule},after row=\midrule},
	every last row/.style={ after row=\bottomrule},
	precision=2, fixed zerofill, column type={r},
	columns/Statistic/.style={string type,column type=l},
	every row/.style={
		column type=r,
		dec sep align,
		fixed,
		fixed zerofill,
	},
	columns={case, polarPower, polarCurrent, cartPower, cartCurrent},
	every odd row/.style={before row={\rowcolor{tablecolor1}}},
    every even row/.style={before row={\rowcolor{tablecolor2}}},
	columns/case/.style={string type, column type={l}, column name={Benchmark}},
	columns/polarPower/.style={column name={Polar-Power}},
	columns/polarCurrent/.style={column name={Polar-Current}},
	columns/cartPower/.style={column name={Cartesian-Power}},
	columns/cartCurrent/.style={column name={Cartesian-Current}},
	empty cells with={---} 
	]\OPFtabletimeKNITROb
\end{table}
\begin{table}[H]
    \footnotesize
	\centering
	\caption{Number of iterations (\KNITRO{} 11). \label{tab:OPFvariantssummaryKNITROb-i}}
	\pgfplotstabletypeset[
	every head row/.style={ before row={\toprule},after row=\midrule},
	every last row/.style={ after row=\bottomrule},
	precision=0, fixed zerofill, column type={r},
	columns/Statistic/.style={string type,column type=l},
	every row/.style={
		column type=r,
		dec sep align,
		fixed,
		fixed zerofill,
	},
	columns={case, polarPower, polarCurrent, cartPower, cartCurrent},
	every odd row/.style={before row={\rowcolor{tablecolor1}}},
    every even row/.style={before row={\rowcolor{tablecolor2}}},
	columns/case/.style={string type, column type={l}, column name={Benchmark}},
	columns/polarPower/.style={column name={Polar-Power}},
	columns/polarCurrent/.style={column name={Polar-Current}},
	columns/cartPower/.style={column name={Cartesian-Power}},
	columns/cartCurrent/.style={column name={Cartesian-Current}},
	empty cells with={---} 
	]\OPFtableitersKNITROb
\end{table}

\begin{table}[H]
	\footnotesize
	\centering
	\caption{Overall time (s) (\KNITRO{} 12). \label{tab:OPFvariantssummaryKNITROb-t}}
	\pgfplotstabletypeset[
	every head row/.style={ before row={\toprule},after row=\midrule},
	every last row/.style={ after row=\bottomrule},
	precision=2, fixed zerofill, column type={r},
	columns/Statistic/.style={string type,column type=l},
	every row/.style={
		column type=r,
		dec sep align,
		fixed,
		fixed zerofill,
	},
	columns={case, polarPower, polarCurrent, cartPower, cartCurrent},
	every odd row/.style={before row={\rowcolor{tablecolor1}}},
	every even row/.style={before row={\rowcolor{tablecolor2}}},
	columns/case/.style={string type, column type={l}, column name={Benchmark}},
	columns/polarPower/.style={column name={Polar-Power}},
	columns/polarCurrent/.style={column name={Polar-Current}},
	columns/cartPower/.style={column name={Cartesian-Power}},
	columns/cartCurrent/.style={column name={Cartesian-Current}},
	empty cells with={---} 
	]\OPFtabletimeKNITROc
\end{table}
\begin{table}[H]
	\footnotesize
	\centering
	\caption{Number of iterations (\KNITRO{} 12). \label{tab:OPFvariantssummaryKNITROb-i}}
	\pgfplotstabletypeset[
	every head row/.style={ before row={\toprule},after row=\midrule},
	every last row/.style={ after row=\bottomrule},
	precision=0, fixed zerofill, column type={r},
	columns/Statistic/.style={string type,column type=l},
	every row/.style={
		column type=r,
		dec sep align,
		fixed,
		fixed zerofill,
	},
	columns={case, polarPower, polarCurrent, cartPower, cartCurrent},
	every odd row/.style={before row={\rowcolor{tablecolor1}}},
	every even row/.style={before row={\rowcolor{tablecolor2}}},
	columns/case/.style={string type, column type={l}, column name={Benchmark}},
	columns/polarPower/.style={column name={Polar-Power}},
	columns/polarCurrent/.style={column name={Polar-Current}},
	columns/cartPower/.style={column name={Cartesian-Power}},
	columns/cartCurrent/.style={column name={Cartesian-Current}},
	empty cells with={---} 
	]\OPFtableitersKNITROc
\end{table}
\subsection{Barrier parameter update rule}

An important aspect of the IP algorithm is the strategy used for modifying the barrier parameter $\mu$. Many IP packages usually provide different strategies, while it is not always obvious which strategy is suited best for the given problem at hand.

\subsection*{\IPOPT{}}
The barrier parameter update strategy based on Mehrotra's probing heuristic seems to outperform other strategies by a significant margin, especially when combined with affine type of corrector steps (strategies 6,8, and 9).
The different strategies are listed in Table~\ref{tab:muruleBeltistos} and the performance profiles for all strategies are shown if Figure~\ref{fig:profileBeltMurulePowerPolarT}.

\pgfplotstableread{
Options	1	2	3	4	5	6	7	8	9
mu\_strategy 	monotone	monotone	adaptive	adaptive	adaptive	adaptive	adaptive	adaptive	adaptive
mu\_oracle	{}	{}	quality-function	probing	loqo	probing	probing	probing	probing
fixed\_mu\_oracle	{}	{}	quality-function	probing	loqo	probing	probing	probing	probing
mu\_init	0.1	0.001	{}	{}	{}	{}	{}	{}	{}
corrector\_type	{}	{}	none	none	none	affine	primal-dual	affine	affine
alpha\_for\_y	primal	primal	primal	primal	primal	primal	primal	safer-min-dual-infeas	primal-and-full
alpha\_for\_y\_tol	{}	{}	{}	{}	{}	{}	{}	{}	1
}\muruleSuccess

\begin{table}[ht!]
	\centering
	\caption{Strategies for updating the barrier parameter in \IPOPT{}. \label{tab:muruleBeltistos}}
	\pgfplotstabletypeset[
	font=\scriptsize,
	every head row/.style={ before row=\toprule,after row=\midrule},
	every last row/.style={ after row=\bottomrule},
	precision=0, fixed zerofill, column type={c},
	columns/Options/.style={string type,column type=l},
	columns/1/.style={string type,column type=l},
	columns/2/.style={string type,column type=l},
	columns/3/.style={string type,column type=C},
	columns/4/.style={string type,column type=l},
	columns/5/.style={string type,column type=l},
	columns/6/.style={string type,column type=l},
	columns/7/.style={string type,column type=l},
	columns/8/.style={string type,column type=C},
	columns/9/.style={string type,column type=C},
	every odd row/.style={before row={\rowcolor{tablecolor1}}},
	every even row/.style={before row={\rowcolor{tablecolor2}}},
	empty cells with={---} 
	]\muruleSuccess
\end{table}

\begin{figure}[ht!]
	\centering
	\input{profilesNEW/barmurule_beltistos_polarPower-t.tex}
	\vspace{-0.2cm}
	\caption{Overall time profile for various $\mu$ update strategies, \IPOPT{}, Polar-Power formulation ($\alpha_{max}=7.2$).\label{fig:profileBeltMurulePowerPolarT} }
\end{figure}

\subsection*{\KNITRO{}}
It seems that a Mehrotra predictor-corrector type rule to determine the barrier parameter without safeguards on the corrector step (\texttt{fullmpc}) is the fastest in case of solving the OPF problems, although it is not as robust as other update rules. Although \texttt{fullmpc} is the fastest, in most of the cases the other strategies do not perform much worse and catch up very fast, e.g. Mehrotra predictor-corrector type rule with safeguards on the corrector step (\texttt{dampmpc}) is also doing very well and the additional safeguards increase also the robustness.
When leaving the decision about the parameter up to \KNITRO{}, is selects \texttt{dampmpc} (``\textit{\KNITRO{} changing \texttt{bar\_murule} from AUTO to 4}'').

Other strategies used for modifying the barrier parameter $\mu$ in the barrier algorithm are adaptive rule based on the complementarity gap (\texttt{adaptive}), using a probing (affine-scaling) step to dynamically determine the barrier parameter (\texttt{probing}), or strategy based on minimization of a quality function at each iteration (\texttt{quality}). The success rate and performance profiles are listed in Table \ref{tab:muruleSuccess} and Figure \ref{fig:muruleSuccess}, respectively.

\pgfplotstableread{
	RuleNo Rule polarPower cartPower polarCurrent cartCurrent Total
	3 probing 24 25 25 24  98
	4 dampmpc 24 25 24 25 98
	6 quality 23 25 24 25 97
	5 fullmpc 24 23 25 24 96
	2 adaptive 21 25 24 25 95
}\muruleSuccess

\begin{table}[ht!]
	\centering
	\caption{Number of solved benchmarks for different $\mu$ update rules and OPF formulations (\KNITRO{}). \label{tab:muruleSuccess}}
	\pgfplotstabletypeset[
	every head row/.style={ before row=\toprule,after row=\midrule},
	every last row/.style={ after row=\bottomrule},
	precision=0, fixed zerofill, column type={c},
	columns={Rule, polarPower, cartPower, polarCurrent, cartCurrent, Total },
	columns/Rule/.style={string type,column type=l},
	every odd row/.style={before row={\rowcolor{tablecolor1}}},
	every even row/.style={before row={\rowcolor{tablecolor2}}},
	columns/polarPower/.style={column name={Polar Power}},
	columns/polarCurrent/.style={column name={Polar Current}},
	columns/cartPower/.style={column name={Cartesian Power}},
	columns/cartCurrent/.style={column name={Cartesian Current}}
	]\muruleSuccess
\end{table}

\begin{figure*}[ht!]
	\centering
	\begin{subfigure}[b]{0.5\columnwidth}
		\centering
		\begin{tikzpicture}


\pgfplotstableread{
ALPHA         SOLVER1         SOLVER2         SOLVER3         SOLVER4         SOLVER5
1.0         0.16000         0.04000         0.36000         0.36000         0.08000
1.2         0.52000         0.60000         0.84000         0.60000         0.52000
1.4         0.76000         0.76000         0.92000         0.64000         0.72000
1.5999999999999999         0.76000         0.92000         0.92000         0.68000         0.80000
1.7999999999999998         0.80000         0.92000         0.96000         0.80000         0.84000
1.9999999999999998         0.80000         0.92000         0.96000         0.80000         0.84000
2.1999999999999997         0.80000         0.96000         0.96000         0.84000         0.88000
2.3999999999999995         0.80000         0.96000         0.96000         0.88000         0.88000
2.5999999999999996         0.80000         0.96000         0.96000         0.88000         0.88000
2.8         0.80000         0.96000         0.96000         0.88000         0.88000
2.9999999999999996         0.80000         0.96000         0.96000         0.88000         0.88000
3.1999999999999993         0.80000         0.96000         0.96000         0.88000         0.88000
3.3999999999999995         0.80000         0.96000         0.96000         0.92000         0.88000
3.5999999999999996         0.80000         0.96000         0.96000         0.92000         0.88000
3.7999999999999994         0.80000         0.96000         0.96000         0.92000         0.88000
3.999999999999999         0.80000         0.96000         0.96000         0.92000         0.88000
4.199999999999999         0.80000         0.96000         0.96000         0.96000         0.88000
4.3999999999999995         0.80000         0.96000         0.96000         0.96000         0.88000
4.6         0.80000         0.96000         0.96000         0.96000         0.88000
4.799999999999999         0.80000         0.96000         0.96000         0.96000         0.88000
4.999999999999999         0.84000         0.96000         0.96000         0.96000         0.88000
5.199999999999999         0.84000         0.96000         0.96000         0.96000         0.88000
5.399999999999999         0.84000         0.96000         0.96000         0.96000         0.88000
5.599999999999999         0.84000         0.96000         0.96000         0.96000         0.88000
5.799999999999999         0.84000         0.96000         0.96000         0.96000         0.88000
5.999999999999999         0.84000         0.96000         0.96000         0.96000         0.88000
6.199999999999999         0.84000         0.96000         0.96000         0.96000         0.88000
6.399999999999999         0.84000         0.96000         0.96000         0.96000         0.88000
6.599999999999999         0.84000         0.96000         0.96000         0.96000         0.88000
6.799999999999999         0.84000         0.96000         0.96000         0.96000         0.88000
6.999999999999998         0.84000         0.96000         0.96000         0.96000         0.88000
7.199999999999998         0.84000         0.96000         0.96000         0.96000         0.88000
7.399999999999999         0.84000         0.96000         0.96000         0.96000         0.88000
7.599999999999999         0.84000         0.96000         0.96000         0.96000         0.88000
7.799999999999999         0.84000         0.96000         0.96000         0.96000         0.88000
7.999999999999998         0.84000         0.96000         0.96000         0.96000         0.88000
8.2         0.84000         0.96000         0.96000         0.96000         0.88000
8.399999999999999         0.84000         0.96000         0.96000         0.96000         0.88000
8.599999999999998         0.84000         0.96000         0.96000         0.96000         0.88000
8.799999999999997         0.84000         0.96000         0.96000         0.96000         0.88000
8.999999999999998         0.84000         0.96000         0.96000         0.96000         0.88000
9.199999999999998         0.84000         0.96000         0.96000         0.96000         0.92000
9.399999999999999         0.84000         0.96000         0.96000         0.96000         0.92000
9.599999999999998         0.84000         0.96000         0.96000         0.96000         0.92000
9.799999999999997         0.84000         0.96000         0.96000         0.96000         0.92000
9.999999999999998         0.84000         0.96000         0.96000         0.96000         0.92000
}\datatable

\begin{axis}[name=symb,
width=\columnwidth,
height=4cm,
ymin=-0.02,
ymax=1.02,
enlarge x limits=0.02,
axis lines*=left, 
ymajorgrids, yminorgrids,
xmajorgrids, xminorgrids,
xticklabel style={rotate=0, xshift=-0.0cm, anchor=north, font=\scriptsize},
ytick={0,  0.2,  0.4, 0.6, 0.8, 1},
xmax=6,
yticklabel style={font=\scriptsize},
ylabel style={font=\scriptsize},
xlabel style={font=\scriptsize},
ylabel={$p_m({\alpha})$},
xlabel={$\alpha$},
legend style={at={(1.0,0.6)},legend cell align=left,align=right,draw=white!85!black,font=\scriptsize,legend columns=2}
]

\addplot[color=mycolor4, loosely dashed, ultra thick] table[x=ALPHA, y=SOLVER1] {\datatable};
\addlegendentry{adaptive}

\addplot[color=mycolor2, densely dashed, ultra thick] table[x=ALPHA, y=SOLVER2] {\datatable}; 
\addlegendentry{probing}

\addplot[color=mycolor1, densely dotted, ultra thick] table[x=ALPHA, y=SOLVER3] {\datatable}; 
\addlegendentry{dampmpc}

\addplot[color=mycolor5, loosely dotted, ultra thick] table[x=ALPHA, y=SOLVER4] {\datatable};
\addlegendentry{fullmpc}

\addplot[color=mycolor3,dashdotted, ultra thick] table[x=ALPHA, y=SOLVER5] {\datatable}; 
\addlegendentry{quality}

\end{axis}

\end{tikzpicture}
		\vspace{-0.2cm}
		\caption{Polar-Power ($\alpha_{max}=10$)\label{fig:profileMurulePowerPolarT} }
	\end{subfigure}
	\begin{subfigure}[b]{0.49\columnwidth}
		\centering
		\begin{tikzpicture}


\pgfplotstableread{
ALPHA         SOLVER1         SOLVER2         SOLVER3         SOLVER4         SOLVER5
1.0         0.04000         0.00000         0.28000         0.64000         0.04000
1.2         0.32000         0.28000         0.56000         0.76000         0.20000
1.4         0.60000         0.64000         0.84000         0.84000         0.60000
1.5999999999999999         0.80000         0.88000         0.92000         0.84000         0.92000
1.7999999999999998         0.96000         0.96000         0.92000         0.88000         1.00000
1.9999999999999998         0.96000         0.96000         0.92000         0.92000         1.00000
2.1999999999999997         0.96000         0.96000         0.96000         0.92000         1.00000
2.3999999999999995         0.96000         0.96000         0.96000         0.92000         1.00000
2.5999999999999996         0.96000         1.00000         0.96000         0.92000         1.00000
2.8         0.96000         1.00000         0.96000         0.92000         1.00000
2.9999999999999996         0.96000         1.00000         0.96000         0.92000         1.00000
3.1999999999999993         0.96000         1.00000         0.96000         0.92000         1.00000
3.3999999999999995         0.96000         1.00000         0.96000         0.92000         1.00000
3.5999999999999996         0.96000         1.00000         0.96000         0.92000         1.00000
3.7999999999999994         0.96000         1.00000         0.96000         0.92000         1.00000
3.999999999999999         0.96000         1.00000         0.96000         0.92000         1.00000
4.199999999999999         0.96000         1.00000         0.96000         0.92000         1.00000
4.3999999999999995         0.96000         1.00000         0.96000         0.92000         1.00000
4.6         1.00000         1.00000         0.96000         0.92000         1.00000
4.799999999999999         1.00000         1.00000         0.96000         0.92000         1.00000
4.999999999999999         1.00000         1.00000         0.96000         0.92000         1.00000
5.199999999999999         1.00000         1.00000         0.96000         0.92000         1.00000
5.399999999999999         1.00000         1.00000         0.96000         0.92000         1.00000
5.599999999999999         1.00000         1.00000         0.96000         0.92000         1.00000
5.799999999999999         1.00000         1.00000         0.96000         0.92000         1.00000
5.999999999999999         1.00000         1.00000         0.96000         0.92000         1.00000
6.199999999999999         1.00000         1.00000         0.96000         0.92000         1.00000
6.399999999999999         1.00000         1.00000         0.96000         0.92000         1.00000
6.599999999999999         1.00000         1.00000         0.96000         0.92000         1.00000
6.799999999999999         1.00000         1.00000         0.96000         0.92000         1.00000
6.999999999999998         1.00000         1.00000         0.96000         0.92000         1.00000
7.199999999999998         1.00000         1.00000         0.96000         0.92000         1.00000
7.399999999999999         1.00000         1.00000         0.96000         0.92000         1.00000
7.599999999999999         1.00000         1.00000         0.96000         0.92000         1.00000
7.799999999999999         1.00000         1.00000         0.96000         0.92000         1.00000
7.999999999999998         1.00000         1.00000         0.96000         0.92000         1.00000
8.2         1.00000         1.00000         0.96000         0.92000         1.00000
8.399999999999999         1.00000         1.00000         0.96000         0.92000         1.00000
8.599999999999998         1.00000         1.00000         1.00000         0.92000         1.00000
8.799999999999997         1.00000         1.00000         1.00000         0.92000         1.00000
8.999999999999998         1.00000         1.00000         1.00000         0.92000         1.00000
9.199999999999998         1.00000         1.00000         1.00000         0.92000         1.00000
9.399999999999999         1.00000         1.00000         1.00000         0.92000         1.00000
9.599999999999998         1.00000         1.00000         1.00000         0.92000         1.00000
9.799999999999997         1.00000         1.00000         1.00000         0.92000         1.00000
9.999999999999998         1.00000         1.00000         1.00000         0.92000         1.00000
}\datatable

\begin{axis}[name=symb,
width=\columnwidth,
height=4cm,
ymin=-0.02,
ymax=1.02,
enlarge x limits=0.02,
axis lines*=left, 
ymajorgrids, yminorgrids,
xmajorgrids, xminorgrids,
xticklabel style={rotate=0, xshift=-0.0cm, anchor=north, font=\scriptsize},
ytick={0,  0.2,  0.4, 0.6, 0.8, 1},
xmax=6,
yticklabel style={font=\scriptsize},
ylabel style={font=\scriptsize},
xlabel style={font=\scriptsize},
ylabel={$p_m({\alpha})$},
xlabel={$\alpha$},
legend style={at={(1.0,0.6)},legend cell align=left,align=right,draw=white!85!black,font=\scriptsize,legend columns=2}
]

\addplot[color=mycolor4, loosely dashed, ultra thick] table[x=ALPHA, y=SOLVER1] {\datatable};
\addlegendentry{adaptive}

\addplot[color=mycolor2, densely dashed, ultra thick] table[x=ALPHA, y=SOLVER2] {\datatable}; 
\addlegendentry{probing}

\addplot[color=mycolor1, densely dotted, ultra thick] table[x=ALPHA, y=SOLVER3] {\datatable}; 
\addlegendentry{dampmpc}

\addplot[color=mycolor5, loosely dotted, ultra thick] table[x=ALPHA, y=SOLVER4] {\datatable};
\addlegendentry{fullmpc}

\addplot[color=mycolor3,dashdotted, ultra thick] table[x=ALPHA, y=SOLVER5] {\datatable}; 
\addlegendentry{quality}

\end{axis}

\end{tikzpicture}
		\vspace{-0.2cm}
		\caption{Rect-Power ($\alpha_{max}=10$)\label{fig:profileMurulePowerCartT} }
	\end{subfigure}
	\begin{subfigure}[b]{0.5\columnwidth}
		\centering
		\input{profilesNEW/barmurule_knitro_polarCurrent-t.tex}
		\vspace{-0.2cm}
		\caption{Polar-Current  ($\alpha_{max}=45$)\label{fig:profileMuruleCurrentPolarT} }
	\end{subfigure}
	\begin{subfigure}[b]{0.49\columnwidth}
		\centering
		\input{profilesNEW/barmurule_knitro_cartCurrent-t.tex}
		\vspace{-0.2cm}
		\caption{Rect-Current ($\alpha_{max}=22$)\label{fig:profileMuruleCurrentCartT} }
	\end{subfigure}
	\caption{\label{fig:muruleSuccess} Overall time profile for various $\mu$ update strategies (\KNITRO{}).}
\end{figure*}

\subsection{Linear solvers \label{sec:LSsolver}}
The robustness and performance of the optimization package also depends on direct sparse solver used to compute the search direction in each IP iteration. Linear systems resulting from the IP methods are known to be very ill-conditioned, especially in the few last iterations before convergence and depending on the benchmark case, can be also very large.

Tables~\ref{tab:mipsLSsummary-iters},~\ref{tab:mipsLSsummary-time},~\ref{tab:mipsLSsummary1-iters} and~\ref{tab:mipsLSsummary1-time} demonstrate the difference between build-in Matlab linear solver (\texttt{mldivide}, also known as the '\textbackslash' operator) and \PARDISO{}, which are currently the two linear solvers supported by \MIPS{} framework. The difference is particularly visible in Table~\ref{tab:mipsLSsummary-time} and~\ref{tab:mipsLSsummary1-time} for the large-scale cases, where \PARDISO{} outperforms MATLAB's backslash operator by a factor up to $6$. \PARDISO{} also reduces the number of iterations until convergence for some cases due to more accurate solution and thus better descent search direction provided to the optimizer.
\pgfplotstableread{
case flat-backslash mpc-backslash pfsolve-backslash flat-pardiso mpc-pardiso pfsolve-pardiso
case1951rte             {} 28  28                        {}   28  28
case2383wp              35  32  35                                35    32  35
case2736sp              30  26  26                                30    26  26
case2737sop             28  28  26                                28    28  26
case2746wop             31  26  29                                31    26  29
case2746wp              31  28  28                                31    28  28
case2868rte             {} 32  31                        {}   32  31
case2869pegase          39  127 32                                39    98  32
case3012wp              44  29  30                                44    29  30
case3120sp              43  111 33                                43    117 33
case3375wp              47  30  30                                47    30  30
case6468rte             {} 47  44                        {}   47  44
case6470rte             {} 43  44                        {}   43  44
case6495rte             {} 74  76                        {}   77  76
case6515rte             {} 59  61                        {}   59  61
case9241pegase          45  68  46                                45    68  46
case\_ACTIVSg2000       35  26  29                                35    26  29
case\_ACTIVSg10k        {} 41  81                        {}   41  81
case13659pegase         65  42  54                                64    42  54
case\_ACTIVSg25k        {} 56  76                        {}   56  76
case\_ACTIVSg70k        {} {} {}       {}   {} {}
case21k                 71  53  53                                71    53  53
case42k                 83  64  64                                83    64  64
case99k                 {} {} {}       97    78  78
case193k                {} {} {}       113   92  92
}\mipsLStableiters

\pgfplotstableread{
case flat-backslash mpc-backslash pfsolve-backslash flat-pardiso mpc-pardiso pfsolve-pardiso
case1951rte             {} 3.51    3.71                  {}   4.54    4.75
case2383wp              5   4.71    5.14                          6.49  6.02    6.6
case2736sp              4.68    4.28    4.44                      6.04  5.47    5.98
case2737sop             4.51    4.55    4.57                      5.7   5.83    5.68
case2746wop             4.99    4.39    4.94                      6.63  5.92    6.4
case2746wp              5.44    4.68    4.85                      6.53  6.1 6.26
case2868rte             {} 5.13    5.28                  {}   6.73    6.82
case2869pegase          7.22    24.83   6.11                      9.04  23.44   7.66
case3012wp              7.81    5.16    5.5                       9.86  6.52    7.13
case3120sp              7.45    21.61   6.36                      9.4   28.23   7.6
case3375wp              8.77    5.72    5.91                      11.23 7.2 7.37
case6468rte             {} 13.95   14.18                 {}   17.25   16.57
case6470rte             {} 13.43   14.16                 {}   16.84   17.31
case6495rte             {} 23.34   24.29                 {}   31.68   30.8
case6515rte             {} 18.42   19.36                 {}   23.51   24.71
case9241pegase          23.84   36.05   25.9                      30.67 47.08   31.91
case\_ACTIVSg2000       4.9 3.84    4.38                          6.41  5.04    5.7
case\_ACTIVSg10k        {} 20.26   45.32                 {}   28.07   60.06
case13659pegase         36.51   24.42   30.53                     51.19 33.55   43.54
case\_ACTIVSg25k        {} 72.43   108.18                {}   93.24   138.79
case\_ACTIVSg70k        {} {} {}       {}   {} {}
case21k                 291.87  238.12  241.2                     135.4 97.4    98
case42k                 2264.12 1925.02 1957.95                   390.42    289.74  303.38
case99k                 {} {} {}       1329.69   1082.13 1081.72
case193k                {} {} {}       3703.2    2953.02 2982.83
}\mipsLStabletime

\begin{table*}[ht!]
    \footnotesize
	\centering
	\caption{Performance of \MIPS{} with various linear solvers - number of iterations. \label{tab:mipsLSsummary-iters}}
	\pgfplotstabletypeset[
	every head row/.style={ before row={%
            \toprule
            & \multicolumn{2}{c}{Flat} & \multicolumn{2}{c}{MPC} & \multicolumn{2}{c}{PF} \\
            \cmidrule(l{4pt}r{2pt}){2-3} \cmidrule(l{4pt}r{2pt}){4-5} \cmidrule(l{4pt}r{2pt}){6-7}
                         },after row=\midrule},
	every last row/.style={ after row=\bottomrule},
	precision=0, fixed zerofill, column type={r},
	every row/.style={
		column type=r,
		dec sep align,
		fixed,
		fixed zerofill,
	},
	every odd row/.style={before row={\rowcolor{tablecolor1}}},
    every even row/.style={before row={\rowcolor{tablecolor2}}},
	columns={case, flat-pardiso, flat-backslash, mpc-pardiso, mpc-backslash, pfsolve-pardiso, pfsolve-backslash},
	columns/case/.style={string type, column type={l}, column name={Benchmark}},
	columns/flat-backslash/.style={column name={Backslash}},
	columns/flat-pardiso/.style={column name={\PARDISO{}}},
	columns/mpc-backslash/.style={column name={Backslash}},
	columns/mpc-pardiso/.style={column name={\PARDISO{}}},
	columns/pfsolve-backslash/.style={column name={Backslash}},
	columns/pfsolve-pardiso/.style={column name={\PARDISO{}}},
	empty cells with={---} 
	]\mipsLStableiters
\end{table*}

\begin{table*}[ht!]
    \footnotesize
	\centering
	\caption{Performance of \MIPS{} with various linear solvers - overall time (s). \label{tab:mipsLSsummary-time}}
	\pgfplotstabletypeset[
	every head row/.style={ before row={%
            \toprule
            & \multicolumn{2}{c}{Flat} & \multicolumn{2}{c}{MPC} & \multicolumn{2}{c}{PF} \\
            \cmidrule(l{4pt}r{2pt}){2-3} \cmidrule(l{4pt}r{2pt}){4-5} \cmidrule(l{4pt}r{2pt}){6-7}
                         },after row=\midrule},
	every last row/.style={ after row=\bottomrule},
	precision=2, fixed zerofill, column type={r},
	every row/.style={
		column type=r,
		dec sep align,
		fixed,
		fixed zerofill,
	},
	every odd row/.style={before row={\rowcolor{tablecolor1}}},
    every even row/.style={before row={\rowcolor{tablecolor2}}},
	columns={case, flat-pardiso, flat-backslash, mpc-pardiso, mpc-backslash, pfsolve-pardiso, pfsolve-backslash},
	columns/case/.style={string type, column type={l}, column name={Benchmark}},
	columns/flat-backslash/.style={column name={Backslash}},
	columns/flat-pardiso/.style={column name={\PARDISO{}}},
	columns/mpc-backslash/.style={column name={Backslash}},
	columns/mpc-pardiso/.style={column name={\PARDISO{}}},
	columns/pfsolve-backslash/.style={column name={Backslash}},
	columns/pfsolve-pardiso/.style={column name={\PARDISO{}}},
	empty cells with={---} 
	]\mipsLStabletime
\end{table*}
\pgfplotstableread{
case polarPower-backslash cartPower-backslash polarCurrent-backslash cartCurrent-backslash     polarPower-pardiso cartPower-pardiso polarCurrent-pardiso cartCurrent-pardiso
case1951rte             28  62  32  106                                   28    61  32  164
case2383wp              35  40  35  46                                    35    40  35  46
case2736sp              26  32  27  46                                    26    32  27  46
case2737sop             26  30  26  30                                    26    30  26  30
case2746wop             29  29  29  31                                    29    29  29  31
case2746wp              28  37  30  37                                    28    37  30  37
case2868rte             31  41  {} 428                           31    41  63  {}
case2869pegase          32  47  36  {}                           32    47  36  {}
case3012wp              30  36  29  37                                    30    36  29  37
case3120sp              33  38  34  36                                    33    38  34  36
case3375wp              30  37  33  40                                    30    37  33  40
case6468rte             44  {} 66  {}                   44    {} 54  {}
case6470rte             44  178 94  {}                           44    172 95  {}
case6495rte             76  {} 131 165                           76    153 116 205
case6515rte             61  160 59  {}                           61    165 67  152
case9241pegase          46  {} 51  {}                   46    218 51  {}
case\_ACTIVSg2000       29  31  52  48                                    29    31  52  48
case\_ACTIVSg10k        81  66  108 83                                    81    67  94  84
case13659pegase         54  {} {} {}           54    {} {} {}
case\_ACTIVSg25k        76  72  78  103                                   76    72  78  103
case\_ACTIVSg70k        {} 96  96  152                           {}   96  92  152
case21k                 53  61  53  61                                    53    61  {} {}
case42k                 64  73  64  74                                    64    73  92  {}
case99k                 {} {} {} {}   78    85  {} {}
case193k                {} {} {} {}   92    101 {} {}
}\mipsLStableitersVariants

\pgfplotstableread{
case polarPower-backslash cartPower-backslash polarCurrent-backslash cartCurrent-backslash     polarPower-pardiso cartPower-pardiso polarCurrent-pardiso cartCurrent-pardiso
case1951rte             3.71    8.47    4.1 14.53                              4.75 10.46   5.69    29.97
case2383wp              5.14    5.51    5.17    6.87                           6.6  7.38    6.55    8.63
case2736sp              4.44    5.23    4.55    7.7                            5.98 6.82    6.05    9.62
case2737sop             4.57    5.55    4.45    5.44                           5.68 6.28    5.91    6.32
case2746wop             4.94    4.83    5.03    5.23                           6.4  6.36    6.84    6.87
case2746wp              4.85    5.85    5.1 6.25                               6.26 7.67    6.66    7.81
case2868rte             5.28    6.76    {} 96.18                      6.82 8.66    13.79   {}
case2869pegase          6.11    9.26    6.63    {}                    7.66 12.18   8.39    {}
case3012wp              5.5 6.38    5.41    6.67                               7.13 8.44    7.12    8.38
case3120sp              6.36    6.94    6.14    6.31                           7.6  8.79    7.7 8.16
case3375wp              5.91    7.17    6.42    7.46                           7.37 8.82    8.11    9.85
case6468rte             14.18   {} 21.91   {}                16.57    {} 19.93   {}
case6470rte             14.16   67.94   33.07   {}                    17.31    81.82   42.32   {}
case6495rte             24.29   {} 42.42   58.21                      30.8 72.56   49.91   94.68
case6515rte             19.36   55.22   19.44   {}                    24.71    79.67   29.13   68.41
case9241pegase          25.9    {} 26.51   {}                31.91    185.65  34.71   {}
case\_ACTIVSg2000       4.38    4.71    8.02    7.42                           5.7  6.29    10.19   9.28
case\_ACTIVSg10k        45.32   34.48   62.39   48.05                          60.06    47.41   69.08   62.79
case13659pegase         30.53   {} {} {}            43.54    {} {} {}
case\_ACTIVSg25k        108.18  92.13   104.14  150.27                         138.79   120.95  138.73  192.83
case\_ACTIVSg70k        {} 404.02  423.96  728.76                     {}  500.27  486.22  877.68
case21k                 241.2   253.8   230.33  250.33                         98   110.55  {} {}
case42k                 1957.95 2023.77 1966.92 2090.61                        303.38   333.91  449.81  {}
case99k                 {} {} {} {}        1081.72  1174.49 {} {}
case193k                {} {} {} {}        2982.83  3235.83 {} {}
}\mipsLStabletimeVariants

\begin{table*}[ht!]
    \footnotesize
	\centering
	\caption{Performance of \MIPS{} with various linear solvers - number of iterations. \label{tab:mipsLSsummary1-iters}}
	\pgfplotstabletypeset[
	every head row/.style={ before row={%
            \toprule
            & \multicolumn{2}{c}{Polar-Power} & \multicolumn{2}{c}{Polar-Current} & \multicolumn{2}{c}{Cartesian-Power} & \multicolumn{2}{c}{Cartesian-Current} \\
            \cmidrule(l{4pt}r{2pt}){2-3} \cmidrule(l{4pt}r{2pt}){4-5} \cmidrule(l{4pt}r{2pt}){6-7} \cmidrule(l{4pt}r{2pt}){8-9}
                         },after row=\midrule},
	every last row/.style={ after row=\bottomrule},
	precision=0, fixed zerofill, column type={r},
	every row/.style={
		column type=r,
		dec sep align,
		fixed,
		fixed zerofill,
	},
	every odd row/.style={before row={\rowcolor{tablecolor1}}},
    every even row/.style={before row={\rowcolor{tablecolor2}}},
	columns={case, polarPower-pardiso, polarPower-backslash, polarCurrent-pardiso, polarCurrent-backslash, cartPower-pardiso, cartPower-backslash, cartCurrent-pardiso, cartCurrent-backslash},
	columns/case/.style={string type, column type={l}, column name={Benchmark}},
	columns/polarPower-backslash/.style={column name={Backslash}},
	columns/polarPower-pardiso/.style={column name={\PARDISO{}}},
	columns/polarCurrent-backslash/.style={column name={Backslash}},
	columns/polarCurrent-pardiso/.style={column name={\PARDISO{}}},
	columns/cartPower-backslash/.style={column name={Backslash}},
	columns/cartPower-pardiso/.style={column name={\PARDISO{}}},
	columns/cartCurrent-backslash/.style={column name={Backslash}},
	columns/cartCurrent-pardiso/.style={column name={\PARDISO{}}},
	empty cells with={---} 
	]\mipsLStableitersVariants
\end{table*}

\begin{table*}[ht!]
    \footnotesize
	\centering
	\caption{Performance of \MIPS{} with various linear solvers - overall time (s). \label{tab:mipsLSsummary1-time}}
	\pgfplotstabletypeset[
	every head row/.style={ before row={%
            \toprule
            & \multicolumn{2}{c}{Polar-Power} & \multicolumn{2}{c}{Polar-Current} & \multicolumn{2}{c}{Cartesian-Power} & \multicolumn{2}{c}{Cartesian-Current} \\
            \cmidrule(l{4pt}r{2pt}){2-3} \cmidrule(l{4pt}r{2pt}){4-5} \cmidrule(l{4pt}r{2pt}){6-7} \cmidrule(l{4pt}r{2pt}){8-9}
                         },after row=\midrule},
	every last row/.style={ after row=\bottomrule},
	precision=2, fixed zerofill, column type={r},
	every row/.style={
		column type=r,
		dec sep align,
		fixed,
		fixed zerofill,
	},
	every odd row/.style={before row={\rowcolor{tablecolor1}}},
    every even row/.style={before row={\rowcolor{tablecolor2}}},
	columns={case, polarPower-pardiso, polarPower-backslash, polarCurrent-pardiso, polarCurrent-backslash, cartPower-pardiso, cartPower-backslash, cartCurrent-pardiso, cartCurrent-backslash},
	columns/case/.style={string type, column type={l}, column name={Benchmark}},
	columns/polarPower-backslash/.style={column name={Backslash}},
	columns/polarPower-pardiso/.style={column name={\PARDISO{}}},
	columns/polarCurrent-backslash/.style={column name={Backslash}},
	columns/polarCurrent-pardiso/.style={column name={\PARDISO{}}},
	columns/cartPower-backslash/.style={column name={Backslash}},
	columns/cartPower-pardiso/.style={column name={\PARDISO{}}},
	columns/cartCurrent-backslash/.style={column name={Backslash}},
	columns/cartCurrent-pardiso/.style={column name={\PARDISO{}}},
	empty cells with={---} 
	]\mipsLStabletimeVariants
\end{table*}

Tables~\ref{tab:ipoptLSsummaryIters}, \ref{tab:ipoptLSsummaryTime}, \ref{tab:ipoptLSsummary1-iters}, and \ref{tab:ipoptLSsummary1-time} demonstrate performance of the \IPOPT{} with two different linear solvers, \PARDISO{} and HSL MA57. MA57 is a robust solver but as the problem size increases it requires significantly higher computational resources than \PARDISO{}. On the other hand, computational times using \IPOPT{} with \PARDISO{} remain feasible also for the large-scale networks and it is thus possible to solve more benchmarks.
\pgfplotstableread{
case flat-pardiso mpc-pardiso pfsolve-pardiso flat-hsl mpc-hsl pfsolve-hsl
case1951rte             238 43  53                             177  36  34
case2383wp              177 51  43                             36   54  39
case2736sp              309 30  22                             30   17  23
case2737sop             34  21  27                             31   20  28
case2746wop             45  18  28                             30   15  29
case2746wp              37  17  26                             33   16  28
case2868rte             136 41  50                             160  41  39
case2869pegase          46  47  72                             36   31  29
case3012wp              41  192 40                             43   49  48
case3120sp              41  44  37                             40   40  38
case3375wp              35  47  86                             34   36  43
case6468rte             414 162 146                            169  38  38
case6470rte             486 284 393                            206  70  68
case6495rte             68  73  73                             224  56  60
case6515rte             274 58  60                             142  62  57
case9241pegase          335 178 65                             41   54  45
case\_ACTIVSg2000       212 26  29                             212  26  29
case\_ACTIVSg10k        287 49  345                            99   31  35
case13659pegase         {} {} 244            225  200 246
case\_ACTIVSg25k        58  62  377                            56   51  50
case\_ACTIVSg70k        102 {} 104                    72   62  65
case21k                 70  145 114                            80   299 {}
case42k                 70  73  180                            181  {} {}
case99k                 91  95  94                             {}  {} {}
case193k                {} {} {}    {}  {} {}
}\ipoptLStableiters

\pgfplotstableread{
case flat-pardiso mpc-pardiso pfsolve-pardiso flat-hsl mpc-hsl pfsolve-hsl
case1951rte             17.380  4.200   5.170                                12.09  3.31    3.4
case2383wp              29.230  9.030   8.080                                3.58   4.78    4.11
case2736sp              30.310  4.210   3.630                                3.73   2.47    3.24
case2737sop             4.830   3.380   4.160                                3.54   2.65    3.64
case2746wop             5.870   3.140   4.170                                3.46   2.32    3.69
case2746wp              4.760   3.020   3.930                                3.73   2.39    3.97
case2868rte             15.210  5.970   6.160                                16.01  4.71    4.77
case2869pegase          9.660   10.220  15.940                               4.68   3.91    3.89
case3012wp              9.110   36.550  9.370                                5.18   5.23    5.57
case3120sp              9.510   9.750   8.650                                4.78   4.66    4.86
case3375wp              8.330   11.010  18.880                               4.27   4.47    5.4
case6468rte             80.930  28.620  27.030                               24.75  6.29    6.65
case6470rte             104.450 62.960  85.710                               33.64  12.18   12.51
case6495rte             15.170  14.320  13.620                               35.49  9.1 10.17
case6515rte             49.620  11.540  11.770                               23.15  10.59   9.86
case9241pegase          235.390 118.080 52.520                               15.05  20.42   17.76
case\_ACTIVSg2000       41.080  3.530   3.990                                44.24  3.44    4
case\_ACTIVSg10k        125.450 19.320  121.120                              38.47  11.58   13.99
case13659pegase         {} {} 180.440              83.52  71.93   86.76
case\_ACTIVSg25k        57.600  55.030  362.330                              54.34  46.75   50.06
case\_ACTIVSg70k        332.800 {} 345.750                      167.8  144.37  152.56
case21k                 161.490 323.110 330.290                              120.91 598.92  {}
case42k                 629.440 850.520 2024.700                             7918.16    {} {}
case99k                 4192.140    4983.650    3852.760                     {}    {} {}
case193k                {} {} {}      {}    {} {}
}\ipoptLStabletime

\begin{table*}[ht!]
    \footnotesize
	\centering
	\caption{Performance of \IPOPT{} with various linear solvers - number of iterations. \label{tab:ipoptLSsummaryIters}}
	\pgfplotstabletypeset[
	every head row/.style={ before row={%
            \toprule
            & \multicolumn{2}{c}{Flat} & \multicolumn{2}{c}{MPC} & \multicolumn{2}{c}{PF} \\
            \cmidrule(l{4pt}r{2pt}){2-3} \cmidrule(l{4pt}r{2pt}){4-5} \cmidrule(l{4pt}r{2pt}){6-7}
                         },after row=\midrule},
	every last row/.style={ after row=\bottomrule},
	precision=0, fixed zerofill, column type={r},
	every row/.style={
		column type=r,
		dec sep align,
		fixed,
		fixed zerofill,
	},
	every odd row/.style={before row={\rowcolor{tablecolor1}}},
    every even row/.style={before row={\rowcolor{tablecolor2}}},
    columns={case, flat-pardiso, flat-hsl, mpc-pardiso, mpc-hsl, pfsolve-pardiso, pfsolve-hsl},
	columns/case/.style={string type, column type={l}, column name={Benchmark}},
	columns/flat-hsl/.style={column name={MA57}},
	columns/flat-pardiso/.style={column name={\PARDISO{}}},
	columns/mpc-hsl/.style={column name={MA57}},
	columns/mpc-pardiso/.style={column name={\PARDISO{}}},
	columns/pfsolve-hsl/.style={column name={MA57}},
	columns/pfsolve-pardiso/.style={column name={\PARDISO{}}},
	empty cells with={---} 
	]\ipoptLStableiters
\end{table*}

\begin{table*}[ht!]
    \footnotesize
	\centering
	\caption{Performance of \IPOPT{} with various linear solvers - overall time (s). \label{tab:ipoptLSsummaryTime}}
	\pgfplotstabletypeset[
	every head row/.style={ before row={%
            \toprule
            & \multicolumn{2}{c}{Flat} & \multicolumn{2}{c}{MPC} & \multicolumn{2}{c}{PF} \\
            \cmidrule(l{4pt}r{2pt}){2-3} \cmidrule(l{4pt}r{2pt}){4-5} \cmidrule(l{4pt}r{2pt}){6-7}
                         },after row=\midrule},
	every last row/.style={ after row=\bottomrule},
	precision=2, fixed zerofill, column type={r},
	every row/.style={
		column type=r,
		dec sep align,
		fixed,
		fixed zerofill,
	},
	every odd row/.style={before row={\rowcolor{tablecolor1}}},
    every even row/.style={before row={\rowcolor{tablecolor2}}},
    columns={case, flat-pardiso, flat-hsl, mpc-pardiso, mpc-hsl, pfsolve-pardiso, pfsolve-hsl},
	columns/case/.style={string type, column type={l}, column name={Benchmark}},
	columns/flat-hsl/.style={column name={MA57}},
	columns/flat-pardiso/.style={column name={\PARDISO{}}},
	columns/mpc-hsl/.style={column name={MA57}},
	columns/mpc-pardiso/.style={column name={\PARDISO{}}},
	columns/pfsolve-hsl/.style={column name={MA57}},
	columns/pfsolve-pardiso/.style={column name={\PARDISO{}}},
	empty cells with={---} 
	]\ipoptLStabletime
\end{table*}
\pgfplotstableread{
case polarPower-pardiso cartPower-pardiso polarCurrent-pardiso cartCurrent-pardiso polarPower-hsl cartPower-hsl polarCurrent-hsl cartCurrent-hsl
case1951rte             53  33  45  48                                    34    33  37  37
case2383wp              43  41  45  45                                    39    41  38  46
case2736sp              22  28  25  33                                    23    28  23  33
case2737sop             27  26  31  29                                    28    26  26  29
case2746wop             28  27  29  48                                    29    27  28  32
case2746wp              26  28  27  27                                    28    28  26  27
case2868rte             50  50  83  41                                    39    50  39  41
case2869pegase          72  33  45  63                                    29    33  34  37
case3012wp              40  38  166 93                                    48    38  48  40
case3120sp              37  44  125 {}                           38    44  40  46
case3375wp              86  37  81  89                                    43    37  40  45
case6468rte             146 41  152 47                                    38    41  39  41
case6470rte             393 85  73  193                                   68    85  67  92
case6495rte             73  64  60  100                                   60    64  64  71
case6515rte             60  60  243 220                                   57    60  72  89
case9241pegase          65  149 122 {}                           45    41  76  47
case\_ACTIVSg2000       29  33  30  38                                    29    33  30  36
case\_ACTIVSg10k        345 36  117 68                                    35    36  36  38
case13659pegase         244 353 226 216                                   246   343 282 374
case\_ACTIVSg25k        377 54  67  {}                           50    155 51  76
case\_ACTIVSg70k        104 70  78  449                                   65    92  64  74
case21k                 114 68  {} 358                           {}   {} 388 {}
case42k                 180 76  388 406                                   {}   {} 343 {}
case99k                 94  94  {} {}                   {}   {} {} {}
case193k                {} {} {} {}   {}   {} {} {}
}\ipoptLStableitersVariants

\pgfplotstableread{
case polarPower-pardiso cartPower-pardiso polarCurrent-pardiso cartCurrent-pardiso polarPower-hsl cartPower-hsl polarCurrent-hsl cartCurrent-hsl
case1951rte             5.17    3.85    6.16    9.34                         3.4    3.44    3.5 4.04
case2383wp              8.08    9.98    6.17    11.18                        4.11   4.36    3.9 5.01
case2736sp              3.63    4.15    4.21    6.14                         3.24   3.7 3.1 4.21
case2737sop             4.16    4.08    4.6 4.87                             3.64   3.55    3.46    3.81
case2746wop             4.17    4.21    4.48    16.36                        3.69   3.81    3.56    4.06
case2746wp              3.93    4.38    4.13    4.14                         3.97   3.75    3.42    3.7
case2868rte             6.16    7.32    11.31   9.7                          4.77   6.59    4.63    5.14
case2869pegase          15.94   10.31   10.98   22.18                        3.89   4.59    4.12    4.97
case3012wp              9.37    11.37   35.96   35.49                        5.57   5.55    5.38    4.91
case3120sp              8.65    12.99   22.4    {}                  4.86   5.65    4.86    5.6
case3375wp              18.88   11.42   25.87   38.24                        5.4    5   5.04    5.7
case6468rte             27.03   10.74   30.63   18.31                        6.65   8.44    6.18    7.96
case6470rte             85.71   19.65   22.06   159.11                       12.51  19.51   11.06   19.73
case6495rte             13.62   13.78   14.41   45.24                        10.17  12.51   10.05   13.18
case6515rte             11.77   13.38   58.29   171.48                       9.86   12.27   11.64   19.43
case9241pegase          52.52   221.52  126.28  {}                  17.76  16.24   30.79   16.8
case\_ACTIVSg2000       3.99    4.81    4.35    8.28                         4  4.45    4.3 4.67
case\_ACTIVSg10k        121.12  16.68   45.53   105                          13.99  14.39   13.26   15.3
case13659pegase         180.44  369.26  189.62  277.08                       86.76  144.02  66.78   137.24
case\_ACTIVSg25k        362.33  52.14   96.87   {}                  50.06  191.89  37.01   117.16
case\_ACTIVSg70k        345.75  214.18  528.73  8500.76                      152.56 339.94  140.52  265.47
case21k                 330.29  227.3   {} 2839.95                  {}    {} 825.15  {}
case42k                 2024.7  1055.39 10250.8 14050.54                     {}    {} 8745.17 {}
case99k                 3852.76 6133.73 {} {}              {}    {} {} {}
case193k                {} {} {} {}      {}    {} {} {}
}\ipoptLStabletimeVariants

\begin{table*}[ht!]
    \footnotesize
	\centering
	\caption{Performance of \IPOPT{} with various linear solvers - number of iterations. \label{tab:ipoptLSsummary1-iters}}
	\pgfplotstabletypeset[
	every head row/.style={ before row={%
            \toprule
            & \multicolumn{2}{c}{Polar-Power} & \multicolumn{2}{c}{Polar-Current} & \multicolumn{2}{c}{Cartesian-Power} & \multicolumn{2}{c}{Cartesian-Current} \\
            \cmidrule(l{4pt}r{2pt}){2-3} \cmidrule(l{4pt}r{2pt}){4-5} \cmidrule(l{4pt}r{2pt}){6-7} \cmidrule(l{4pt}r{2pt}){8-9}
                         },after row=\midrule},
	every last row/.style={ after row=\bottomrule},
	precision=0, fixed zerofill, column type={r},
	every row/.style={
		column type=r,
		dec sep align,
		fixed,
		fixed zerofill,
	},
	every odd row/.style={before row={\rowcolor{tablecolor1}}},
    every even row/.style={before row={\rowcolor{tablecolor2}}},
	columns={case, polarPower-pardiso, polarPower-hsl, polarCurrent-pardiso, polarCurrent-hsl, cartPower-pardiso, cartPower-hsl, cartCurrent-pardiso, cartCurrent-hsl},
	columns/case/.style={string type, column type={l}, column name={Benchmark}},
	columns/polarPower-hsl/.style={column name={MA57}},
	columns/polarPower-pardiso/.style={column name={\PARDISO{}}},
	columns/polarCurrent-hsl/.style={column name={MA57}},
	columns/polarCurrent-pardiso/.style={column name={\PARDISO{}}},
	columns/cartPower-hsl/.style={column name={MA57}},
	columns/cartPower-pardiso/.style={column name={\PARDISO{}}},
	columns/cartCurrent-hsl/.style={column name={MA57}},
	columns/cartCurrent-pardiso/.style={column name={\PARDISO{}}},
	empty cells with={---} 
	]\ipoptLStableitersVariants
\end{table*}

\begin{table*}[ht!]
    \footnotesize
	\centering
	\caption{Performance of \IPOPT{} with various linear solvers - overall time (s). \label{tab:ipoptLSsummary1-time}}
	\pgfplotstabletypeset[
	every head row/.style={ before row={%
            \toprule
            & \multicolumn{2}{c}{Polar-Power} & \multicolumn{2}{c}{Polar-Current} & \multicolumn{2}{c}{Cartesian-Power} & \multicolumn{2}{c}{Cartesian-Current} \\
            \cmidrule(l{4pt}r{2pt}){2-3} \cmidrule(l{4pt}r{2pt}){4-5} \cmidrule(l{4pt}r{2pt}){6-7} \cmidrule(l{4pt}r{2pt}){8-9}
                         },after row=\midrule},
	every last row/.style={ after row=\bottomrule},
	precision=2, fixed zerofill, column type={r},
	every row/.style={
		column type=r,
		dec sep align,
		fixed,
		fixed zerofill,
	},
	every odd row/.style={before row={\rowcolor{tablecolor1}}},
    every even row/.style={before row={\rowcolor{tablecolor2}}},
	columns={case, polarPower-pardiso, polarPower-hsl, polarCurrent-pardiso, polarCurrent-hsl, cartPower-pardiso, cartPower-hsl, cartCurrent-pardiso, cartCurrent-hsl},
	columns/case/.style={string type, column type={l}, column name={Benchmark}},
	columns/polarPower-hsl/.style={column name={MA57}},
	columns/polarPower-pardiso/.style={column name={\PARDISO{}}},
	columns/polarCurrent-hsl/.style={column name={MA57}},
	columns/polarCurrent-pardiso/.style={column name={\PARDISO{}}},
	columns/cartPower-hsl/.style={column name={MA57}},
	columns/cartPower-pardiso/.style={column name={\PARDISO{}}},
	columns/cartCurrent-hsl/.style={column name={MA57}},
	columns/cartCurrent-pardiso/.style={column name={\PARDISO{}}},
	empty cells with={---} 
	]\ipoptLStabletimeVariants
\end{table*}

\subsection{Optimization frameworks performance}

In this section we evaluate high-performance  nonlinear  optimizers that are supported by \MATPOWER{}. These include \BELTISTOS{}, \IPOPT{}, \FMINCON{} 2018b, \KNITRO{} 11,  and  \MATPOWER{}'s  included  default  solver, \MIPS{}. 
The performance profiles all optimization packages for overall timing are shown in Figures~\ref{fig:profilePowerPolarT}--\ref{fig:profileCurrentCartT}.  The performance profiles clearly indicate that \BELTISTOSOPF{} and \KNITRO{} optimizers converged to the optimal solution for majority of the benchmark cases with the best performance among the optimizers considered in this study. On the other hand, the performance of \FMINCON{} was consistently inferior. In terms of the memory efficiency, the best results are achieved for \MIPS{}, while \FMINCON{} requires up to 2.5 times more memory, as shown in Figure \ref{fig:profilePowerPolarM}.
\pgfplotstableread{
case	MIPSsc	MIPSscpardiso	IPOPTpardiso	IPOPTma57	BELTISTOS	FMINCON	KNITRO12
case1951rte	3.71	4.75	5.17	3.4	3.3	9.29	3.08
case2383wp	5.14	6.6	8.08	4.11	4.56	13.25	3.9
case2736sp	4.44	5.98	3.63	3.24	2.95	10.7	3.28
case2737sop	4.57	5.68	4.16	3.64	2.99	7.59	3.39
case2746wop	4.94	6.4	4.17	3.69	3.38	9.28	3.4
case2746wp	4.85	6.26	3.93	3.97	3.7	10.47	3.74
case2868rte	5.28	6.82	6.16	4.77	4.77	22.09	5.34
case2869pegase	6.11	7.66	15.94	3.89	4	8.95	3.95
case3012wp	5.5	7.13	9.37	5.57	5.99	{}	3.93
case3120sp	6.36	7.6	8.65	4.86	5.65	14	4.35
case3375wp	5.91	7.37	18.88	5.4	5.45	{}	4.12
case6468rte	14.18	16.57	27.03	6.65	7.45	{}	6.6
case6470rte	14.16	17.31	85.71	12.51	8.79	110.68	7.63
case6495rte	24.29	30.8	13.62	10.17	11.22	25.56	8.52
case6515rte	19.36	24.71	11.77	9.86	15.63	35.61	8.36
case9241pegase	25.9	31.91	52.52	17.76	12.68	85.86	28.36
case\_ACTIVSg2000	4.38	5.7	3.99	4	3.06	9.44	3.29
case\_ACTIVSg10k	45.32	60.06	121.12	13.99	10.28	38.61	10.78
case13659pegase	30.53	43.54	180.44	86.76	25.51	3975.64	134.32
ACTIVSg25k	108.18	138.79	362.33	50.06	43.04	183.57	45.79
ACTIVSg70k	{}	{}	345.75	152.56	152.38	644.3	187.65
case21k	241.2	98	330.29	{}	65.27	1548.28	47.05
case42k	1957.95	303.38	2024.7	{}	220.83	{}	161.59
case99k	{}	1081.72	3852.76	{}	1231.9	17325.98	753.07
case193k	{}	2982.83	{}	{}	4360.57	{}	1889.51
}\OPFbenchmarksTimePolPow

\pgfplotstableread{
	case	MIPSsc	MIPSscpardiso	IPOPTpardiso	IPOPTma57	BELTISTOS	FMINCON	KNITRO12
case1951rte	8.47	10.46	3.85	3.44	4.19	8.29	3.28
case2383wp	5.51	7.38	9.98	4.36	6.21	53.41	4.26
case2736sp	5.23	6.82	4.15	3.7	4.5	12.01	3.45
case2737sop	5.55	6.28	4.08	3.55	4.93	9.18	3.66
case2746wop	4.83	6.36	4.21	3.81	3.9	8.5	3.82
case2746wp	5.85	7.67	4.38	3.75	5.58	11.75	3.63
case2868rte	6.76	8.66	7.32	6.59	4.73	23.72	4.65
case2869pegase	9.26	12.18	10.31	4.59	6.6	13.56	4.49
case3012wp	6.38	8.44	11.37	5.55	7.53	16.63	4.7
case3120sp	6.94	8.79	12.99	5.65	5.01	17.38	4.81
case3375wp	7.17	8.82	11.42	5	6.67	22.84	4.77
case6468rte	{}	{}	10.74	8.44	9.27	34.39	8.05
case6470rte	67.94	81.82	19.65	19.51	13.38	69.07	8.86
case6495rte	{}	72.56	13.78	12.51	13.7	29.01	10.2
case6515rte	55.22	79.67	13.38	12.27	13.3	61.98	10.16
case9241pegase	{}	185.65	221.52	16.24	14.4	73.96	15.94
case\_ACTIVSg2000	4.71	6.29	4.81	4.45	4.32	11.11	3.6
case\_ACTIVSg10k	34.48	47.41	16.68	14.39	12.27	58.12	13.82
case13659pegase	{}	{}	369.26	144.02	40.24	{}	100.02
ACTIVSg25k	92.13	120.95	52.14	191.89	53.71	388.48	55.72
ACTIVSg70k	404.02	500.27	214.18	339.94	209.02	1619.39	218.07
case21k	253.8	110.55	227.3	{}	68	640.1	64.65
case42k	2023.77	333.91	1055.39	{}	216.82	1421.21	175.47
case99k	{}	1174.49	6133.73	{}	1226.06	6061.47	892.22
case193k	{}	3235.83	{}	{}	3662.7	16911.01	1899.8
}\OPFbenchmarksTimeCarPow

\pgfplotstableread{
	case	MIPSsc	MIPSscpardiso	IPOPTpardiso	IPOPTma57	BELTISTOS	FMINCON	KNITRO12
case1951rte	4.1	5.69	6.16	3.5	2.75	72.99	4.93
case2383wp	5.17	6.55	6.17	3.9	4.48	16.34	3.79
case2736sp	4.55	6.05	4.21	3.1	2.73	10.81	3.31
case2737sop	4.45	5.91	4.6	3.46	2.97	8.75	3.44
case2746wop	5.03	6.84	4.48	3.56	4.18	9.93	3.31
case2746wp	5.1	6.66	4.13	3.42	3.13	9.64	3.49
case2868rte	{}	13.79	11.31	4.63	4	310.19	5.49
case2869pegase	6.63	8.39	10.98	4.12	5.25	8.7	4.01
case3012wp	5.41	7.12	35.96	5.38	6.53	{}	4.04
case3120sp	6.14	7.7	22.4	4.86	4.89	12.45	4.17
case3375wp	6.42	8.11	25.87	5.04	5.8	109.52	4.14
case6468rte	21.91	19.93	30.63	6.18	7.07	81.76	7.37
case6470rte	33.07	42.32	22.06	11.06	7.91	424.97	10.35
case6495rte	42.42	49.91	14.41	10.05	11.34	164.49	9.32
case6515rte	19.44	29.13	58.29	11.64	9.22	118.64	14.25
case9241pegase	26.51	34.71	126.28	30.79	15.48	204.36	18.6
case\_ACTIVSg2000	8.02	10.19	4.35	4.3	3.67	9.98	3.46
case\_ACTIVSg10k	62.39	69.08	45.53	13.26	14.29	46.11	10.84
case13659pegase	{}	{}	189.62	66.78	17.61	4894.19	896.73
ACTIVSg25k	104.14	138.73	96.87	37.01	29.57	226.34	46.66
ACTIVSg70k	423.96	486.22	528.73	140.52	152.48	1269.91	182.97
case21k	230.33	{}	{}	825.15	63.54	{}	49.06
case42k	1966.92	449.81	10250.8	8745.17	198.52	{}	158.84
case99k	{}	{}	{}	{}	953.69	{}	620.7
case193k	{}	{}	{}	{}	3592.19	{}	1111.63
}\OPFbenchmarksTimePolCurr

\pgfplotstableread{
	case	MIPSsc	MIPSscpardiso	IPOPTpardiso	IPOPTma57	BELTISTOS	FMINCON	KNITRO12
case1951rte	14.53	29.97	9.34	4.04	3.96	86.23	5.19
case2383wp	6.87	8.63	11.18	5.01	5.61	19.24	4.29
case2736sp	7.7	9.62	6.14	4.21	4.69	12.47	3.56
case2737sop	5.44	6.32	4.87	3.81	5.47	9.77	3.65
case2746wop	5.23	6.87	16.36	4.06	4.5	9.64	3.77
case2746wp	6.25	7.81	4.14	3.7	4.25	11.52	3.71
case2868rte	96.18	{}	9.7	5.14	4.37	593.45	6.56
case2869pegase	{}	{}	22.18	4.97	4.32	13.61	5.36
case3012wp	6.67	8.38	35.49	4.91	5.82	18.98	4.42
case3120sp	6.31	8.16	{}	5.6	6.41	23.29	4.68
case3375wp	7.46	9.85	38.24	5.7	6.37	22.09	5.05
case6468rte	{}	{}	18.31	7.96	9.3	52.77	10.12
case6470rte	{}	{}	159.11	19.73	14.89	1599.78	14.84
case6495rte	58.21	94.68	45.24	13.18	14.09	115.46	13.97
case6515rte	{}	68.41	171.48	19.43	29.28	184.59	13.99
case9241pegase	{}	{}	{}	16.8	19.73	139.45	33.14
case\_ACTIVSg2000	7.42	9.28	8.28	4.67	4.08	12.17	3.95
case\_ACTIVSg10k	48.05	62.79	105	15.3	15.93	105.04	13.98
case13659pegase	{}	{}	277.08	137.24	68.2	{}	{}
ACTIVSg25k	150.27	192.83	{}	117.16	59.44	291.82	55.06
ACTIVSg70k	728.76	877.68	8500.76	265.47	201.65	2657.26	235.23
case21k	250.33	{}	2839.95	{}	76.82	419.82	66.25
case42k	2090.61	{}	14050.54	{}	216.45	1693.47	202.64
case99k	{}	{}	{}	{}	1150.35	11369.02	729.39
case193k	{}	{}	{}	{}	3049.82	14333.75	1778.61
}\OPFbenchmarksTimeCarCurr

\pgfplotstableread{
case	MIPSsc	MIPSscpardiso	IPOPTpardiso	IPOPTma57	BELTISTOS	FMINCON	KNITRO12
case1951rte	28	28	53	34	25	48	21
case2383wp	35	35	43	39	25	62	29
case2736sp	26	26	22	23	15	46	17
case2737sop	26	26	27	28	15	31	18
case2746wop	29	29	28	29	18	39	19
case2746wp	28	28	26	28	16	45	19
case2868rte	31	31	50	39	28	62	27
case2869pegase	32	32	72	29	19	34	22
case3012wp	30	30	40	48	36	{}	23
case3120sp	33	33	37	38	29	56	25
case3375wp	30	30	86	43	31	{}	23
case6468rte	44	44	146	38	29	{}	27
case6470rte	44	44	393	68	34	67	29
case6495rte	76	76	73	60	45	60	37
case6515rte	61	61	60	57	65	70	31
case9241pegase	46	46	65	45	26	46	64
case\_ACTIVSg2000	29	29	29	29	17	41	17
case\_ACTIVSg10k	81	81	345	35	22	52	24
case13659pegase	54	54	244	246	75	77	132
ACTIVSg25k	76	76	377	50	36	65	43
ACTIVSg70k	{}	{}	104	65	41	83	53
case21k	53	53	114	{}	53	415	38
case42k	64	64	180	{}	59	{}	44
case99k	{}	78	94	{}	70	348	50
case193k	{}	92	{}	{}	96	{}	57
}\OPFbenchmarksIters

\pgfplotstableread{
case	MIPSsc	MIPSscpardiso	IPOPTpardiso	IPOPTma57	BELTISTOS	FMINCON	KNITRO12
case1951rte	737.9	786.39	716.97	714.79	736.3	772.11	740.12
case2383wp	746.83	807.86	735.6	725.59	774.33	780.51	773.12
case2736sp	750.65	811.62	748.94	744.44	763.43	802.21	786.54
case2737sop	765.03	818.9	746.71	747.18	757.31	820.25	783.04
case2746wop	746.45	817.03	738.57	736.18	782.96	792.38	787.43
case2746wp	758.39	805.16	740.73	742	773.53	804.03	781.61
case2868rte	765.57	825.27	747.02	744.03	789.24	799.06	765.39
case2869pegase	780.68	836.65	742.2	738.77	773.98	809.27	788.59
case3012wp	748.04	813.38	742.23	739.73	803.89	{}	794.94
case3120sp	760.32	841.39	751.66	747.22	780.27	816	776.61
case3375wp	773.06	861.74	751.33	748.24	810.01	{}	794.8
case6468rte	868.32	1105.98	808.76	797.8	853.83	{}	858.33
case6470rte	858.26	1103.44	840.05	791.12	918.89	914.06	858.48
case6495rte	841.17	1266.79	784	803.68	918.91	910.62	852.35
case6515rte	850.7	1163.16	791.75	789.01	1035.17	917.28	868.42
case9241pegase	965.5	1297.51	894.79	868.82	1010.56	1069.83	976.2
case\_ACTIVSg2000	739.16	799.79	726.63	743.18	748.17	835.91	778.87
case\_ACTIVSg10k	969.37	1727.16	969.05	907.3	1037.09	1170.01	1021.93
case13659pegase	1030.89	1687.03	958.08	918.37	1683.52	1179.02	1060.05
ACTIVSg25k	1348.77	2976.82	1350.35	1275.96	1727.64	1757.24	1569.16
ACTIVSg70k	{}	{}	2372.09	2255.15	3728.51	3163.79	2889.87
case21k	1810.5	2252.1	1345.51	{}	1922.93	1908.41	1588.84
case42k	4529.98	4183.77	1979.29	{}	3337.98	{}	2533.07
case99k	{}	10632.97	3422.38	{}	7812.75	8418.02	5381.1
case193k	{}	6542.07	{}	{}	9622.57	{}	9097.03
}\OPFbenchmarksMemory

\begin{table}[ht!]
    \footnotesize
	\centering
	\caption{Time to solution (s), Polar-Power formulation. \label{tab:OPFbenchmarks-t}}
	\pgfplotstabletypeset[
	every head row/.style={ before row={\toprule},after row=\midrule},
	every last row/.style={ after row=\bottomrule},
	precision=2, fixed zerofill, column type={r},
	columns/Statistic/.style={string type,column type=l},
	every row/.style={
		column type=r,
		dec sep align,
		fixed,
		fixed zerofill,
	},
	columns={case, MIPSsc, MIPSscpardiso, IPOPTpardiso, IPOPTma57, BELTISTOS, FMINCON, KNITRO12},
	every odd row/.style={before row={\rowcolor{tablecolor1}}},
    every even row/.style={before row={\rowcolor{tablecolor2}}},
	columns/case/.style={string type, column type={l}, column name={Benchmark}},
	columns/MIPSsc/.style={column name={\MIPS{}-'\textbackslash'}},
	columns/MIPSscpardiso/.style={column name={\MIPS{}-\PARDISO{}}},
	columns/IPOPTpardiso/.style={column name={\IPOPT{}-\PARDISO{}}},
	columns/IPOPTma57/.style={column name={\IPOPT{}-MA57}},
	columns/BELTISTOS/.style={column name={\BELTISTOS{}}},
	columns/FMINCON/.style={column name={\FMINCON{}}},
	columns/KNITRO12/.style={column name={\KNITRO{}}},
	empty cells with={---} 
	]\OPFbenchmarksTimePolPow
\end{table}

\begin{table}[ht!]
	\footnotesize
	\centering
	\caption{Time to solution (s), Cartesian-Power formulation. \label{tab:OPFbenchmarks-t01}}
	\pgfplotstabletypeset[
	every head row/.style={ before row={\toprule},after row=\midrule},
	every last row/.style={ after row=\bottomrule},
	precision=2, fixed zerofill, column type={r},
	columns/Statistic/.style={string type,column type=l},
	every row/.style={
		column type=r,
		dec sep align,
		fixed,
		fixed zerofill,
	},
	columns={case, MIPSsc, MIPSscpardiso, IPOPTpardiso, IPOPTma57, BELTISTOS, FMINCON, KNITRO12},
	every odd row/.style={before row={\rowcolor{tablecolor1}}},
	every even row/.style={before row={\rowcolor{tablecolor2}}},
	columns/case/.style={string type, column type={l}, column name={Benchmark}},
	columns/MIPSsc/.style={column name={\MIPS{}-'\textbackslash'}},
	columns/MIPSscpardiso/.style={column name={\MIPS{}-\PARDISO{}}},
	columns/IPOPTpardiso/.style={column name={\IPOPT{}-\PARDISO{}}},
	columns/IPOPTma57/.style={column name={\IPOPT{}-MA57}},
	columns/BELTISTOS/.style={column name={\BELTISTOS{}}},
	columns/FMINCON/.style={column name={\FMINCON{}}},
	columns/KNITRO12/.style={column name={\KNITRO{}}},
	empty cells with={---} 
	]\OPFbenchmarksTimeCarPow
\end{table}

\begin{table}[ht!]
	\footnotesize
	\centering
	\caption{Time to solution (s), Polar-Current formulation. \label{tab:OPFbenchmarks-t10}}
	\pgfplotstabletypeset[
	every head row/.style={ before row={\toprule},after row=\midrule},
	every last row/.style={ after row=\bottomrule},
	precision=2, fixed zerofill, column type={r},
	columns/Statistic/.style={string type,column type=l},
	every row/.style={
		column type=r,
		dec sep align,
		fixed,
		fixed zerofill,
	},
	columns={case, MIPSsc, MIPSscpardiso, IPOPTpardiso, IPOPTma57, BELTISTOS, FMINCON, KNITRO12},
	every odd row/.style={before row={\rowcolor{tablecolor1}}},
	every even row/.style={before row={\rowcolor{tablecolor2}}},
	columns/case/.style={string type, column type={l}, column name={Benchmark}},
	columns/MIPSsc/.style={column name={\MIPS{}-'\textbackslash'}},
	columns/MIPSscpardiso/.style={column name={\MIPS{}-\PARDISO{}}},
	columns/IPOPTpardiso/.style={column name={\IPOPT{}-\PARDISO{}}},
	columns/IPOPTma57/.style={column name={\IPOPT{}-MA57}},
	columns/BELTISTOS/.style={column name={\BELTISTOS{}}},
	columns/FMINCON/.style={column name={\FMINCON{}}},
	columns/KNITRO12/.style={column name={\KNITRO{}}},
	empty cells with={---} 
	]\OPFbenchmarksTimePolCurr
\end{table}

\begin{table}[ht!]
	\footnotesize
	\centering
	\caption{Time to solution (s), Cartesian-Current formulation. \label{tab:OPFbenchmarks-t11}}
	\pgfplotstabletypeset[
	every head row/.style={ before row={\toprule},after row=\midrule},
	every last row/.style={ after row=\bottomrule},
	precision=2, fixed zerofill, column type={r},
	columns/Statistic/.style={string type,column type=l},
	every row/.style={
		column type=r,
		dec sep align,
		fixed,
		fixed zerofill,
	},
	columns={case, MIPSsc, MIPSscpardiso, IPOPTpardiso, IPOPTma57, BELTISTOS, FMINCON, KNITRO12},
	every odd row/.style={before row={\rowcolor{tablecolor1}}},
	every even row/.style={before row={\rowcolor{tablecolor2}}},
	columns/case/.style={string type, column type={l}, column name={Benchmark}},
	columns/MIPSsc/.style={column name={\MIPS{}-'\textbackslash'}},
	columns/MIPSscpardiso/.style={column name={\MIPS{}-\PARDISO{}}},
	columns/IPOPTpardiso/.style={column name={\IPOPT{}-\PARDISO{}}},
	columns/IPOPTma57/.style={column name={\IPOPT{}-MA57}},
	columns/BELTISTOS/.style={column name={\BELTISTOS{}}},
	columns/FMINCON/.style={column name={\FMINCON{}}},
	columns/KNITRO12/.style={column name={\KNITRO{}}},
	empty cells with={---} 
	]\OPFbenchmarksTimeCarCurr
\end{table}

\begin{table}[ht!]
    \footnotesize
	\centering
	\caption{Number of iterations. Polar-Power formulation \label{tab:OPFbenchmarks-i}}
	\pgfplotstabletypeset[
	every head row/.style={ before row={\toprule},after row=\midrule},
	every last row/.style={ after row=\bottomrule},
	precision=0, fixed zerofill, column type={r},
	columns/Statistic/.style={string type,column type=l},
	every row/.style={
		column type=r,
		dec sep align,
		fixed,
		fixed zerofill,
	},
	columns={case, MIPSsc, MIPSscpardiso, IPOPTpardiso, IPOPTma57, BELTISTOS, FMINCON, KNITRO12},
	every odd row/.style={before row={\rowcolor{tablecolor1}}},
    every even row/.style={before row={\rowcolor{tablecolor2}}},
	columns/case/.style={string type, column type={l}, column name={Benchmark}},
	columns/MIPSsc/.style={column name={\MIPS{}-'\textbackslash'}},
	columns/MIPSscpardiso/.style={column name={\MIPS{}-\PARDISO{}}},
	columns/IPOPTpardiso/.style={column name={\IPOPT{}-\PARDISO{}}},
	columns/IPOPTma57/.style={column name={\IPOPT{}-MA57}},
	columns/BELTISTOS/.style={column name={\BELTISTOS{}}},
	columns/FMINCON/.style={column name={\FMINCON{}}},
	columns/KNITRO12/.style={column name={\KNITRO{}}},
	empty cells with={---} 
	]\OPFbenchmarksIters
\end{table}

\begin{table}[ht!]
    \footnotesize
	\centering
	\caption{Maximum memory requirements (MB). Polar-Power formulations. \label{tab:OPFbenchmarks-m}}
	\pgfplotstabletypeset[
	every head row/.style={ before row={\toprule},after row=\midrule},
	every last row/.style={ after row=\bottomrule},
	precision=2, fixed zerofill, column type={r},
	columns/Statistic/.style={string type,column type=l},
	every row/.style={
		column type=r,
		dec sep align,
		fixed,
		fixed zerofill,
	},
	columns={case, MIPSsc, MIPSscpardiso, IPOPTpardiso, IPOPTma57, BELTISTOS, FMINCON, KNITRO12},
	every odd row/.style={before row={\rowcolor{tablecolor1}}},
    every even row/.style={before row={\rowcolor{tablecolor2}}},
	columns/case/.style={string type, column type={l}, column name={Benchmark}},
	columns/MIPSsc/.style={column name={\MIPS{}-'\textbackslash'}},
	columns/MIPSscpardiso/.style={column name={\MIPS{}-\PARDISO{}}},
	columns/IPOPTpardiso/.style={column name={\IPOPT{}-\PARDISO{}}},
	columns/IPOPTma57/.style={column name={\IPOPT{}-MA57}},
	columns/BELTISTOS/.style={column name={\BELTISTOS{}}},
	columns/FMINCON/.style={column name={\FMINCON{}}},
	columns/KNITRO12/.style={column name={\KNITRO{}}},
	empty cells with={---} 
	]\OPFbenchmarksMemory
\end{table}

\pgfplotstableread{
	Optimizer   SuccessAll  SolvedAll  SuccessLarge  SolvedLarge
	BELTISTOS                      100     25/25       100    6/6
	{KNITRO 12}                   100      25/25       100    6/6
	{KNITRO 11}                   96      24/25       83    5/6	
	MIPS-PARDISO               96      24/25       83     5/6
	IPOPT-PARDISO             96     24/25        83    5/6
	{MIPS-'\textbackslash'}  88      22/25       50     3/6
	IPOPT-MA57                  84      21/25       33     2/6
	{FMINCON 2018b}          80      20/25      66     4/6
}\OPFtable

\begin{table}[ht!]
\centering
\caption{Ratio of the solved benchmark cases (Polar-Power formulation). \label{tab:OPFsummary}}
\pgfplotstabletypeset[
    columns={Optimizer,SuccessAll,SuccessLarge},
	every odd row/.style={before row={\rowcolor{tablecolor1}}},
    every even row/.style={before row={\rowcolor{tablecolor2}}},
    every head row/.style={ before row=\toprule,after row=\midrule},
    every last row/.style={ after row=\bottomrule},
    precision=0,
    columns/Optimizer/.style={string type,column type=l},
    columns/SuccessAll/.style={
        column type=c,
        dec sep align,
        postproc cell content/.append code={
            \ifnum1=\pgfplotstablepartno
                \pgfkeysalso{@cell content/.add={}{\%}}%
            \fi
        },
        fixed,
        fixed zerofill,
        column name={All benchmarks}
    },
    columns/SuccessLarge/.style={
        column type=c,
        dec sep align,
        postproc cell content/.append code={
            \ifnum1=\pgfplotstablepartno
                \pgfkeysalso{@cell content/.add={}{\%}}%
            \fi
        },
        fixed,
        fixed zerofill,
        column name={Large-scale benchmarks}
    }
]\OPFtable
\end{table}

\begin{figure*}[ht!]
	\centering
	\begin{subfigure}[b]{0.5\columnwidth}
		\centering
		\input{profilesNEW/OPFbenchmarksPowerPolar-all-t.tex}
		\vspace{-0.2cm}
		\caption{Polar-Power ($\alpha_{max}=33$)\label{fig:profilePowerPolarT} }
	\end{subfigure}
	\begin{subfigure}[b]{0.49\columnwidth}
		\centering
		\begin{tikzpicture}


\pgfplotstableread{
ALPHA         SOLVER1         SOLVER2         SOLVER3         SOLVER4         SOLVER5         SOLVER6         SOLVER7
1.0         0.00000         0.00000         0.04000         0.08000         0.16000         0.00000         0.72000
1.5         0.24000         0.04000         0.48000         0.68000         0.84000         0.00000         0.96000
2.0         0.48000         0.52000         0.52000         0.72000         1.00000         0.00000         0.96000
2.5         0.52000         0.64000         0.72000         0.76000         1.00000         0.04000         1.00000
3.0         0.60000         0.68000         0.76000         0.76000         1.00000         0.16000         1.00000
3.5         0.60000         0.72000         0.76000         0.76000         1.00000         0.32000         1.00000
4.0         0.64000         0.76000         0.80000         0.84000         1.00000         0.40000         1.00000
4.5         0.64000         0.76000         0.80000         0.84000         1.00000         0.44000         1.00000
5.0         0.64000         0.76000         0.80000         0.84000         1.00000         0.52000         1.00000
5.5         0.68000         0.76000         0.80000         0.84000         1.00000         0.60000         1.00000
6.0         0.68000         0.76000         0.80000         0.84000         1.00000         0.60000         1.00000
6.5         0.68000         0.76000         0.84000         0.84000         1.00000         0.64000         1.00000
7.0         0.68000         0.76000         0.88000         0.84000         1.00000         0.68000         1.00000
7.5         0.68000         0.80000         0.88000         0.84000         1.00000         0.72000         1.00000
8.0         0.72000         0.84000         0.88000         0.84000         1.00000         0.80000         1.00000
8.5         0.72000         0.84000         0.88000         0.84000         1.00000         0.84000         1.00000
9.0         0.72000         0.84000         0.88000         0.84000         1.00000         0.88000         1.00000
9.5         0.72000         0.88000         0.92000         0.84000         1.00000         0.88000         1.00000
10.0         0.72000         0.88000         0.92000         0.84000         1.00000         0.92000         1.00000
10.5         0.72000         0.88000         0.92000         0.84000         1.00000         0.92000         1.00000
11.0         0.72000         0.88000         0.92000         0.84000         1.00000         0.92000         1.00000
11.5         0.72000         0.88000         0.92000         0.84000         1.00000         0.92000         1.00000
12.0         0.76000         0.88000         0.92000         0.84000         1.00000         0.92000         1.00000
12.5         0.76000         0.88000         0.92000         0.84000         1.00000         0.92000         1.00000
13.0         0.76000         0.92000         0.92000         0.84000         1.00000         0.96000         1.00000
13.5         0.76000         0.92000         0.92000         0.84000         1.00000         0.96000         1.00000
14.0         0.76000         0.92000         0.92000         0.84000         1.00000         0.96000         1.00000
14.5         0.76000         0.92000         0.92000         0.84000         1.00000         0.96000         1.00000
15.0         0.76000         0.92000         0.92000         0.84000         1.00000         0.96000         1.00000
15.5         0.76000         0.92000         0.96000         0.84000         1.00000         0.96000         1.00000
16.0         0.76000         0.92000         0.96000         0.84000         1.00000         0.96000         1.00000
16.5         0.76000         0.92000         0.96000         0.84000         1.00000         0.96000         1.00000
17.0         0.76000         0.92000         0.96000         0.84000         1.00000         0.96000         1.00000
17.5         0.76000         0.92000         0.96000         0.84000         1.00000         0.96000         1.00000
18.0         0.76000         0.92000         0.96000         0.84000         1.00000         0.96000         1.00000
18.5         0.76000         0.92000         0.96000         0.84000         1.00000         0.96000         1.00000
19.0         0.76000         0.92000         0.96000         0.84000         1.00000         0.96000         1.00000
19.5         0.76000         0.92000         0.96000         0.84000         1.00000         0.96000         1.00000
20.0         0.76000         0.92000         0.96000         0.84000         1.00000         0.96000         1.00000
}\datatable

\begin{axis}[name=symb,
width=\columnwidth,
height=4cm,
ymin=-0.02,
ymax=1.02,
enlarge x limits=0.02,
axis lines*=left, 
ymajorgrids, yminorgrids,
xmajorgrids, xminorgrids,
xtick = {1,5,10,15,20,25,30,35,40},
xticklabel style={rotate=0, xshift=-0.0cm, anchor=north, font=\scriptsize},
ytick={0,  0.2,  0.4, 0.6, 0.8, 1},
xmax=15,
yticklabel style={font=\scriptsize},
ylabel style={font=\scriptsize},
xlabel style={font=\scriptsize},
ylabel={$p_m({\alpha})$},
xlabel={$\alpha$},
legend style={at={(1.0,0.4)},legend cell align=left,align=right,draw=white!85!black,font=\scriptsize,legend columns=1}
]

\addplot[forget plot, color=mycolor4, loosely dashed, ultra thick] table[x=ALPHA, y=SOLVER1] {\datatable};

\addplot[forget plot, color=mycolor2, densely dashed, ultra thick] table[x=ALPHA, y=SOLVER2] {\datatable}; 

\addplot[color=mycolor4, densely dotted, ultra thick] table[x=ALPHA, y=SOLVER3] {\datatable}; 
\addlegendentry{IPOPT-PARDISO}

\addplot[color=mycolor2, loosely dotted, ultra thick] table[x=ALPHA, y=SOLVER4] {\datatable};
\addlegendentry{IPOPT-MA57}

\addplot[color=mycolor3,dashdotted, ultra thick] table[x=ALPHA, y=SOLVER5] {\datatable}; 

\addplot[color=mycolor1,  thick] table[x=ALPHA, y=SOLVER6] {\datatable}; 

\addplot[color=mycolor7, dotted, ultra thick] table[x=ALPHA, y=SOLVER7] {\datatable}; 

\end{axis}

\end{tikzpicture}
		\vspace{-0.2cm}
		\caption{Rect-Power ($\alpha_{max}=16$)\label{fig:profilePowerCartT} }
	\end{subfigure}
	\begin{subfigure}[b]{0.5\columnwidth}
		\centering
		\input{profilesNEW/OPFbenchmarksCurrentPolar-all-t.tex}
		\vspace{-0.2cm}
		\caption{Polar-Current  ($\alpha_{max}=79$)\label{fig:profileCurrentPolarT} }
	\end{subfigure}
	\begin{subfigure}[b]{0.49\columnwidth}
		\centering
		\input{profilesNEW/OPFbenchmarksCurrentCart-all-t.tex}
		\vspace{-0.2cm}
		\caption{Rect-Current ($\alpha_{max}=137$)\label{fig:profileCurrentCartT} }
	\end{subfigure}
	\caption{Overall time profile for OPF formulations considering all benchmarks.}
\end{figure*}

\begin{figure}[ht!]
	\centering
	\input{profilesNEW/OPFbenchmarksPowerPolar-all-m.tex}
	\vspace{-0.4cm}
	\caption{Memory efficiency profile for Polar-Power formulation ($\alpha_{max}=3.2$). \label{fig:profilePowerPolarM} }
\end{figure}

%
%



\clearpage
\subsection{Summary of the optimal points}
We conclude our study by reporting the optimal solutions found by all the optimization frameworks and for all initial guesses. We consider the different solutions with the same objective function value to be equivalent. We report that all optimizers found the same solution up to a relative difference $10^{-5}$.

The OPF problems are non-convex and thus different local minimums  can be reached from different starting points. We report that the relative difference $(\max(f)-\min(f))/\max(f)$ between the solutions for any given optimizer using different initial guesses is up to $10^{-5}$. The optimizers thus converged to the same solution, no matter which starting point was used (but it is the case that for poor initial guess the optimizer might not converge at all as can be seen in Tables~\ref{tab:OPFobjectiveStart4}, \ref{tab:OPFobjectiveStart5}, \ref{tab:OPFobjectiveStart6}, \ref{tab:OPFobjectiveStart7}, \ref{tab:OPFobjectiveStart8}, \ref{tab:OPFobjectiveStart9} and \ref{tab:OPFobjectiveStart10}).
\pgfplotstableread{
case start1 start2 start3 reldiff
case1951rte             81737.6751  81737.6751  81737.6751  0
case2383wp              1868170.428 1868170.428 1868170.428 0
case2736sp              1308014.964 1308014.964 1308014.964 0
case2737sop             777727.6693 777727.6693 777727.6693 0
case2746wop             1208258.467 1208258.467 1208258.467 0
case2746wp              1631707.88  1631707.88  1631707.88  0
case2868rte             79794.6792  79794.6782  79794.6782  1.25322E-08
case2869pegase          133999.2878 133999.2878 133999.2878 0
case3012wp              2591706.499 2591706.499 2591706.499 0
case3120sp              2142703.719 2142703.719 2142703.719 0
case3375wp              7412072.168 7412072.168 7412072.168 0
case6468rte             86829.0187  86829.0187  86829.0187  0
case6470rte             98345.4922  98345.4922  98345.4922  0
case6495rte             106283.3729 106283.3717 106283.3717 1.12906E-08
case6515rte             109804.2302 109804.2302 109804.2302 0
case9241pegase          315913.21   315911.5551 315911.5558 5.23846E-06
case\_ACTIVSg2000       1228892.065 1228892.065 1228892.065 0
case\_ACTIVSg10k        2485898.704 2485898.704 2485898.704 0
case13659pegase         386106.5599 386106.5599 386106.5381 5.64611E-08
case\_ACTIVSg25k        6017830.462 6017830.462 6017830.462 1.66173E-11
case\_ACTIVSg70k        16439499.16 16439499.16 16439499.16 9.73266E-11
case21k                 2592559.438 2592559.438 2592559.438 0
case42k                 2592935.865 2592935.865 2592935.865 0
case99k                 2594671.48  2594671.48  2594671.48  0
case193k                2596369.957 2596369.957 2596369.957 0
}\OPFtableStartBeltistos

\pgfplotstableread{
case start1 start2 start3 reldiff
case1951rte             81737.6751  81737.6751  81737.6751  0
case2383wp              1868170.428 1868170.428 1868170.428 0
case2736sp              1308014.964 1308014.964 1308014.965 2.29355E-10
case2737sop             777727.6693 777727.6693 777727.6693 0
case2746wop             1208258.467 1208258.467 1208258.467 0
case2746wp              1631707.881 1631707.88  1631707.88  4.28998E-10
case2868rte             79794.6792  79794.6782  79794.6782  1.25322E-08
case2869pegase          133999.2878 133999.2878 133999.2878 0
case3012wp              2591706.499 2591706.499 2591706.499 0
case3120sp              2142703.719 2142703.719 2142703.719 0
case3375wp              7412072.168 7412072.168 7412072.168 0
case6468rte             86829.0187  86829.0187  86829.0187  0
case6470rte             98345.4922  98345.4922  98345.4922  0
case6495rte             106283.3717 106283.3717 106283.3717 0
case6515rte             109804.2302 109804.2302 109804.2302 0
case9241pegase          315912.1625 315911.5558 315911.5551 1.92269E-06
case\_ACTIVSg2000       1228892.065 1228892.065 1228892.065 0
case\_ACTIVSg10k        2485898.704 2485898.704 2485898.704 0
case13659pegase         {} {} 386106.5381 0
case\_ACTIVSg25k        6017830.462 6017830.462 6017830.462 1.66173E-11
case\_ACTIVSg70k        16439499.16 {} 16439499.16 9.73266E-11
case21k                 2592559.442 2592559.437 2592559.438 1.81288E-09
case42k                 2592935.866 2592935.865 2592935.865 6.94194E-10
case99k                 2594671.48  2594671.48  2594671.48  0
case193k                {} {} {} 0
}\OPFtableStartIpoptPardiso

\pgfplotstableread{
case start1 start2 start3 reldiff
case1951rte            81737.6751   81737.6751  81737.6751  0
case2383wp             1868170.428  1868170.428 1868170.428 0
case2736sp             1308014.964  1308014.964 1308014.964 0
case2737sop            777727.6693  777727.6693 777727.6693 0
case2746wop            1208258.467  1208258.467 1208258.467 0
case2746wp             1631707.881  1631707.881 1631707.881 6.12855E-11
case2868rte            79794.6782   79794.6782  79794.6782  0
case2869pegase         133999.2878  133999.2878 133999.2878 0
case3012wp             2591706.499  2591706.499 2591706.499 0
case3120sp             2142703.719  2142703.719 2142703.719 0
case3375wp             7412072.168  7412072.168 7412072.168 0
case6468rte            86829.0187   86829.0187  86829.0187  0
case6470rte            98345.4922   98345.4922  98345.4922  0
case6495rte            106283.3717  106283.3717 106283.3717 0
case6515rte            109804.2302  109804.2302 109804.2302 0
case9241pegase         315911.5556  315911.5551 315911.5559 2.53235E-09
case\_ACTIVSg2000      1228892.065  1228892.065 1228892.065 0
case\_ACTIVSg10k       2485898.704  2485898.704 2485898.704 0
case13659pegase        386106.5599  386106.5522 386106.5161 1.1344E-07
case\_ACTIVSg25k       6017830.462  6017830.462 6017830.462 1.66173E-11
case\_ACTIVSg70k       16439499.16  16439499.16 16439499.16 9.73266E-11
case21k                2592559.438  2592559.438 {} 0
case42k                2592935.857  {} {} 0
case99k                {}  {} {} 0
case193k               {}  {} {} 0
}\OPFtableStartIpoptHSL

\pgfplotstableread{
case start1 start2 start3 reldiff
case1951rte            {}  81737.6756  81737.6756  0
case2383wp             1868170.494  1868170.493 1868170.493 5.35283E-11
case2736sp             1308014.997  1308014.997 1308014.997 2.29355E-10
case2737sop            777727.6848  777727.6847 777727.6848 1.2858E-10
case2746wop            1208258.503  1208258.503 1208258.503 1.65527E-10
case2746wp             1631707.935  1631707.935 1631707.935 6.12855E-11
case2868rte            {}  79794.6795  79794.6795  0
case2869pegase         133999.2881  133999.2881 133999.2881 0
case3012wp             2591706.566  2591706.566 2591706.566 7.71691E-11
case3120sp             2142703.763  2142703.765 2142703.766 1.02674E-09
case3375wp             7412072.199  7412072.199 7412072.199 1.34915E-11
case6468rte            {}  86829.0191  86829.0191  0
case6470rte            {}  98345.4935  98345.4935  0
case6495rte            {}  106283.3727 106283.3727 0
case6515rte            {}  109804.2417 109804.2417 0
case9241pegase         315912.4336  315911.5707 315911.5757 2.73145E-06
case\_ACTIVSg2000      1228892.076  1228892.076 1228892.076 0
case\_ACTIVSg10k       {}  2485898.751 2485898.751 4.02268E-11
case13659pegase        386117.3613  386117.4576 386113.5176 1.02041E-05
case\_ACTIVSg25k       {}  6017830.612 6017830.612 1.66173E-11
case\_ACTIVSg70k       {}  {} {} 0
case21k                2592559.524  2592559.522 2592559.522 6.94295E-10
case42k                2592935.95   2592935.948 2592935.948 5.78495E-10
case99k                {}  {} {} 0
case193k               {}  {} {} 0
}\OPFtableStartMIPSsc

\pgfplotstableread{
case start1 start2 start3 reldiff
case1951rte            {}  81737.6756  81737.6756  0
case2383wp             1868170.494  1868170.493 1868170.493 5.35283E-11
case2736sp             1308014.997  1308014.997 1308014.997 2.29355E-10
case2737sop            777727.6848  777727.6847 777727.6848 1.2858E-10
case2746wop            1208258.503  1208258.503 1208258.503 1.65527E-10
case2746wp             1631707.935  1631707.935 1631707.935 6.12855E-11
case2868rte            {}  79794.6795  79794.6795  0
case2869pegase         133999.2881  133999.2881 133999.2881 0
case3012wp             2591706.566  2591706.566 2591706.566 7.71691E-11
case3120sp             2142703.763  2142703.765 2142703.766 1.02674E-09
case3375wp             7412072.199  7412072.199 7412072.199 1.34915E-11
case6468rte            {}  86829.0191  86829.0191  0
case6470rte            {}  98345.4935  98345.4935  0
case6495rte            {}  106283.3727 106283.3727 0
case6515rte            {}  109804.2417 109804.2417 0
case9241pegase         315912.4336  315911.5707 315911.5757 2.73145E-06
case\_ACTIVSg2000      1228892.076  1228892.076 1228892.076 0
case\_ACTIVSg10k       {}  2485898.751 2485898.751 4.02268E-11
case13659pegase        386114.0263  386117.4576 386113.5176 1.02041E-05
case\_ACTIVSg25k       {}  6017830.612 6017830.612 1.66173E-11
case\_ACTIVSg70k       {}  {} {} 0
case21k                2592559.524  2592559.522 2592559.522 6.94295E-10
case42k                2592935.963  2592935.948 2592935.943 7.559E-09
case99k                2594671.538  2594671.538 2594671.536 6.55189E-10
case193k               2596370.012  2596370.011 2596370.011 2.31092E-10
}\OPFtableStartMIPSscPardiso

\pgfplotstableread{
case start1 start2 start3 reldiff
case1951rte            {}    81737.679   81737.7618  1.013E-06
case2383wp             1868170.494  1868170.587 1868170.587 4.97813E-08
case2736sp             1308014.999  1308015.007 1308015.007 6.34549E-09
case2737sop            777727.6862  777727.6862 777727.8223 1.74997E-07
case2746wop            1208258.516  1208258.569 1208258.569 4.40303E-08
case2746wp             1631707.95   1631707.95  1631707.935 9.19282E-09
case2868rte            {}    79794.683   79794.7004  2.1806E-07
case2869pegase         133999.2883  133999.2889 133999.2883 4.47764E-09
case3012wp             2591706.566  {}   {}   0
case3120sp             2142703.837  {}   2142703.78  2.63686E-08
case3375wp             7412072.326  {}   {}   0
case6468rte            {}    {}   {}   0
case6470rte            {}    98345.4921  98345.5255  3.39619E-07
case6495rte            {}    106283.5523 106283.4079 1.35863E-06
case6515rte            {}    109804.238  109804.238  0
case9241pegase         315912.3379  315911.5567 315911.5679 2.47284E-06
case\_ACTIVSg2000      1228892.166  1228892.094 1228892.076 7.36436E-08
case\_ACTIVSg10k       {}    2485898.829 2485899.14  1.24864E-07
case13659pegase        386107.7082  386106.5137 386106.5179 3.0937E-06
case\_ACTIVSg25k       6017831.823  6017831.584 6017831.823 3.97485E-08
case\_ACTIVSg70k       16439502.42  16439499.94 16439500.35 1.51318E-07
case21k                2592559.502  {}   2592559.532 1.17259E-08
case42k                2592935.925  2592936.002 {}   2.97732E-08
case99k                {}    2594671.563 2594686.05  5.58349E-06
case193k               2596370.287  {}   {}   0
}\OPFtableStartFmincon

\pgfplotstableread{
case start1 start2 start3 reldiff
case1951rte            81737.6756   81737.6756  81737.6756  0
case2383wp             1868170.493  1868170.493 1868170.493 0
case2736sp             1308014.996  1308014.996 1308014.996 3.82258E-10
case2737sop            777727.6848  777727.6848 777727.6848 0
case2746wop            1208258.503  1208258.503 1208258.503 8.27638E-11
case2746wp             1631707.935  1631707.934 1631707.935 4.28998E-10
case2868rte            79794.6801   79794.6786  79794.6787  1.87982E-08
case2869pegase         133999.2881  133999.2881 133999.2881 0
case3012wp             2591706.566  2591706.566 2591706.566 3.85845E-11
case3120sp             2142703.765  2142703.765 2142703.765 4.66701E-11
case3375wp             7412072.2    7412072.198 7412072.2   2.02373E-10
case6468rte            86829.0191   86829.0191  86829.0191  0
case6470rte            98345.493    98345.493   98345.493   0
case6495rte            106283.3727  106283.3727 106283.3727 0
case6515rte            109804.2312  109804.2312 109804.2312 0
case9241pegase         315912.4058  315911.5561 315911.5561 2.68967E-06
case\_ACTIVSg2000      1228892.076  1228892.075 1228892.075 2.44122E-10
case\_ACTIVSg10k       2485898.751  2485898.753 2485898.753 6.83857E-10
case13659pegase        386107.515   386106.5528 386106.5136 2.59358E-06
case\_ACTIVSg25k       6017830.612  6017830.613 6017830.612 4.98518E-11
case\_ACTIVSg70k       16439499.83  16439499.83 16439499.83 7.29949E-11
case21k                2592559.511  2592559.509 2592559.506 1.77431E-09
case42k                2592935.979  {} 2592935.936 1.66607E-08
case99k                2594671.556  2594671.601 {} 1.72662E-08
case193k               2596451.881  2596451.83  2596454.235 9.26186E-07
}\OPFtableStartKNITRO

\begin{table}[H]
    \footnotesize
	\centering
	\caption{Final objective function value (\$/h) for \IPOPT{}-\PARDISO{}. \label{tab:OPFobjectiveStart4}}
	\pgfplotstabletypeset[
	every head row/.style={ before row={\toprule},after row=\midrule},
	every last row/.style={ after row=\bottomrule},
	precision=2, fixed, fixed zerofill, column type={r},
	columns/Statistic/.style={string type,column type=l},
	every row/.style={
		column type=r,
		dec sep align,
		fixed,
		fixed zerofill,
	},
	columns={case, start1, start2, start3, reldiff},
	every odd row/.style={before row={\rowcolor{tablecolor1}}},
    every even row/.style={before row={\rowcolor{tablecolor2}}},
	columns/case/.style={string type, column type={l}, column name={Benchmark}},
	columns/start1/.style={column name={Flat start}},
	columns/start2/.style={column name={\matpower{} case data}},
	columns/start3/.style={column name={Power flow solution}},
	columns/reldiff/.style={sci, column name={Rel. error}},
	empty cells with={---} 
	]\OPFtableStartIpoptPardiso
\end{table}
\begin{table}[H]
    \footnotesize
	\centering
	\caption{Final objective function value (\$/h) for \BELTISTOS{}. \label{tab:OPFobjectiveStart5}}
	\pgfplotstabletypeset[
	every head row/.style={ before row={\toprule},after row=\midrule},
	every last row/.style={ after row=\bottomrule},
	precision=2, fixed, fixed zerofill, column type={r},
	columns/Statistic/.style={string type,column type=l},
	every row/.style={
		column type=r,
		dec sep align,
		fixed,
		fixed zerofill,
	},
	columns={case, start1, start2, start3, reldiff},
	every odd row/.style={before row={\rowcolor{tablecolor1}}},
    every even row/.style={before row={\rowcolor{tablecolor2}}},
	columns/case/.style={string type, column type={l}, column name={Benchmark}},
	columns/start1/.style={column name={Flat start }},
	columns/start2/.style={column name={\matpower{} case data}},
	columns/start3/.style={column name={Power flow solution}},
	columns/reldiff/.style={sci, column name={Rel. error}},
	empty cells with={---} 
	]\OPFtableStartBeltistos
\end{table}
\begin{table}[H]
    \footnotesize
	\centering
	\caption{Final objective function value (\$/h) for \IPOPT{}-MA57. \label{tab:OPFobjectiveStart6}}
	\pgfplotstabletypeset[
	every head row/.style={ before row={\toprule},after row=\midrule},
	every last row/.style={ after row=\bottomrule},
	precision=2, fixed, fixed zerofill, column type={r},
	columns/Statistic/.style={string type,column type=l},
	every row/.style={
		column type=r,
		dec sep align,
		fixed,
		fixed zerofill,
	},
	columns={case, start1, start2, start3, reldiff},
	every odd row/.style={before row={\rowcolor{tablecolor1}}},
    every even row/.style={before row={\rowcolor{tablecolor2}}},
	columns/case/.style={string type, column type={l}, column name={Benchmark}},
	columns/start1/.style={column name={Flat start}},
	columns/start2/.style={column name={\matpower{} case data}},
	columns/start3/.style={column name={Power flow solution}},
	columns/reldiff/.style={sci, column name={Rel. error}},
	empty cells with={---} 
	]\OPFtableStartIpoptHSL
\end{table}
\begin{table}[H]
    \footnotesize
	\centering
	\caption{Final objective function value (\$/h) for \MIPS{}-MATLAB'\textbackslash'. \label{tab:OPFobjectiveStart7}}
	\pgfplotstabletypeset[
	every head row/.style={ before row={\toprule},after row=\midrule},
	every last row/.style={ after row=\bottomrule},
	precision=2, fixed, fixed zerofill, column type={r},
	columns/Statistic/.style={string type,column type=l},
	every row/.style={
		column type=r,
		dec sep align,
		fixed,
		fixed zerofill,
	},
	columns={case, start1, start2, start3, reldiff},
	every odd row/.style={before row={\rowcolor{tablecolor1}}},
    every even row/.style={before row={\rowcolor{tablecolor2}}},
	columns/case/.style={string type, column type={l}, column name={Benchmark}},
	columns/start1/.style={column name={Flat start}},
	columns/start2/.style={column name={\matpower{} case data}},
	columns/start3/.style={column name={Power flow solution}},
	columns/reldiff/.style={sci, column name={Rel. error}},
	empty cells with={---} 
	]\OPFtableStartMIPSsc
\end{table}
\begin{table}[H]
    \footnotesize
	\centering
	\caption{Final objective function value (\$/h) for \MIPS{}-\PARDISO{}. \label{tab:OPFobjectiveStart8}}
	\pgfplotstabletypeset[
	every head row/.style={ before row={\toprule},after row=\midrule},
	every last row/.style={ after row=\bottomrule},
	precision=2, fixed, fixed zerofill, column type={r},
	columns/Statistic/.style={string type,column type=l},
	every row/.style={
		column type=r,
		dec sep align,
		fixed,
		fixed zerofill,
	},
	columns={case, start1, start2, start3, reldiff},
	every odd row/.style={before row={\rowcolor{tablecolor1}}},
    every even row/.style={before row={\rowcolor{tablecolor2}}},
	columns/case/.style={string type, column type={l}, column name={Benchmark}},
	columns/start1/.style={column name={Flat start}},
	columns/start2/.style={column name={\matpower{} case data}},
	columns/start3/.style={column name={Power flow solution}},
	columns/reldiff/.style={sci, column name={Rel. error}},
	empty cells with={---} 
	]\OPFtableStartMIPSscPardiso
\end{table}
\begin{table}[H]
    \footnotesize
	\centering
	\caption{Final objective function value (\$/h) for \FMINCON{} 2018b. \label{tab:OPFobjectiveStart9}}
	\pgfplotstabletypeset[
	every head row/.style={ before row={\toprule},after row=\midrule},
	every last row/.style={ after row=\bottomrule},
	precision=2, fixed, fixed zerofill, column type={r},
	columns/Statistic/.style={string type,column type=l},
	every row/.style={
		column type=r,
		dec sep align,
		fixed,
		fixed zerofill,
	},
	columns={case, start1, start2, start3, reldiff},
	every odd row/.style={before row={\rowcolor{tablecolor1}}},
    every even row/.style={before row={\rowcolor{tablecolor2}}},
	columns/case/.style={string type, column type={l}, column name={Benchmark}},
	columns/start1/.style={column name={Flat start}},
	columns/start2/.style={column name={\matpower{} case data}},
	columns/start3/.style={column name={Power flow solution}},
	columns/reldiff/.style={sci, column name={Rel. error}},
	empty cells with={---} 
	]\OPFtableStartFmincon
\end{table}
\begin{table}[H]
    \footnotesize
	\centering
	\caption{Final objective function value (\$/h) for \KNITRO{} 11. \label{tab:OPFobjectiveStart10}}
	\pgfplotstabletypeset[
	every head row/.style={ before row={\toprule},after row=\midrule},
	every last row/.style={ after row=\bottomrule},
	precision=2, fixed, fixed zerofill, column type={r},
	columns/Statistic/.style={string type,column type=l},
	every row/.style={
		column type=r,
		dec sep align,
		fixed,
		fixed zerofill,
	},
	columns={case, start1, start2, start3, reldiff},
	every odd row/.style={before row={\rowcolor{tablecolor1}}},
    every even row/.style={before row={\rowcolor{tablecolor2}}},
	columns/case/.style={string type, column type={l}, column name={Benchmark}},
	columns/start1/.style={column name={Flat start}},
	columns/start2/.style={column name={\matpower{} case data}},
	columns/start3/.style={column name={Power flow solution}},
	columns/reldiff/.style={sci, column name={Rel. error}},
	empty cells with={---} 
	]\OPFtableStartKNITRO
\end{table}


\section{Numerical results - MPOPF}
\label{sec:ResultsMPOPF}


The benchmarks in this section focus on performance related to the KKT solution (factorization and solution phases) since these represent the bottleneck of the IP method for large-scale MPOPF problems. We do not consider common factors for IP methods, including computation such as assembly of the MPOPF case or symbolic factorization routines, which are performed only once per the optimization algorithm, or other routines (vector updates, evaluations of stopping criteria, etc.) performed in each IP iteration.

\begin{table}[t!]
    \centering
    \footnotesize
    \caption{Selected MPOPF benchmark statistics including the number of time periods $N$, storage devices $\Ns$ and the number linear constraints $|A(x)|$. \label{tab:benchmarksMPOPF}}
\begin{filecontents*}{./data/CasesTableMPOPF.dat}
CaseStatistics  n   ns      nb  ng nl  nlc  ndc nvar nnle   nnli   nlin
case118 600 10  118 54  186 0   0   230400 141600 223200 6000
1369pegase  600 10  1354    260 1991    1432    0   1960800   1624800   2389200   6000
2869pegase  600 10  2869    510 4582    2743    0   4078800   3442800   5498400   6000
case118 1200    10  118 54  186 0   0   460800 283200 446400 12000
1369pegase  1200    10  1354    260 1991    1432    0   3921600   3249600   4778400   12000
2869pegase  1200    10  2869    510 4582    2743    0   8157600   6885600   10996800  12000
case118 4800    10  118 54  186 0   0   1843200   1132800   1785600   48000
1369pegase  4800    10  1354    260 1991    1432    0   15686400  12998400  19113600  48000
2869pegase  4800    10  2869    510 4582    2743    0   32630400  27542400  43987200  48000
\end{filecontents*}
\pgfplotstabletypeset[
    columns={CaseStatistics,n, ns,nvar,nnle,nnli,nlin},
    every head row/.style={ before row=\toprule,after row=\midrule},
    precision=0, fixed, fixed zerofill, column type={r},
    columns/CaseStatistics/.style={column name=Benchmark,string type,column type=l},
    columns/n/.style={column name=$N$},
    columns/ns/.style={column name=$\Ns$},
    columns/nb/.style={column name=$n_b$},
    columns/nvar/.style={column name=$|x|$},
    columns/ng/.style={column name=$n_g$},
    columns/nl/.style={column name=$n_l$},
    columns/nnle/.style={column name=$|g(x)|$},
    columns/nnli/.style={column name=$|h(x)|$},
    columns/nlin/.style={column name=$|A(x)|$},
    every last row/.style={ after row=\bottomrule},
    every row no 0/.style={ before row={\rowcolor{tablecolor2}}},
    every row no 1/.style={ before row={\rowcolor{tablecolor2}}},
    every row no 2/.style={ before row={\rowcolor{tablecolor2}}},
    every row no 3/.style={ before row={\rowcolor{tablecolor1}}},
    every row no 4/.style={ before row={\rowcolor{tablecolor1}}},
    every row no 5/.style={ before row={\rowcolor{tablecolor1}}},
    every row no 6/.style={ before row={\rowcolor{tablecolor2}}},
    every row no 7/.style={ before row={\rowcolor{tablecolor2}}},
    every row no 8/.style={ before row={\rowcolor{tablecolor2}}},
    every row no 9/.style={ before row={\rowcolor{tablecolor2}}},
    every nth row={19}{before row=\midrule \multicolumn{6}{l}{Large-scale benchmarks}\\ \midrule \rowcolor{tablecolor1}},
]
{./data/CasesTableMPOPF.dat}
\end{table}
\subsection{Benchmarks setup}
For all MPOPF simulations, the length of the time period $\delta t$ is set to one hour. The  load scaling profile
was used for all benchmarks to generate a time dependent load as a multiplicator of the nominal loads given in the \matpower{} case files. The profile is based on a load pattern for the Swiss canton Ticino\footnote{Data available at \url{https://www.swissgrid.ch/en/home/operation/grid-data/generation.html}.}.
In each \matpower{} case file, the energy storage devices are located at the first $\Ns$ buses sorted according to the largest positive active load demand specified in the case file.  The storage size $\bepsilonmaxS$ is chosen to contain up to two hours of the nominal active power demand $P^d_{k}$ of the bus $k$.  If the storage device $i$ is located at bus $k$, then $\bepsilonmaxSi = 2\,P^d_{k}$ and $\bepsilonminSi = 0$.  The initial state of charge is $70$\%, which represents $\bepsilon_0 = 0.7\bepsilonmaxS$.  The storage device power ratings are limited to allow a  complete discharging and charging within three hours and two hours, respectively. Therefore, ${P_g^{sd}}^{\max} = \frac{1}{3}\bepsilonmaxS$ and ${P_g^{sc}}^{\min} = -\frac{1}{2}\bepsilonmaxS$.  All storage device discharging and charging efficiencies are chosen as $\eta_{\text{d}} = 0.97 $ and $\eta_{\text{c}} = 0.95 $. 
Further changes with respect to the original \matpower{} cases include line flow limits for every line such that the imposed limit is greater than the line flow in the OPF solution considering the nominal load. In this way, the limits are not restricting the feasible solution but the computational complexity better reflects the real-life conditions. Additionally, the minimum generator output limits were relaxed.

The benchmark cases include one small and two medium size power networks, listed in table \ref{tab:benchmarksMPOPF}. The table illustrates MPOPF problem sizes for different number of time periods and number of storage devices.

\begin{figure}[!t]
\centering
\begin{tikzpicture}
\begin{axis}[
width=0.8\columnwidth,
height=0.15\columnwidth,
xlabel/.append style={font=\scriptsize},
ylabel/.append style={font=\scriptsize},
at={(0.0in,0.0in)},
scale only axis,
separate axis lines,
every outer x axis line/.append style={black},
every x tick label/.append style={font=\scriptsize},
ymajorgrids,
xmin=1,
xmax=240,
xtick={1,50,100,150,200,240},
xtick pos=left,
ytick pos=left,
xlabel={Time $n$ (hours)},
every outer y axis line/.append style={black},
every y tick label/.append style={font=\scriptsize},
ymode=normal,
ymin=0.2,
ymax=1.1,
ytick={0,0.2,0.4,0.6,0.8,1.0,1.2,1.4,1.6,1.8,2.0},
yminorticks=true,
ylabel={Load scaling coefficient},
axis background/.style={fill=none},
legend style={at={(0.25,1.)},legend cell align=left,align=left,draw=none,fill=none, font=\scriptsize},
]

\addplot[color=mycolor5,line width=1.5pt,mark size=1.0pt,mark=none,mark options={fill=mycolor5}]
 table[x=Time,y=LoadScaled]{data/InputData.dat};

\end{axis}

\end{tikzpicture}%
\caption{Hourly demand scaling coefficient $k^\text{L}$ over a 10 day period.}
\label{fig:loadProfile}
\end{figure}

\subsection{Number of time periods and storage devices \label{sec:resultsTimePeriods}}
We investigate the performance of \IPOPT{}, \KNITRO{} and \BELTISTOSMP{} for increasing problem sizes, changing the number of time periods and storage devices. The investigation is focused on the bottleneck of the IP method, which is the solution of the KKT system in each IP iteration.
Since the KKT matrix changes numerically but not structurally at every iteration, it is therefore reasonable to separate the symbolic factorisation phase that determines a sparsity preserving pivot order from the numerical factorisation phase. The symbolic factorisation phase only needs to be done once at the beginning of the IP algorithm.

The average time of the KKT system solution, consisting of the numerical factorization and forward-backward substitutions (excluding the symbolic factorization) is shown in Figure \ref{fig:increasingN_timeStacked} for increasing number of time periods. Different power grids are shown, each containing ten storage devices. It is evident that \BELTISTOSMP{} outperforms the black-box solution approaches, providing orders of magnitude faster solution times. Figure \ref{fig:increasingN_timeStacked} also provides comparison of the factorization and forward-backward substitution phases, illustrated by the horizontal line inside the bars (note the logarithmic scale on y-axis). 
The factorization phase clearly dominates for the \IPOPT{} black-box approach, therefore the forward-backward substitution phase is no visible in the figure. The factorization and the backsubstitution phases are comparable for the \BELTISTOSMP{}, with the backsubstitution phase taking more time for the memory efficient version of \BELTISTOSMP{} - \BELTISTOSMEM{}, since some portion of the computation needs to be recomputed redundantly in order to reduce the memory requirements of the algorithm. However, the performance benefit is still significant, compared to the standard solution method. We observe that the performance gap between \BELTISTOSMP{} and all other black-box solvers increases with increasing values of $N$. For $N=4800$ such solvers failed due to exceeding the memory limit during the symbolic factorization of matrix. With the increasing benchmark size, the failure was observed also for $N=2400$ or $N=1200$ due to exceeding the available memory or the time limit. \BELTISTOSMP{} requires approximately 1\% of the time needed by the best competitor for the smallest problem, with increasing performance benefit for larger problems. 

Furthermore, we would like to investigate the performance of \BELTISTOSMP{} for increasing number of storage devices. This is particularly important for storage sizing and placement problems, where the optimal size and location of the storage devices are sought for. The number of storage devices increases the number of coupling variables and thus the size of the Schur complement, therefore also posing a bottleneck for large problems. Figure \ref{fig:increasingNs_timeStacked} shows the average solution time of the KKT system for increasing number of storage devices $N_s$ for different power grids, considering $N=600$ time periods. In case of the case1354pegase benchmark, the 10 fold increase in the storage devices resulted in 81 times longer computation for the \IPOPT{}, while only 3.7 times increase for the computation time of \BELTISTOSMP{}.

\begin{figure*}[t!]
	\centering
	\input{fig_AvgTimePerIterTime_stacked.tex}
	\vspace{-0.4cm}
	\caption{ Average numerical  factorization  and forward-backward substitutions time ($N_s=10$).} 
	\label{fig:increasingN_timeStacked}
\end{figure*}

\begin{figure*}[t!]
	\centering
	\input{fig_AvgTimePerIterTimeNs_stacked.tex}
	\vspace{-0.4cm}
	\caption{ Average numerical  factorization  and forward-backward substitutions time ($N=600$).} 
	\label{fig:increasingNs_timeStacked}
\end{figure*}

\begin{figure*}[t!]
	\centering
	\newcommand\offone{-8.0ex}
\newcommand\offtwo{-10.5ex}
\newcommand\offthree{-4.5ex}


\begin{tikzpicture}[scale=1.0]

\pgfplotstableread{
N	OPF	MPOPF	MPOPFmem KNITRO
600	4050.191406	750.9659082	85.2627832 10922
1200	10292.804688	1501.934072	169.4340723  49457
2400	35070.82031	3003.87041	337.7766602  268754
4800	3e6	6007.743096	674.4618457  3e6
}\datatableA 

\pgfplotstableread{
N	OPF	MPOPF	MPOPFmem KNITRO
600	66081.78516	7416.937832	834.3616602  32343
1200	108758.82031	14833.8789	1657.748037  84764
2400	3e6	29667.76104	3304.520801  397968
4800	3e6	59335.52532	6598.066338  3e6
}\datatableB 

\pgfplotstableread{
N	OPF	MPOPF	MPOPFmem KNITRO
600	107606.19531	16751.9172	1872.257041  53191
1200	3e6	33503.84057	3719.688223  3e6
2400	3e6	67017.06328	7414.550586  3e6
4800	3e6	134034.1318	14804.27532   3e6
}\datatableC 

\begin{groupplot}[max space between ticks=15, group style = {group size = 3 by 1, horizontal sep = 20pt}, width=\textwidth, height = 5.0cm]
    \nextgroupplot[ title = {Case 118},
    width=0.36\textwidth,
    height=0.2\textwidth,
    separate axis lines,
    ymajorgrids,
    yminorgrids,
    minor tick num=1,
    thick,
    every axis plot post/.style={/pgf/number format/fixed},
    ybar,
    bar width=3pt,
    ymin=1,
    ymax = 1e6,
    xtick=data,
    ytick={1,1e1,1e2,1e3,1e4,1e5,1e6},
    ymode=log,
    log basis y=10,
    yminorticks=true,
    enlarge x limits=0.1,
    separate axis lines,
    every outer x axis line/.append style={black},
    every x tick label/.append style={font=\scriptsize},
    every y tick label/.append style={font=\scriptsize},
    symbolic x coords={600,	1200,	2400,	4800},
    visualization depends on=rawy\as\rawy, 
    after end axis/.code={ 
            \draw [ultra thick, white, decoration={snake, amplitude=1pt}, decorate] (rel axis cs:0,1.04) -- (rel axis cs:1,1.04);
        },
    axis lines*=left,
    clip=false,
    axis background/.style={fill=none},
    legend style={at={(1.0,1.25)},legend cell align=left,align=left,draw=none,fill=none,font=\scriptsize,legend columns=4, legend to name = MPOPFgrouplegendMEM},
    xlabel/.append style={font=\scriptsize},
    ylabel/.append style={font=\scriptsize},
    title style={yshift=-1ex,},
    xlabel={$N$},
    ylabel={Memory (MB)} ]
    
    \addplot[fill=mycolor2] table [x=N, y=MPOPF]\datatableA;
    \addlegendentry{\BELTISTOSMP{}}
    
    \addplot[fill=mycolor3] table [x=N, y=MPOPFmem]\datatableA;
    \addlegendentry{\BELTISTOSMEM{}mem}
    
    \addplot[fill=mycolor1] table [x=N, y=OPF]\datatableA;
    \addlegendentry{\IPOPT{}}
    
    \addplot[fill=mycolor4] table [x=N, y=KNITRO]\datatableA;
    \addlegendentry{\KNITRO{}}

    \nextgroupplot[ title = {1354 Pegase},
    width=0.36\textwidth,
    height=0.2\textwidth,
    separate axis lines,
    ymajorgrids,
    yminorgrids,
    minor tick num=1,
    thick,
    every axis plot post/.style={/pgf/number format/fixed},
    ybar,
    bar width=3pt,
    ymin=1,
    ymax = 1e6,
    xtick=data,
    ytick={1,1e1,1e2,1e3,1e4,1e5,1e6},
    ymode=log,
    log basis y=10,
    yminorticks=true,
    enlarge x limits=0.1,
    separate axis lines,
    every outer x axis line/.append style={black},
    every x tick label/.append style={font=\scriptsize},
    every y tick label/.append style={font=\scriptsize},
    symbolic x coords={600,	1200,	2400,	4800},
    visualization depends on=rawy\as\rawy, 
    after end axis/.code={ 
            \draw [ultra thick, white, decoration={snake, amplitude=1pt}, decorate] (rel axis cs:0,1.04) -- (rel axis cs:1,1.04);
        },
    axis lines*=left,
    clip=false,
    axis background/.style={fill=none},
    legend style={at={(1.0,1.15)},legend cell align=left,align=left,draw=none,fill=none,font=\scriptsize,legend columns=3},
    xlabel/.append style={font=\scriptsize},
    ylabel/.append style={font=\scriptsize},
    title style={yshift=-1ex,},
    xlabel={$N$},    
    ]
    
    \addplot[fill=mycolor2] table [x=N, y=MPOPF]\datatableB;
    
    \addplot[fill=mycolor3] table [x=N, y=MPOPFmem]\datatableB;
    
    \addplot[fill=mycolor1] table [x=N, y=OPF]\datatableB;
    
    \addplot[fill=mycolor4] table [x=N, y=KNITRO]\datatableB;
    \nextgroupplot[ title = {2869 Pegase},
    width=0.36\textwidth,
    height=0.2\textwidth,
    separate axis lines,
    ymajorgrids,
    yminorgrids,
    minor tick num=1,
    thick,
    every axis plot post/.style={/pgf/number format/fixed},
    ybar,
    bar width=3pt,
    ymin=1,
    ymax = 1e6,
    xtick=data,
    ytick={1,1e1,1e2,1e3,1e4,1e5,1e6},
    ymode=log,
    log basis y=10,
    yminorticks=true,
    enlarge x limits=0.1,
    separate axis lines,
    every outer x axis line/.append style={black},
    every x tick label/.append style={font=\scriptsize},
    every y tick label/.append style={font=\scriptsize},
    symbolic x coords={600,	1200,	2400,	4800},
    visualization depends on=rawy\as\rawy, 
    after end axis/.code={ 
            \draw [ultra thick, white, decoration={snake, amplitude=1pt}, decorate] (rel axis cs:0,1.04) -- (rel axis cs:1,1.04);
        },
    axis lines*=left,
    clip=false,
    axis background/.style={fill=none},
    legend style={at={(1.0,1.15)},legend cell align=left,align=left,draw=none,fill=none,font=\tiny,legend columns=3},
    xlabel/.append style={font=\scriptsize},
    ylabel/.append style={font=\scriptsize},
    title style={yshift=-1ex,},
    xlabel={$N$},    
    ]
    
    \addplot[fill=mycolor2] table [x=N, y=MPOPF]\datatableC;
    
    \addplot[fill=mycolor3] table [x=N, y=MPOPFmem]\datatableC;
    
    \addplot[fill=mycolor1] table [x=N, y=OPF]\datatableC;
    
    \addplot[fill=mycolor4] table [x=N, y=KNITRO]\datatableC;
    
\end{groupplot}
\node at ($(group c2r1) + (0,-2cm)$) {\ref*{MPOPFgrouplegendMEM}}; 

\end{tikzpicture}
	\vspace{-0.4cm}
	\caption{ Memory requirements for the optimizer ($N_s=10$).}
	\label{fig:MPOPFmemoryComplexity}
\end{figure*}

\subsection{Memory complexity}
In this section we examine the memory efficiency of the solvers. Figure \ref{fig:MPOPFmemoryComplexity} illustrates the memory requirements of each solver on the same set of MPOPF benchmarks as in Figure \ref{fig:increasingN_timeStacked}. Clearly, the black-box approaches have the most extensive memory requirements, related to storing the factors for the full KKT systems. The structure-exploiting approach implemented in \BELTISTOSMP{} reduces the memory requirements by more than one order of magnitude by usage of efficient linear algebra components adopted for the particular structure of the KKT system. Additionally, \BELTISTOSMP{} can be executed using the memory saving computation model, represented by \BELTISTOSMEM{}. The memory requirements can be further reduced by releasing the memory required to store the $L$ factors of the diagonal blocks, and recomputing the factorization during different phases of the Schur algorithm. Obviously, the redundant computations in the memory efficient algorithm are reflected in increased execution time. However, the execution time is still orders of magnitude faster than the general purpose approaches.

\section*{Acknowledgement}
This project is carried out within the frame of the Swiss Centre for Competence in Energy Research on the Future Swiss Electrical Infrastructure (SCCER-FURIES) with the financial support of the Swiss Innovation Agency (Innosuisse - SCCER program).

\bibliographystyle{plain}
\bibliography{OPFbibliography}

\end{document}